%% file: main.tex
\documentclass[12pt]{article}
 
\usepackage{myarticle}
\usepackage[colorinlistoftodos,prependcaption,textsize=tiny,backgroundcolor=black!5!white,bordercolor=red]{todonotes}

\usepackage{pgfplots}
\pgfplotsset{
  tick label style={font=\footnotesize},
  label style={font=\footnotesize},
  legend style={font=\footnotesize}
}
\input{commands}

\graphicspath{{./figures/}}
\begin{document}

\thispagestyle{empty}
\begin{center}
\vspace*{0.03\paperheight}
{\Large\bf Fixed-Point Automatic Differentiation of \\[0.8ex] Forward--Backward Splitting Algorithms for \\[0.8ex] Partly Smooth Functions \par}
% {\Large\bf Fixed-Point Automatic Differentiation for Partly \\[0.6ex] Smooth Functions:
% Improving Convergence and \\[1ex] Memory Requirements} \\
\bigskip
\bigskip
{\large
Sheheryar Mehmood$^\star$ and Peter Ochs$^\dagger$ \\ \medskip
%\today 
{\small
$^\star$~University of T\"{u}bingen, T\"{u}bingen, Germany \\
$^\dagger$~Saarland University, Saarbr\"{u}cken, Germany \\
}
}
\end{center}
\bigskip

% ********************
% >>>>> ABSTRACT <<<<<
% ********************
\begin{abstract}
A large class of non-smooth practical optimization problems can be written as minimization of a sum of smooth and partly smooth functions. We examine such structured problems which also depend on a parameter vector and study the problem of differentiating its solution mapping with respect to the parameter which has far reaching applications in sensitivity analysis and parameter learning problems. Under partial smoothness and other mild assumptions, we apply Implicit (ID) and Automatic Differentiation (AD) to the fixed-point iterations of proximal splitting algorithms. We show that AD of the sequence generated by these algorithms converges (linearly under further assumptions) to the derivative of the solution mapping. For a variant of automatic differentiation, which we call Fixed-Point Automatic Differentiation (FPAD), we remedy the memory overhead problem of the Reverse Mode AD and moreover provide faster convergence theoretically. We numerically illustrate the convergence and convergence rates of AD and FPAD on Lasso and Group Lasso problems and demonstrate the working of FPAD on prototypical image denoising problems by learning the regularization term.
\end{abstract}

%\makekeywords

%\tableofcontents

% ************************
% >>>>> INTRODUCTION <<<<<
% ************************
\section{Introduction} \label{sec:intro}
Given a parameter $\u\in\spPrm$, we consider an additive composite optimization problem of the form
\begin{equation} \label{prob:min} \tag{\mbox{$\mathcal P$}}
    \min_{\x\in\spVar} F(\x, \u) \,, \quad F\coloneqq f+g \,,
\end{equation}
where $\map{f}{\spVar\times\spPrm}{\R}$ is $C^2$-smooth and $\map{g}{\spVar\times\spPrm}{\eR}$ is possibly non-smooth and extended-valued with a simple proximal mapping when the second argument is fixed. Let $\smap\Argmap\spPrm\spVar$ be such that $\Argmap(\u)\coloneqq \Argmin_{\x\in\spVar} F(\x, \u)$ gives us the set of solutions of \eqref{prob:min} depending on $\u$. Often, we are interested in computing the variation of $\Argmap$, governed by (generalized) differentiation, which has many far-reaching applications, for example, in Machine Learning, as we discuss in Section~\ref{ssec:intro:apps}.

In general, it is quite challenging to ensure the differentiability of $\Argmap$ on some open set $\neighbourhoodPRM\subset\spPrm$. However, in special situations, suitable approaches exist. While there are analysis tools for generalized differentiation, from a practical and computational perspective the most frequently used approach is to identify points or neighbourhoods where classical derivatives are well-defined and efficiently computable. In fact, for locally Lipschitz continuous functions, differentiability is asserted almost everywhere thanks to Rademacher's Theorem which is the foundation for constructing a commonly used generalized derivative, the Clarke subdifferential \cite{Cla83}.

In this paper, we consider computational approaches for computing classical derivatives under special conditions within a generally non-smooth setup and thereby broadening the understanding by providing strong theoretical guarantees.

For illustration, let us first consider a locally smooth setup where $g=0$ in \eqref{prob:min} and $\Hess[\x] f (\x^*, \u)$ is positive definite for $\set{\xmin}\subset\Argmap (\u)$, we can employ the Implicit Function Theorem on the necessary optimality condition $\grad[\x] f (\x^*, \u) = 0$ to establish the local differentiability of the solution mapping around $\u$. An alternative to the so-called Implicit Differentiation (ID) approach is to iteratively construct a sequence of derivatives with increasing accuracy to the derivative of the solution mapping. In particular, we run an algorithm for say $K\in\N$ iterations to approach the critical point of interest $\xK\approx\xmin$ with the desired accuracy (called forward pass) and then compute the derivative $D_{\u} \xK$. Good convergence, stability or robustness properties of the optimization algorithm can often be transferred to the AD algorithm which is essentially computed through the classical chain rule for differentiation along the composition of iteration mapping in the forward pass. Nevertheless, on the practical side, AD must be used cleverly to circumvent the computational bottleneck (for instance, the memory overhead). Both convergence and memory issues have been addressed previously \cite{Gil92, Chr94, GBC+93, OBP15, MO20, Rii20, BKM+21} in standard setups.

However, when $g\neq0$ in \eqref{prob:min}, which is mostly the case in practical applications in Machine Learning when low complexity (for example, sparse or low-rank) solutions are sought through the formulation of regularization \cite{Tib96, YL06, Faz02, RFP10}, using ID in the same way as the smooth case is not possible theoretically. Similarly, using AD on forward--backward splitting algorithms, which are perfectly suited for problems of the form \eqref{prob:min}, are problematic theoretically and practically. Such methods involve a proximal mapping which in itself is an optimization problem and is not differentiable in general. For example, consider Proximal Gradient Descent (PGD) \cite{LM79} which is an iterative algorithm and has the update step
\begin{equation} \label{eq:prox:grad}
    \xkp \coloneqq \argmin_{\x\in\spVar} g(\x, \u) + \scal{\grad[\x] f(\xk, \u)}{\x - \xk} + \frac{1}{2\alpha} \norm{\x - \xk}^2 \,,
\end{equation}
where $\alpha$ denotes the step size. Computing $D_{\u}\xkp$ in terms of $D_{\u}\xk$ requires differentiating the solution $\xkp$ with respect to the parameters $(\xk, \u)$, which clearly depends on the non-smoothness properties of $g$.

Nevertheless, non-smoothness in practical problems appears in a structured manner. In particular, low-complexity regularization is governed by the low-dimensional activity of the constraint sets associated to the optimal solution $\xmin$ which, under some assumptions, are stable under minor perturbations. This activity often reveals some hidden smoothness, as it usually defines a set that is a smooth manifold $\manif$ containing $\xmin$, that is, the problem restricted to $\manif$ is smooth. Fortunately, under some assumptions, this activity is identified by PGD in a finite number of iterations \cite{LFP14, LFP17}. That is, the sequence generated by the algorithm eventually lies in $\manif$. This is exactly the starting point for the research in this paper on various theoretical and computational aspects of AD for such algorithms. The key to enjoy this property is provided by the choice of partly-smooth functions\cite{Lew02}. As an example, we consider the Lasso problem (see Equation~\ref{prob:lasso}) \cite{Tib96}. Solving this problem yields a sparse solution for large enough regularization parameter $\lambda$, a simple but illustrative example of a low-complexity regularizer. The activity in this case is the support of $\xmin$ given by $\supp {\x^*} \coloneqq \{ i : x^*_i \neq 0 \} $. It remains fixed if we change $\lambda$ slightly and is identified by forward--backward splitting algorithms \cite{LFP14, LFP17} such as PGD in a finite number of iterations. This simple idea is lifted to a broad theoretical and computational grounding with the framework of partly smooth functions. This powerful framework gives a more general Implicit Function Theorem for \eqref{prob:min} that is applicable in some non-smooth scenarios, that is, when $g$ and hence $F$ is partly smooth (see Theorem~\ref{thm:IFT}). This not only facilitates ID but also paves way for AD since the update step of PGD is also an optimization problem of the form \eqref{eq:prox:grad} with partly smooth objective (see Section~\ref{ssec:PDM:AD} for more detail).

\textbf{Contributions:} Our main contributions are the following:
\begin{enumerate}[label=(\roman*)]
    \item We apply the Implicit Function Theorem from \cite{Lew02} to the fixed-point equations of Proximal Gradient Descent (PGD) and Accelerated Proximal Gradient (APG) method (also known as FISTA \cite{BT09}) to obtain the derivative of the solution mapping of \eqref{prob:min} when $F$ belongs to a large subclass of convex partly-smooth functions which encompasses various practical problems \cite{VDP+17}. We also analyze Implicit Differentiation (ID) of these equations near the fixed-point.
    \item We provide convergence guarantees for the sequences obtained by Automatic Differentiation (AD) applied to APG method to the derivative of the solution mapping of \eqref{prob:min} for the same class of functions. We show that the convergence rate of the iterates of PGD is simply mirrored in the AD sequence.
    \item We show that Fixed-Point Automatic Differentiation (FPAD), a previously known Neumann series based method, borrows all the nice properties of both ID and AD to estimate the derivative of the solution mapping of \eqref{prob:min}.
    \item We provide experimental demonstration of convergence and convergence rates of FPAD and AD applied to PGD and APG for different practical applications, including a bilevel formulation for learning the regularizer of an image denoising problem \cite{ROF92} beyond the smooth setting.
\end{enumerate}

\subsection{Applications} \label{ssec:intro:apps}

Consider bilevel optimization \cite{DKP+15}, which can be seen as a formulation of parameter learning in Machine Learning, where the objective of one optimization problem (the upper or outer problem) depends on the solution of another problem (the lower or inner problem). The gradient based approaches require computing the derivative of the solution mapping of the inner problem. We will discuss two of the major applications here. For other applications like Meta-Learning and the so-called Learning-to-Optimize or L2O, we request the reader to go through some recent surveys, for example, \cite{CCC+21, HAMA20}.

\textbf{Hyperparameter Learning:}
Hyperparameters in neural network or any machine learning model
can be learned through grid search \cite{BB12} or Bayesian Optimization \cite{Moc75, BCF10}. More recently the trend has shifted to gradient-based approach to compensate for the growing number of hyperparameters \cite{Ben00, Dom12, KP13, DVFP14, MDA15a, Ped16}. This approach involves differentiating the validation loss (upper-level problem) with respect to the hyperparamters by differentiating the solution mapping of the training loss (lower-level problem) through Implicit or Automatic Differentiation and has been put to test to solve many practical problems like image denoising \cite{Dom12, KP13}, image segmentation \cite{ORBP16} and data cleaning \cite{FDF+17}.

\textbf{Implicit Models:} In recent years, there has been a tendency towards embedding a function defined by an implicit equation inside the prediction function of the machine learning model. Training such models is equivalent to solving a bilevel problem with the implicit equation being the lower problem and the training of the model being the upper problem. Examples of implicit models include Deep Equilibrium Models \cite{BKK19, BKK20, WK20}, Neural ODE's \cite{CRBD18} and Optimization Layers \cite{AK17, MB18, AAB+19, BBT+20}.
When the implicit equation in such a model represents some inherent structure of the practical application, the learned model may provide better interpretability and generalization.

\subsection{Related Work} \label{ssec:intro:RW}

The problem of computing the derivative of the solution mapping of a structured problem in \eqref{prob:min} has been studied extensively in a bilevel setup. We discuss here some previous works whose methodologies apply to certain classes of non-smooth problems.

The non-smoothness of \eqref{prob:min} was circumvented in \cite{DVFP14, ORBP15, ORBP16} by applying AD on the algorithms which are differentiable (almost) everywhere without providing the convergence guarantees for the derivative sequence. Deledalle et al.~\cite{DVFP14} utilize weak-differentiablity of PGD in \eqref{eq:prox:grad} thanks to the Rademacher's theorem. Ochs et al.~\cite{ORBP15, ORBP16}, on the other hand, use a smooth algorithm which was obtained by carefully selecting the kernel generating distance function for Bregman Proximal Gradient \cite{Bre67, BBT17}. Bolte et al.~\cite{BLPS21} used the conservative calculus developed in \cite{BP20, BP21} to rigorously study the implicit differentiation performed by using standard autograd packages like TensorFlow \cite{ABC+16}, PyTorch \cite{PGM+19}, JAX \cite{BFH+18}, etc. on possibly non-smooth implicit equations. Thereby, they provide a mathematical framework for the types of ``derivatives'' that are used in many practical applications. The same tools in \cite{BP20, BP21} were applied to study the application of autograd packages to differentiate the output of a non-smooth iterative algorithm in \cite{BPV22}.

Christof \cite{Chr20} replaced \eqref{prob:min} with an elliptic variational inequality \cite[\DRSecs{2.1}{2A}]{DR09}. They showed that the solution mapping for this type of inequality is Fr{\'e}chet differentiable and the derivative can be obtained by solving an analogous elliptic variational equality. Chambolle and Pock~\cite{CP21} learn a linear operator in a saddle point problem with convex and $C^2$-smooth objective. They show that the gradient of the upper-level objective with respect to the linear operator is obtained by solving an adjoint-state equation which is equivalent to solving another saddle point problem. The authors showed that the two problems can be solved together in a Piggyback style \cite{GF03}. More recently, the convergence guarantees for this technique are given for $C^1$-smooth functions \cite{BCP22}.

In \cite{BKM+21, Rii20}, the restricted smoothness of various practical non-smooth regularizers was leveraged to make PGD smooth on a neighbourhood in $\spPrm$ and to show convergence of AD applied to PGD. In \cite{BKM+21, Ber21}, only additive separable non-smooth functions were considered which do not include many interesting practical regularizers like $ \ell_{\infty}$ norm, $\ell_{2,1}$ group norm and nuclear norm \cite{VDP+17}. In \cite{Rii20}, the class of non-smooth functions was broadened to partly-smooth functions. However, the memory-overhead of Reverse Mode AD were not addressed and a restrictive assumption of strongly convex objective in \eqref{prob:min} was considered which fails on many practical problems like the Foreground/Background separation problem \cite{BSJ+16}.

\subsection{Notation \& Arrangement of the Paper} \label{ssec:intro:Not}

In the rest of the paper, we use $\N_0 \coloneqq \N\cup \set{0}$ and define $[K]\coloneqq\set{0,\ldots,K}$ for any $K\in\N$. We denote the set of extended real numbers by $\eR \coloneqq \R\cup\{-\infty, +\infty\} $ and assume that $\spVar$ and $\spPrm$ are Euclidean spaces of dimensions $N$ and $P$, respectively, each equipped with an inner product $\scal\cdot\cdot$ and an induced norm $\norm\cdot$. The corresponding dual spaces are written as $\spVar^*$ and $\spPrm^*$. $\neighbourhood (\x)$ denotes the set of all the neighbourhoods of a point $\x$. The set of all linear maps from $\spVar$ to $\spPrm$ is denoted by $\spLin (\spVar, \spPrm)$. The elements of $\eR$, $ \spVar$ $($and $\spPrm)$ and $\spLin(\spVar, \spPrm)$ are denoted by lower-case normal (for example, $t$), lower-case bold (for example, $\x$) and upper-case normal (for example, $A$) letters, respectively. For $\x\in\R^N$ we write $\x \coloneqq (x_1 , \ldots, x_N ) = (x_i)_{i=1,\ldots,N}$ and for $A\in\R^{N\times P}$ we write $A =(a_{i,j})_{i=1,\ldots,N,j=1,\ldots,P}$. For a linear operator defined on a subspace $V$ of $\spVar$, say $A\in\spLin (V, \spVar)$, we may extend its domain to $\spVar$ by composing it with the orthogonal projection onto V, that is, $A \Pi_{V}$ and abuse the notation by denoting the resulting operator, also, by $A$.

For the subdifferential of a proper and lower semi-continuous function, we follow the terminology of \cite{RW98}. However, our setting simplifies to using the standard convex (Fr{\'e}chet) subdifferential because we deal with regular functions \cite[Definition~7.25]{RW98} of the form $\map{g}{\spVar\times\spPrm}{\eR}$ which are convex in the first argument and smooth in the second. For such $g$, thanks to \cite[Corollary~10.11]{RW98}, $\partial g (\x, \u)$ simplifies to the Cartesian product of the convex subdifferential of $ g(\cdot, \u)$ at $\x$ and a set containing only the gradient of $g(\x, \cdot)$ at $\u$.
\begin{lemma} \label{lem:subdiff:cartesian}
    Let $\map{g}{\spVar\times\spPrm}{\eR}$ be proper and lower semi-continuous, and $(\x, \u)\in\spVar\times\spPrm$ be a point. If $g$ is regular at $(\x, \u)$ and $g(\x, \cdot)$ is differentiable at $\u$, then we have:
    \begin{equation} \label{eq:subdiff:cartesian}
        \partial g (\x, \u) = \partial_{\x} g (\x, \u) \times \set{\grad[\u] g(\x, \u)} \,.
    \end{equation}
\end{lemma}
\findProofIn{lem:subdiff:cartesian}
\textbf{Arrangement of the Paper:} The rest of the paper is arranged as follows. We recall some classical results regarding parametric fixed-point iterations and differentiating the fixed-point with respect to the parameter by using ID, AD and FPAD in Section~\ref{ssec:intro:FPI}. Section~\ref{sec:Pre} provides a quick review of partial smoothness and some basic terminologies from Riemannian geometry. We introduce the class of functions, on which, our analysis will be applicable, in Section~\ref{sec:PS} and present Algorithm~\ref{alg:APG} --- the classical algorithm to solve \eqref{prob:min} for such functions --- in Section~\ref{sec:APG}. Following the outline of the classical results in Section~\ref{ssec:intro:FPI}, we analyze the differentiation of the solution mapping of \eqref{prob:min} in Section~\ref{sec:PDM} using ID, AD and FPAD applied to Algorithm~\ref{alg:APG} from Section~\ref{sec:APG} for the problem class from Section~\ref{sec:PS}. Finally, in Section~\ref{sec:Exp}, we demonstrate our theory on two sets of experiments; we verify the convergence rates on the lasso and the group lasso problems (see Section~\ref{ssec:Exp:CRD}) and use FPAD to differentiate the lower level problem in a bilevel optimization setup for learning the TV regularization operator (see Section~\ref{ssec:Exp:ID}).
\subsection{Approaches for Differentiating Fixed-Point Mappings} \label{ssec:intro:FPI}
Before we move on to the main part of our paper, we recall the classical results on differentiating an implicit map obtained by solving a non-linear parametric equation. In particular, given $\map{\fixMap}{\spVar\times\spPrm}{\spVar}$ and $\u\in\spPrm$, we consider the problem of solving
\begin{equation} \label{eq:FPE} \tag{FPE}
    \x = \fixMap (\x, \u) \,,
\end{equation}
with respect to $\x$. The corresponding solution is a fixed-point of $\fixMap(\cdot, \u)$ and under certain assumptions, is guaranteed to exist and is obtainable by the following so-called fixed-point iterations
\begin{equation} \label{itr:FPI} \tag{FPI}
    \xkp (\u) \coloneqq \fixMap(\xk (\u), \u) \,,
\end{equation}
for $k\geq0$ where $\xz\in\spVar$ is the initial point. Notice that both the sequence $\xk (\u)$ and the fixed-point depend on $\u$ due to the dependence of $\fixMap$ on $\u$. We are interested in computing the derivative of the fixed-point as a function depending on $\u$ through ID and AD. In this section, we recall the assumptions under which, the iterates from \eqref{itr:FPI} converge to the fixed-point \cite[Section~2.1.2, Theorem~1]{Pol87}, and two differentiation techniques that provide a good estimate of the derivative \cite{Gil92, Chr94}. Furthermore, we also present FPAD, a memory-efficient variant of AD for fixed-point iterations originally proposed in \cite{Chr94}. We will apply these results for differentiating the solution of \eqref{prob:min}  (with respect to $\u$) through the fixed-point iterations of APG.

Beck~\cite{Bec94} provided a general treatment on the implicit and automatic differentiation of the fixed-point when $\fixMap$ in \eqref{itr:FPI} is replaced by an iteration dependent mapping $\map{\fixMap[k]}{\spVar\times\spPrm}{\spVar}$ where $\fixMap[k]$ satisfies some boundedness assumptions and $\fixMap[k]\to\fixMap$ pointwise. Such a setting encapsulates the cases when a changing step size or inertial parameter is employed in an algorithm, for instance, the Nesterov's accelerated gradient method \cite{Nes83} and APG \cite{BT09} which \eqref{itr:FPI} does not. However, the convergence rate guarantees for AD requires more assumptions on $\fixMap[k]$ which may not be practical, for example, linear convergence of $\fixMap[k]\to\fixMap$.

The roots of correctness and effectiveness of all the differentiation techniques can be traced back to the following underlying condition on $\fixMap$ which readily generalizes to a whole neighbourhood as shown in Lemma~\ref{lem:FPE:neighbourhood} below.
\begin{assumption} \label{ass:C1:contraction}
    $\fixMap$ is $C^1$-smooth on a neighbourhood of $(\xmin, \givenPrm)$ with $\xmin = \fixMap (\xmin, \givenPrm)$ and $\rho (D_{\x} \fixMap (\xmin, \givenPrm)) < 1$.
\end{assumption}
\begin{lemma} \label{lem:FPE:neighbourhood}
    Let $(\xmin, \givenPrm)$ be such that Assumption~\ref{ass:C1:contraction}, is satisfied for $\fixMap$. Then $\exists \ \neighbourhoodFM\in\neighbourhood (\xmin, \givenPrm)$ such that $\fixMap$ is $C^1$-smooth on $\neighbourhoodFM$ and $\rho (D_{\x} \fixMap (\x, \u)) < 1$ for all $(\x, \u) \in \neighbourhoodFM$.
\end{lemma}

\subsubsection{Implicit Differentiation} \label{sssec:FPI:ID}
The first approach for computing the derivative of the fixed-point map is Implicit Differentiation. In order to formulate the main theorem, we invoke the following lemma.
\begin{lemma} \label{lem:spec_rad:inv_op}
    For any linear operator $\map{B}{\spVar}{\spVar}$ with $\rho (B) < 1$, the linear operator $\opid - B$ is invertible.
\end{lemma}
\findProofIn{lem:spec_rad:inv_op}
From Assumption~\ref{ass:C1:contraction} and Lemma~\ref{lem:spec_rad:inv_op}, we obtain the following well known result \cite[Theorem~1B.1]{DR09}.
\begin{theorem}[Implicit Function Theorem] \label{thm:IFT:FPE}
    Let $(\xmin, \givenPrm)$ be such that Assumption~\ref{ass:C1:contraction}, is satisfied for $\fixMap$. Then $\exists \ \neighbourhoodPRM\in\neighbourhood (\givenPrm)$ and a $C^1$-smooth mapping $\map{\argmap}{\neighbourhoodPRM}{\spVar}$ such that $\forall \u\in \neighbourhoodPRM$, $\argmap(\u) = \fixMap(\argmap(\u), \u)$, and
    \begin{equation} \label{eq:FPE:IFT}
        D\argmap(\u) = \big( \opid - D_{\x} \fixMap (\argmap(\u), \u) \big)^{-1} D_{\u} \fixMap (\argmap(\u), \u) \,.
    \end{equation}
\end{theorem}
Using Lemmas~\ref{lem:FPE:neighbourhood} and \ref{lem:spec_rad:inv_op}, we define a mapping $\map{\phi}{\neighbourhoodFM}{\spLin (\spVar, \spPrm)}$ on $\neighbourhoodFM$ from Lemma~\ref{lem:FPE:neighbourhood} by
\begin{equation} \label{eq:FPE:IFT:approx}
    \phi (\x, \u) \coloneqq \big( \opid - D_{\x} \fixMap (\x, \u) \big)^{-1} D_{\u} \fixMap (\x, \u) \,.
\end{equation}
In practice, for a given $\u\in \neighbourhoodPRM$, we run $\eqref{itr:FPI}$ for a sufficiently large number of iterations to obtain an estimated fixed-point $\x\approx\argmap (\u)$. If the estimate $\x$, is close enough to $\argmap (\u)$, then $(\x, \u)\in\neighbourhoodFM$ and we may compute $\phi (\x, \u)$ as an estimate of $D \argmap (\u) = \phi (\argmap (\u), \u)$. Finally, mostly in practice, once we obtain the estimate $\x$ of the fixed-point, we use it to evaluate a $C^1$-smooth scalar-valued objective say $\map{\ell}{\spVar}{\R}$. Our goal then is to differentiate $\ell (\x)$ with respect to $\u$. When $\xminb \coloneqq D \ell (\x)$, $ \xminb\phi (\x, \u)$ provides an estimate of the derivative $D_{\u}\ell (\x)$ from the chain rule. However, in this situation, the product $\xminb (\opid - D_{\x}\fixMap (\x, \u))^{-1}$ is computed first, possibly through a linear solver, and the result is then multiplied by $D_{\u}\fixMap (\x, \u)$. 
Implicit differentiation is a useful practical tool for differentiating the solution mapping of a fixed-point equation because it does not have a memory overhead and it provides flexibility as compared to automatic differentiation. However, this flexibility comes at a price of custom implementation which increases the workload in practice.

\subsubsection{Automatic Differentiation} \label{sssec:FPI:AD}

AD is a fundamental tool today for automating the derivative computation of various practical functions. A brief account on the two modes of AD is given in \refSupp{sec:AD:supp}. For further reading, the reader is requested to look into \cite{GW08,BPR+18}. AD is another contender in practice for estimating the derivative of the solution mapping \cite{Dom12} and was first studied by \cite{Gil92} with convergence guarantees. The finite sequence $\xk_{k\in[K]}$ from \eqref{itr:FPI} is unrolled (see Figure~\ref{fig:FPI}) and then we simply use the chain rule to compute the derivative of $\xK$ with respect to $\u$.
\begin{figure}[ht]
	\centering
	\tikzset{font={\fontsize{8pt}{12}\selectfont}}
	\begin{tikzpicture}[->, semithick, node distance=4cm, >=stealth', auto]
	\tikzstyle{place} = [circle, thick, draw=black!75, fill=black!10, minimum size=6mm]
	\begin{scope}
	
	\node [place] (xz) [label=above:$\xz$] {};
	\node [place] (x1) [right of=xz, label=above:$\x^{(1)}$] {};
	\node [place] (xKm) [right of=x1, label=above:$\x^{(K-1)}$] {};
	\node [place] (xK) [right of=xKm, label=above:$\xK$] {};
	\node [place] (u) [below of=xz, label=above left:$\u$] {};
	
	\path[every node/.style={sloped,anchor=south,auto=false}]
	(xz) edge node [above, pos=0.5] {$D_{\x} \fixMap (\xz, \u)$} (x1)
	(xKm) edge node [above, pos=0.5] {$D_{\x} \fixMap (\x^{(K-1)}, \u)$} (xK)
	(x1) edge [dotted] node {} (xKm)
        (u) edge [dotted] node [above, pos=0.6] {$D_{\u} \xz$} (xz)
	(u) edge [bend right] node [above, pos=0.6] {$D_{\u} \fixMap (\xz, \u)$} (x1)
	(u) edge [bend right] node [above, pos=0.6] {$D_{\u} \fixMap (\x^{(K-2)}, \u)$} (xKm)
	(u) edge [bend right] node [above, pos=0.6] {$D_{\u} \fixMap (\x^{(K-1)}, \u)$} (xK);
	
	\end{scope}
	\end{tikzpicture}
	\caption{Depiction of Unrolling of \eqref{itr:FPI}. $\xz$ need not be independent of $\u$.}
	\label{fig:FPI}
\end{figure}
For the forward mode, given $\xzd \coloneqq 0\in\spVar$ and $\ud\in\spPrm$, we perform the following iterations for $k = 0, \ldots, K-1$
\begin{equation} \label{itr:FPI:FAD} \tag{FPI-F}
    \xkpd \coloneqq D_{\x} \fixMap (\xk, \u) \xkd + D_{\u} \fixMap (\xk, \u) \ud \,,
\end{equation}
to obtain the forward derivative $\xKd$. We similarly apply reverse mode AD to \eqref{itr:FPI} by starting with $\xKb\coloneqq \xminb\in\spVar^*$, and $\uzb\coloneqq 0\in\spPrm^*$ and computing for $n = 0,\ldots,K-1$,
\begin{equation} \label{itr:FPI:RAD} \tag{FPI-R}
    \begin{aligned}
        \xkb &\coloneqq \xkpb D_{\x} \fixMap (\xk, \u) \\
        \unpb &\coloneqq \unb + \xkpb D_{\u} \fixMap (\xk, \u) \,,
    \end{aligned}
\end{equation}
where $k \coloneqq K - n - 1$. Notice that in \eqref{itr:FPI:RAD}, we introduce a new index $n$ which increases as we move from right (output) to left (input) along the graph in Figure~\ref{fig:FPI} to distinguish it from $k$. The subscript $K$ in the reverse mode derivatives here is due to the fact that we are computing the derivative of $\xK$. The output of reverse mode AD applied to \eqref{itr:FPI} is $\ubK \coloneqq \uKb$.
\begin{remark}
    \begin{enumerate}[label=(\roman*)]
        \item We omit the dependence of $\xk$, $\xkd$ and $\unb$ on $\u$ in \eqref{itr:FPI:FAD} and \eqref{itr:FPI:RAD} for brevity.
        \item Even though the reverse mode has a memory overhead unlike the forward mode, it is computationally efficient for most practical applications and is preferred in practice.
        \item We may assume that $\xz$ is independent of $\u$ which implies that the starting iterates are $\xzd=0$ and $\uzb=0$. However, the dependence of $\xz$ on $\u$ does not have any effect on the convergence of the derivative sequence as established by Theorem~\ref{thm:LGFPI:conv}.
    \end{enumerate}
\end{remark}
\textbf{Linear Convergence of Iterative Process:} Since the convergence of derivatives is closely linked to that of the iterates, we first recall the following result from \cite[Section~2.1.2, Theorem~1]{Pol87}.
\begin{theorem}[Local Linear Convergence of Iterates] \label{thm:FPI}
    Let $(\xmin, \givenPrm)$ be such that Assumption~\ref{ass:C1:contraction} is satisfied for $\fixMap$. Then the sequence $\xk (\givenPrm)$ generated by \eqref{itr:FPI} converges locally linearly to $\xmin$ and for all $\delta \in (0, 1 - \rho)$ there exist $\epsilon > 0$ and $C(\delta)$ such that for all $k\in\N$
    \begin{equation*}
        \norm{\xk (\givenPrm) - \xmin} \leq C(\delta) (\rho + \delta)^k \,,
    \end{equation*}
    for $\norm{\xz - \xmin} \leq \epsilon$ where $\rho\coloneqq \rho (D_{\x} \fixMap (\xmin, \givenPrm))$.
\end{theorem}
\begin{remark} \label{rem:FPI}
    \begin{enumerate}[label=(\roman*)]
        \item Alternatively, instead of assuming that $\xz$ lies near $\xmin$, we may assume that the sequence $\xk (\givenPrm)\to\xmin$ and then prove its asymptotic linear convergence.
        \item \label{itm:FPI:conv:neigh} In Theorem~\ref{thm:FPI}, the linear convergence is implied for any $\u$ in the neighbourhood $\neighbourhoodPRM$ of $\givenPrm$ defined in Theorem~\ref{thm:IFT:FPE} provided that the sequence $\xk (\u)\to\argmap (\u)$. The rate of convergence is then $\rho (D_{\x} \fixMap (\argmap (\u), \u))$.
    \end{enumerate}
\end{remark}

The following result is derived from \cite[Section~2.1.2, Theorem~1]{Pol87} and is useful for proving the (linear) convergence of the derivative iterates.
\begin{theorem} \label{thm:LGFPI:conv}
Let $( B_k )_{k\in\N_0}$ and $ ( \bk)_{k\in\N_0}$ be sequences in $\spLin(\spVar, \spVar)$ and $\spVar$ with limits $B$ and $\b$, respectively. If $\rho=\rho (B) < 1$, the sequence $( \xk )_{k\in\N}$, with $\xz\in \spVar$, generated by
\begin{equation} \label{itr:LGFPI}
    \xkp \coloneqq B_k \xk + \bk \,,
\end{equation}
converges to $\x^*\coloneqq(\opid - B)^{-1} \b$.
% The convergence is linear when $( B_k)_{k\in\N_0}$ and $(\bk)_{k\in\N_0}$ converge linearly with rates $q_B$ and $q_{\b}$, respectively. In fact, for all $\delta \in (0, 1 - \rho)$, there exists $C(\delta)$, such that for all $k\in\N$, we have
% \begin{equation} \label{eq:LGFPI:conv_rates}
%     \norm{\xk - \x^*} \leq \begin{cases} C(\delta) k q^k,& \text{if } \rho+\delta = \max(q_B, q_{\b})\,; \\ C(\delta) q^k,& \text{if } \rho+\delta \neq \max(q_B, q_{\b})\,, \end{cases}
% \end{equation}
% where $q\coloneqq \max (\rho + \delta, q_B, q_{\b})$.
If the convergence of $( B_k)_{k\in\N_0}$ and $(\bk)_{k\in\N_0}$ is linear and for all $\delta \in (0, 1 - \rho)$ and for some $c(\delta)$ and $K\in\N$,
\begin{equation} \label{ineq:LGFPI:B:b:conv_rates}
    \begin{aligned}
        \norm{B_k - B} &\leq c(\delta) (\rho + \delta)^{k-K} \\
        \norm{\bk - \b} &\leq c(\delta) (\rho + \delta)^{k-K} \,,
    \end{aligned}
\end{equation}
for all $k\geq K$, then we have
\begin{equation} \label{ineq:LGFPI:conv_rates}
    \norm{\xk - \x^*} \leq c_1(\delta) (k-K) (\rho + \delta)^{k-K} + c_2(\delta) (\rho + \delta)^{k-K} \,,
\end{equation}
for all $k\geq K$ and for some $c_1(\delta)$ and $c_2(\delta)$.
\end{theorem}
\findProofIn{thm:LGFPI:conv}
\begin{remark}
    The convergence rate in Theorem~\ref{thm:LGFPI:conv} is non-asymptotic when $K=0$.
\end{remark}
\textbf{Local (Linear) Convergence of AD Process:} For a given $\u$, when $\xk (\u)$ is close enough to the fixed-point, the update step in \eqref{itr:FPI:FAD} furnishes well-defined iterates $\xkd$ which converge to the derivative of the fixed-point. The convergence is asymptotically linear under an additional local Lipschitz continuity assumption.
\begin{corollary}[Convergence of Derivative Iterates] \label{cor:FPI:FAD}
    Let $(\xmin, \givenPrm)$ be such that Assumption~\ref{ass:C1:contraction} is satisfied for $\fixMap$. If the sequence $\xk (\givenPrm)$ generated by \eqref{itr:FPI} is close enough to $\xmin$, then, for any $\ud\in\spVar$, the sequence $\xkd (\givenPrm)$ generated by \eqref{itr:FPI:FAD} converges to $D\phi (\givenPrm)\ud$. When $D\fixMap$ is additionally locally Lipschitz continuous near $(\xmin, \givenPrm)$, the convergence is linear and for all $\delta > 0$, there exist $C(\delta)$ and $K\in\N$ such that for all $k\geq K$, we have
\begin{equation} \label{eq:FPI:FAD:conv_rates}
    \norm{\xkd (\givenPrm) - \xmind} \leq C(\delta) (k-K) (\rho + \delta)^{k-K} \,,
\end{equation}
where $\rho \coloneqq \rho (D_{\x} \fixMap (\xmin, \givenPrm))$.
\end{corollary}
\findProofIn{cor:FPI:FAD}
\begin{remark}
    \begin{enumerate}[label=(\roman*)]
        \item Our discussion in Remark~\ref{rem:FPI}\ref{itm:FPI:conv:neigh} similarly extends to the convergence of $\xkd (\u)$ to $D\phi (\u)\ud$ for any $\u\in \neighbourhoodPRM$ defined in Theorem~\ref{thm:IFT:FPE}.
        \item The additional term $k$ in the rate of convergence of $\xkd$ as compared to that of $\xk$ contributes to the curse of unrolling \cite{SGB+22}.
        \item The convergence guarantees for reverse mode are the same as those of the forward because the two modes essentially compute the same quantity in the end.
        \item The convergence rate is non-asymptotic when $D\fixMap$ has global Lipschitz continuity.
    \end{enumerate}
\end{remark}
\textbf{Global Convergence of the AD Process:} If $\fixMap$ is $C^1$-smooth on the whole domain, we do not require all the iterates $\xk$ to be close enough to the fixed-point. The derivative iterates converge regardless of the initial point $\xz$ as established by the following result.
\begin{corollary}[Global Convergence of Derivative Iterates] \label{cor:FPI:FAD:Global}
    Let $(\xmin, \givenPrm)$ be such that Assumption~\ref{ass:C1:contraction} is satisfied for $\fixMap$. If $\fixMap$ is $C^1$-smooth on the whole space and $\xk (\givenPrm)\to\xmin$, then the conclusion of Corollary~\ref{cor:FPI:FAD} holds for any $\xz\in\spVar$.
\end{corollary}
\findProofIn{cor:FPI:FAD:Global}
\textbf{Late-Starting the AD Process:} When $\fixMap$ is $C^1$-smooth only near the fixed-point and the initial iterates $\xk$ are too far away, the iterations in \eqref{itr:FPI:FAD} may not be well-defined for small $k$. In this situation, we may rely on late-starting, that is, we start the AD process after performing the updates in \eqref{itr:FPI} for sufficiently large number of iterations, say $N\in\N$. More formally, for $k = 0, \ldots, K-1$, the forward mode update is given by
\begin{equation} \label{itr:FPI:FAD:Late}
    \xkpd \coloneqq D_{\x} \fixMap (\x\iter{k+N}, \u) \xkd + D_{\u} \fixMap (\x\iter{k+N}, \u) \ud \,.
\end{equation}
For reverse mode, when performing the backward pass, we stop at $\x\iter N$ instead of $\xz$, that is,
\begin{equation} \label{itr:FPI:RAD:Late}
    \begin{aligned}
        \xkb &\coloneqq \xkpb D_{\x} \fixMap (\x\iter{k+N}, \u) \\
        \unpb &\coloneqq \unb + \xkpb D_{\u} \fixMap (\x\iter{k+N}, \u) \,,
    \end{aligned}
\end{equation}
for $n = 0,\ldots,K-1$, where $k \coloneqq K - n - 1$.
\begin{corollary}[Convergence with Late-Start] \label{cor:FPI:FAD:Late}
    Let $(\xmin, \givenPrm)$ be such that Assumption~\ref{ass:C1:contraction} is satisfied for $\fixMap$. If $N$ is large enough and $\xk (\givenPrm)\to\xmin$, then the conclusion of Corollary~\ref{cor:FPI:FAD} also holds when $\xkd (\givenPrm)$ is generated by $\eqref{itr:FPI:FAD:Late}$.
\end{corollary}
\findProofIn{cor:FPI:FAD:Late}
The reverse mode AD is more common in practice. However, the huge drawback of memory overhead can make it impractical for various applications involving differentiation of the solution mapping. Nevertheless, the ease of implementation and the translation of convergence rate still make it a strong competitor of ID.

\subsubsection{Fixed-Point Automatic Differentiation} \label{sssec:FPI:FPAD}
The ease of using AD on Algorithm~\ref{alg:APG} combined with the results shown in Section~\ref{sssec:FPI:AD} provide a powerful tool for estimating the derivative of the fixed-point mapping. However, the reverse mode AD has a memory overhead and therefore can become impractical when $\xk$ changes slowly in the beginning and is required to run for a larger $K$. Several techniques have been proposed to get around this problem, for instance, Checkpointing \cite{VO85, DH06} and Truncating \cite{WP90} to reduce the memory usage. However, this is done at the cost of either reduced accuracy of the derivative as in Truncating or increased computation time as in Checkpointing. Also, with these techniques, the memory still scales with $K$, albeit not linearly. Fortunately, in our context, a better alternative exists \cite{Gil92, Chr94, Azm97, Bar98, SWGH08} which easily lifts the curse of memory without compromising the accuracy of the derivative.

For a given $\u$, we replace the sequence $\xk (\u)$ in \eqref{itr:FPI:FAD} and \eqref{itr:FPI:RAD} with the (approximate) fixed-point $\x\coloneqq\xK$ to obtain the update steps for a modified AD approach which we call Fixed-Point Automatic Differentiation or FPAD. Here we assume that $(\x, \u) \in\neighbourhoodFM$ where $\neighbourhoodFM$ is the neighbourhood from Lemma~\ref{lem:FPE:neighbourhood}. For forward mode FPAD, the update step reads
\begin{equation} \label{itr:FPI:FPFAD} \tag{\mbox{FPI-\^F}}
    \xkph \coloneqq D_{\x} \fixMap (\x, \u) \xkh + D_{\u} \fixMap (\x, \u) \ud \,,
\end{equation}
for $k = 0, \ldots, K-1$, $\xzh\in\spVar$, and $\ud\in\spPrm$. Similarly for reverse mode, we set $\xKt\coloneqq \xminb\in\spVar^*$, and $\uzt\coloneqq 0\in\spPrm^*$ and get the following modified iterations of \eqref{itr:FPI:RAD} for $n = 0,\ldots, K-1$
\begin{equation} \label{itr:FPI:FPRAD} \tag{\mbox{FPI-\~R}}
    \begin{aligned}
        \xkt &\coloneqq \xkpt D_{\x} \fixMap (\x, \u) \\
        \unpt &\coloneqq \unt + \xkpt D_{\u} \fixMap (\x, \u) \,,
    \end{aligned}
\end{equation}
where $k \coloneqq K - n - 1$. Because the update steps are modified, we also change the notation for forward and reverse mode FPAD, that is, $\xkh$ instead of $\xkd$ for the forward mode and $\unt$ instead of $\unb$ for the reverse mode.
\begin{remark}
    The iterates $\xkh$ and $\unb$ depend on both the approximate fixed-point $\x$ and the parameter $\u$.
\end{remark}
\textbf{Resemblance of FPAD and AD:} The similarity between the AD update steps, that is, \eqref{itr:FPI:FAD} and \eqref{itr:FPI:RAD}, and the FPAD update steps, that is, \eqref{itr:FPI:FPFAD} and \eqref{itr:FPI:FPRAD} is already visible. Moreover, if we perform AD through late start, that is, \eqref{itr:FPI:FAD:Late} with $N$ sufficiently large and $(\x\iter N, \u) \in\neighbourhoodFM$, we get a computational graph with nodes $\x\iter{N} \approx \x\iter{N+1} \approx \ldots \approx \x\iter{N+K-1} \approx \x\iter{N+K}$. Performing forward and reverse mode AD on such a graph is almost like performing FPAD with $\x \coloneqq \x\iter{N}$ because then $D\fixMap (\x, \u) \coloneqq D\fixMap (\x\iter{N}, \u) \approx \ldots \approx D\fixMap (\x\iter{N+K}, \u)$.

We prove the following result which is crucial in understanding the convergence and convergence rate of FPAD and its connection with ID.
\begin{lemma} \label{lem:spec_rad:neu_ser}
    For any linear operators $B\in\spLin (\spVar, \spVar)$ and $C\in\spLin (\spPrm, \spVar)$ with $\rho\coloneqq\rho (B) < 1$, the Neumann series $\sum_{n=0}^{\infty}B^nC$ converges to $X_* \coloneqq (\opid - B)^{-1} C$ and for any $r\in\N_0$, the partial sums $X_r \coloneqq \sum_{i=0}^{r-1} B^i C$ can be expressed recursively in two different ways, that is,
    \begin{equation} \label{itr:Neu:fwd}
        X_{r+1} \coloneqq B X_r + C \,,
    \end{equation}
    and
    \begin{equation} \label{itr:Neu:bwd}
        \begin{aligned}
            Y_{r+1} &\coloneqq Y_r B \\
            X_{r+1} &\coloneqq X_r + Y_r C \,,
        \end{aligned}
    \end{equation}
    with $X_0\in\spLin (\spPrm, \spVar)$ and $Y_0\coloneqq \opid\in\spLin (\spVar, \spVar)$. Moreover for any $\delta\in(0, 1-\rho)$, there exists $C(\delta)$ such that for all $k\in\N$, we have
    \begin{equation}
        \norm{X_{r} - X_*} \leq C(\delta) (\rho + \delta)^r \,.
    \end{equation}
\end{lemma}
\findProofIn{lem:spec_rad:neu_ser}
\textbf{Resemblance of FPAD and ID:} In Lemma~\ref{lem:spec_rad:neu_ser}, if we set $B\coloneqq D_{\x} \fixMap (\x, \u)$ and $C\coloneqq D_{\u} \fixMap (\x, \u)$, we easily see the connection between FPAD and ID by comparing \eqref{itr:Neu:fwd} with \eqref{itr:FPI:FPFAD} and \eqref{itr:Neu:bwd} with \eqref{itr:FPI:FPRAD}.
\begin{remark} \label{rem:FPRAD}
Lemma~\ref{lem:spec_rad:neu_ser} suggests that \eqref{itr:FPI:FPFAD} and \eqref{itr:FPI:FPRAD} can be performed indefinitely and the total number of iterations can be decided independently of the forward pass. Also, due to $\rho (D_{\x} \fixMap (\x, \u)) < 1$, $\tilde\x^{(K-n-1)} \to 0$ as $n \to \infty$, we get a natural stopping criterion $\norm{\tilde\x^{(K-n-1)}} < \epsilon$ for reverse mode FPAD.
\end{remark}
The connection of FPAD with AD and ID already exists in literature and although not explicitly, was first made in \cite[see Section~3]{Chr94}. Different names of FPAD in literature include Post-Convergence Automatic Differentiation \cite{Azm97}, Post Differentiation \cite{Bar98} and Neumann Series method (or some variation of it). FPAD combines the best of both worlds. Unlike AD, it does not have a memory overhead and has a linear convergence rate regardless of local Lipschitz continuity (see Corollary~\ref{cor:FPI:FAD}). Unlike ID, it is very easy to implement using a standard autograd package. It is due to this ease of implementation many implementations of DEQ \cite{BKK19} just compute the partial sum of the Neumann Series rather than using a linear solver. In the appendix, we provide a PyTorch code of a simple example for computing the derivative using FPAD in Section~\ref{ssec:FPAD:Imp:supp} and a comparison of time and memory complexity for different methods of derivative computation in Section~\ref{ssec:T&MC:supp}.

\section{Preliminaries} \label{sec:Pre}
In the following subsections, we provide a gentle introduction to the fundamentals of our contribution to make the results widely accessible.

\subsection{Riemannian Geometry} \label{ssec:Pre:RG}
We recap here a few definitions and results from Riemannian Geometry and refer the reader to \cite{Lee03, Cha06} for further study.
\begin{definition}[Manifold, Tangent and Normal Spaces]
    For any $k\in\N$, we say that $\manif\subset\spVar$ is a $C^k$-smooth $m$-dimensional submanifold of $\spVar$, if for every $\x\in\manif$ there exist an open set $X\subset\spVar$ and a $C^k$-smooth map $\map{\Phi}{X}{\R^{N-m}}$ such that $\x\in X$, the derivative $D \Phi (\x)$ is surjective, and $X\cap\manif = \Phi^{-1} (0) = \set{\y \in X \setsep \Phi (y) = 0}$. We call $\Phi$, the local defining function of $\manif$ at $\x$. We say that $\v\in\spVar$ is a tangent vector of $\manif$ at $\x$ if there exist $\epsilon>0$ and a $C^1$-smooth curve $\map{\gamma}{(-\epsilon, \epsilon)}{\manif}$ on $\manif$ with $\gamma (0) = \x$ and $\dot\gamma (0) = \v$. The set of all tangent vectors of $\manif$ at $\x$ constitute $\spTan\x$, the tangent space of $\manif$ at $\x$. We define the normal space of $\manif$ at $\x$ by $\spNor\x \coloneqq (\spTan\x)^\perp$, the orthogonal complement of $\spTan\x$.
\end{definition}
For any $\x\in\manif$, we have $\spTan\x = \ker D\Phi(\x)$ where $\Phi$ is the local defining function for $\manif$ at $\x$. In the rest of the paper we will refer a $C^k$-smooth $m$-dimensional submanifold by simply $C^k$-smooth manifold. Furthermore the natural embedding of $\manif$ in $\spVar$ allows us to define a Riemannian metric on $\manif$, making it a Riemannian manifold. We define the functions $ \map{\projTan, \projNor}{\manif}{\spLin (\spVar, \spVar)} $ that provides for any $\x\in\manif$, the projection onto the tangent space $\projTan (\x) \coloneqq \proj{\spTan{\x}}$ and the projection onto the normal space $\projNor (\x) \coloneqq \proj{\spNor{\x}}$.
\begin{definition}[Riemannian Gradient]
    For any $k\in\N$, let $\manif\subset\spVar$ be a $C^k$-smooth manifold, $\map{f}{\manif}{\R}$ be a function and $\x\in\manif$. We say that $f$ is $C^k$-smooth at $\x$ if there exist a neighbourhood $X\subset\spVar$ of $\x$ and a $C^k$-smooth function $\map{\tilde{f}}{X}{\R}$ such that $\tilde{f}$ agrees with $f$ on $\manif \cap X$. In this case, we call $\tilde{f}$ a smooth extension of $f$ around $\x$. We call $\grad[\manif] f(\x) \in\spTan\x$, the Riemannian gradient of $f$ at $\x$ if for all $\v\in\spTan\x$, $\scal{\grad[\manif] f(\x)}{\v} = (f \circ \gamma )^\prime (0)$, where $\map{\gamma}{(-\epsilon, \epsilon)}{\manif}$ is any $C^1$-curve with $\gamma (0) = \x$ and $\dot\gamma (0) = \v$.
\end{definition}
The Riemannian Gradient of $f$ can alternatively be expressed in terms of the gradient of the smooth extension $\tilde{f}$ of $f$ by $\grad[\manif] f (\x) = \projTan (\x) \grad \tilde{f} (\x)$. Note that this gradient neither depends on the choice of the curve $\gamma$ nor the smooth extension $\tilde{f}$.
\begin{definition}[Riemannian Hessian]
    Let $\manif\subset\spVar$ be a $C^2$-smooth manifold, $\map{f}{\manif}{\R}$ be a $C^2$-smooth function and $\x\in\manif$. We call $\map{\Hess[\manif] f(\x)}{\spTan\x}{\spTan\x}$, the Riemannian Hessian of $f$ at $\x$ if for all $\v\in\spTan\x$, $\scal{\Hess[\manif] f(\x) \v}{\v} = (f \circ \gamma )^{\prime\prime} (0)$, where $\map{\gamma}{(-\epsilon, \epsilon)}{\manif}$ is any $C^1$-curve with $\gamma (0) = \x$ and $\dot\gamma (0) = \v$.
\end{definition}
We can similarly express the Riemannian Hessian $\Hess[\manif] f(\x) \in \spLin (\spTan\x, \spTan\x)$ by using the smooth extension $\tilde{f}$, that is,
\begin{equation} \label{eq:Riem:hess}
    \Hess[\manif] f (\x) = \projTan (\x) \Hess \tilde{f} (\x) + \wein[\x]{\cdot}{\projNor (\x) \grad \tilde{f} (\x)} \,,
\end{equation}
which requires the mapping $\wein[\x]{\cdot}{\w} \in \spLin (\spTan\x, \spTan\x)$ for $\w\in\spNor\x$ is called the Weingarten map and is defined by $\v \mapsto \wein[\x]{\v}{\w} \coloneqq - \projTan (\x) \mathrm{d} W [\v]$, where $W$ is a local extension of $\w$ to a normal vector field on $\manif$. $\wein[\x]{\cdot}{\w}$ is independent of the choice of normal field $W$ \cite[Proposition~II.2.1]{Cha06} and caters for the change of the tangent space as we move away from $\x$. It vanishes when the manifold is affine, that is, $\manif = \x + \spTan\x $ reducing the Hessian expression to $\Hess[\manif] f (\x) = \projTan (\x) \Hess \tilde{f} (\x)$.

\subsection{Partial Smoothness} \label{ssec:PS}
Partial smoothness \cite{Lew02} defines a large class of functions including the loss and regularization functions \cite[Section~3]{VDP+17} that are typically used in Machine Learning. The definition lifts key properties for handling constraint sets (for example, active set approaches) to the world of non-smooth functions in more generality.
\begin{definition}[Partial Smoothness] \label{def:part:C2}
    Let $\map{f}{\spVar}{\eR}$ be proper and lower semi-continuous and $\manif\subset\spVar$ be a set. We say that $f$ is partly smooth at a point $\x\in\manif$ relative to $\manif$ if the following conditions hold:
    \begin{enumerate}[label=(\roman*)]
        \item (Regularity:) $f$ is regular at every point close to $\x$ and $\partial f (\x) \neq \emptyset$.
        \item (Smoothness:) $\manif$ is a $C^2$-smooth manifold and $\rstDom f \manif$ is $C^2$-smooth around $\x$.
        \item (Sharpness:) $\spNor\x = \spPar \partial f(\x)$.
        \item (Continuity:) $\partial f$ is continuous at $\x$ relative to $\manif$.
    \end{enumerate}
    We call $f$ partly smooth relative to $\manif$ if $f$ is partly smooth at every $\x\in\manif$ relative to $\manif$.
\end{definition}
\begin{remark}
The first condition in Definition~\ref{def:part:C2} is automatically satisfied when $f$ is convex on $\ri (\dom f)$ and $\x \in \ri (\dom f)$. However when $f$ is non-convex, this assumption needs to be verified. For more on regularity, we refer the reader to \cite[Chapter~8]{RW98}.
\end{remark}

\section{Problem Setting} \label{sec:PS}

We aim for a theoretically sound and efficient strategies for computing the derivative of the solution of \eqref{prob:min} for partly smooth functions with respect to the parameter $\u$. In this section we define the class of functions and the necessary assumptions they must satisfy to obtain good estimates of derivative by applying the applying the methodologies from Section~\ref{ssec:intro:FPI} to the fixed-point iterations of APG. We also recall the Implicit Function Theorem of \cite{Lew02} which, although is not a practical way of differentiating the solution mapping, but is the backbone of our analysis.

Let $\manif\subset\spVar$ be a $C^2$-smooth manifold and $\setPrm\subset\spPrm$ be an open set. We consider the composite optimization problem \eqref{prob:min} where $\map{f, g}{\spVar\times\spPrm}{\eR}$ satisfy the following assumption. Examples of problems which satisfy this assumption are given in \cite{VDP+17}. A few examples are also provided in Section~\ref{sec:Exp} for numerical demonstration of our results.
\begin{assumption}[Convex Partly Smooth Objective] \label{ass:CPSO}
$\map{f}{\spVar\times\spPrm}{\R}$ is $C^2$-smooth, $\map{g}{\spVar\times\spPrm}{\eR}$ is partly smooth relative to $\manif\times\setPrm$ and for every $\u\in\setPrm$, $f(\cdot, \u)$ and $g(\cdot, \u)$ are convex and $f(\cdot, \u)$ has an $L$-Lipschitz continuous gradient.
\end{assumption}
\begin{remark}
\begin{enumerate}[label=(\roman*)]
    \item Since $\setPrm$ is open, the partial smoothness of $g$ relative to $\manif\times\setPrm$ implies $C^2$-smoothness of $g(\x, \cdot)$ for any $\x\in\manif$ \cite[Exercise~3.67]{Bou23} on $\setPrm$.
    \item The Lipschitz constant $L$ can be made dependent on $\u$ without changing the following analysis. We assume it to be independent of $\u$ for brevity.
    \item Assumption~\ref{ass:CPSO} encapsulates a large class of problems. However, it must be pointed out that many practical problems presented in Section~\ref{ssec:intro:apps} are not convex.
    \item Replacing $D_{\u}$, $\grad[\u]$ and $\grad[\u]^2$ with $D_{\Omega}$, $\grad[\Omega]$ and $\grad[\Omega]^2$ respectively due to the smoothness of $F$ w.r.t. $\u$ and using \cite[Exercise~3.67 and Example~5.19]{Bou23}, we express the Riemannian gradient and Hessian of $F$ at $(\x, \u)$ by
    \begin{equation*}
        \begin{aligned}
            \grad[\manif\times\Omega] F (\x, \u) &= \begin{bmatrix}
                \grad[\manif] F (\x, \u) \\ \grad[\u] F (\x, \u)
            \end{bmatrix} \,, \\
            \grad[\manif\times\Omega]^2 F (\x, \u) &= \begin{bmatrix}
                \Hess[\manif] F (\x, \u) & D_{\u} \grad[\manif] F (\x, \u) \\
                D_{\manif} \grad[\u] F (\x, \u) & \grad[\u]^2 F (\x, \u)
            \end{bmatrix} \,.
        \end{aligned}
    \end{equation*}
\end{enumerate}
\end{remark}

The next assumption that we require is reminiscent of conditions for the Implicit Function Theorem \cite[Theorem~1B.1]{DR09}. Let $\xmin \in \manif$ be a minimizer of $F(\cdot, \givenPrm)$ for $\givenPrm\in\setPrm$. A key requirement for using the Implicit Function Theorem for obtaining the derivative of the solution mapping at $\givenPrm$ is the positive definiteness of $\Hess[\manif] F (\xmin, \givenPrm)$. However, this alone will not suffice for proving the convergence of the derivative iterates of (Accelerated) Proximal-Gradient Descent on $\spTan\xmin$ (see Theorem~\ref{thm:APG:FAD}).
We additionally need positive definiteness of $\Hess[\x] f (\xmin, \u)$. For this purpose, we impose the following assumption on $f$ and $g$ at some $(\xmin, \givenPrm) \in \manif\times\setPrm$.
\begin{assumption}[Restrictive Positive Definiteness] \label{ass:RPD} Let $(\xmin, \givenPrm) \in  \manif\times\setPrm$ be a given point.
\begin{enumerate}[label=(\roman*)]
    \item \label{itm:RPD-i}The Hessian $\Hess[\manif] F (\xmin, \givenPrm)$ is positive definite on $\spTan\xmin $, that is,
    \begin{equation} \tag{RPD-i} \label{eq:RPD-i}
        \grad[\manif] F (\xmin, \givenPrm) \succ 0 \,.
    \end{equation}
    \item \label{itm:RPD-ii}Moreover, $\projTan (\xmin) \Hess[\x] f (\xmin, \givenPrm)$ is positive definite on $\spTan{\xmin}$, that is,
    \begin{equation} \tag{RPD-ii} \label{eq:RPD-ii}
        \rstDom{\projTan (\xmin) \Hess[\x] f (\xmin, \givenPrm)}{\spTan\xmin} \succ 0 \,.
    \end{equation}
\end{enumerate}
\end{assumption}
\begin{remark} \label{rem:RPD}
    \begin{enumerate}[label=(\roman*)]
        \item \label{itm:strong:min} The condition \eqref{eq:RPD-i} is equivalent to $\xmin$ being a strong local minimizer of $\rstDom{F(\cdot, \u)}{\manif}$ \cite[Definition~7, Lemma~4]{VDP+17}.
        \item We may need only one or both of the restricted positive-definiteness assumptions, that is, \eqref{eq:RPD-i} and \eqref{eq:RPD-ii} to hold at a point $(\xmin, \givenPrm)$. When both the assumptions are required to hold, we will simply say ``Assumption~\ref{ass:RPD} is satisfied at $(\xmin, \givenPrm)$''.
    \end{enumerate}
\end{remark}
Finally the solution of \eqref{prob:min} must change in a stable manner as the parameter $\u$ changes. This is guaranteed (see Theorem~\ref{thm:IFT}) by the following assumption on $f$ and $g$ at some $\givenPrm\in\setPrm$ and $\xmin \in \argmin_{\x\in\spVar} F(\x, \givenPrm)$, which is a standard condition when working with partly-smooth functions.
\begin{assumption}[Non-degeneracy] \label{ass:ND}
For a given point $(\xmin, \givenPrm)\in\manif\times\setPrm$, the non-degeneracy condition is satisfied, that is,
\begin{equation} \tag{ND} \label{eq:ND}
    0 \in \ri \partial_{\x} F (\xmin, \givenPrm) \,.
\end{equation}
\end{assumption}
\begin{remark} \label{rem:ND}
\begin{enumerate}[label=(\roman*)]
    \item In Assumption~\ref{ass:RPD}, $\xmin$ can be any arbitrary point. But since we always use Assumption~\ref{ass:RPD} together with Assumption~\ref{ass:ND}, we have $\xmin \in \argmin_{\x} F(\x, \givenPrm)$.
    \item The requirement in \eqref{eq:ND} is stronger than the Fermat's rule, that is, $0 \in \partial_{\x} F (\xmin, \givenPrm)$
    % \[0 \in \partial_{\x} F (\xmin, \givenPrm) \,,\]
    and therefore automatically implies optimality of $\xmin$.
    \item Assumptions~\ref{ass:RPD} and \ref{ass:ND} are required to be satisfied by $f$ and $g$ only at a specific point in $\manif\times\setPrm$, unlike Assumption~\ref{ass:CPSO}.
    \item \label{itm:strong:crit} Lewis~\cite{Lew02} calls $\xmin\in\manif$ a strong critical point of $F(\cdot, \u)$ relative to $\manif$ when \eqref{eq:RPD-i} and \eqref{eq:ND} are satisfied by $F$ at $(\xmin, \givenPrm)$.
    \item When $\manif$ is an affine manifold, any point $(\xmin, \givenPrm)$ that satisfies Assumption~\ref{ass:RPD}\ref{itm:RPD-ii} also satisfies Assumption~\ref{ass:RPD}\ref{itm:RPD-i}. However, when $\manif$ is not affine, then it is possible that a point satisfies only one of these two assumptions.
\end{enumerate}
\end{remark}
Under these assumptions, we can differentiate the solution mapping of \eqref{prob:min} on an open subset of $\setPrm$. The following result is derived from \cite{Lew02, VDP+17, Rii20}.
\begin{theorem}[Differentiation of Solution Map] \label{thm:IFT}
Let $f$ and $g$ satisfy Assumption~\ref{ass:CPSO} and $(\xmin, \givenPrm)\in\manif\times\setPrm$ be such that Assumptions~\ref{ass:RPD}\ref{itm:RPD-i} and \ref{ass:ND} are satisfied. Then there exist an open neighbourhood $\neighbourhoodPRM\subset\setPrm$ of $\givenPrm$ and a continuously differentiable mapping $\map{\argmap}{\neighbourhoodPRM}{\manif}$ such that for all $\u\in \neighbourhoodPRM$,
\begin{enumerate}[label=(\roman*)]
    \item $\argmap(\u)$ is the unique minimizer of $\rstDom{F(\cdot, \u)}{\manif}$, \label{itm:IFT:i}
    \item \eqref{eq:ND} and \eqref{eq:RPD-i} are satisfied at $(\argmap (\u), \u)$, and \label{itm:IFT:ii}
    \item the derivative of $\argmap$ is given by
    \begin{equation} \label{eq:IFT}
        D \argmap (\u) = -\Hess[\manif] F (\argmap (\u), \u)^{\dagger} D_{\u} \grad[\manif] F (\argmap (\u), \u) \,,
    \end{equation}
    where $\Hess[\manif] F (\argmap (\u), \u)^{\dagger}$ denotes the pseudoinverse of $\Hess[\manif] F (\argmap (\u), \u)$. \label{itm:IFT:iii}
\end{enumerate}
\end{theorem}
\findProofIn{thm:IFT}
\begin{remark} \label{rem:IFT:ND}
Assumption~\ref{ass:ND} is necessary for the differentiability of the solution mapping in general. One can still study the sensitivity analysis of the solution mapping when \eqref{eq:ND} is violated \cite{FMP18, Rii20}, which requires other assumptions. However it must be noted that when $f$, $g$ and the set of all possible partial smoothness active manifolds associated to $g$ are definable in an o-minimal structure \cite{van98,Cos00}, the set where \eqref{eq:ND} does not hold, is of Lebesgue measure zero \cite[Theorem~3]{VDP+17}.
\end{remark}
Theorem~\ref{thm:IFT} elegantly extends Implicit Differentiation to the setting when the objective is partly smooth and provides a theoretically justified way to computing the derivative of the solution mapping. However, it is not practical in the sense that the Riemannian Gradient and Hessian in general cannot be obtained through autograd packages and computing them by-hand can become cumbersome. The remedy is to benefit from the resemblance of the update step of PGD or APG --- the algorithm being used to solve \eqref{prob:min} --- with \eqref{prob:min} and apply Theorem~\ref{thm:IFT} to obtain their derivative expressions. These expressions can be used to lay down a solid theoretical foundation for ID, AD and FPAD applied to PGD and APG for a large class of functions. In practice, the closed-form solutions of proximal mappings can be exploited by computing the derivatives directly through standard autograd packages \cite{ABC+16, PGM+19, BFH+18} without going into the domains of Riemannian geometry.

% ***************************
% >>>>> PROBLEM SETTING <<<<<
% ***************************
\section{Accelerated Proximal Gradient (APG)} \label{sec:APG}
As mentioned above, AD and FPAD rely on an iterative procedure, for which we rely on the state-of-the-art algorithm for solving \eqref{prob:min}. Introduced in \cite{BT09} as Fast Iterative Shrinkage/Thresholding Algorithm or FISTA, the algorithm is the generalization of the celebrated Nesterov Acceleration \cite{Nes83} to Proximal Gradient Descent \cite{LM79}. For convex problems of the form \eqref{prob:min}, APG in general exhibits convergence like $\O (1/k^2)$ in objective values as compared to PGD which converges like $\O(1/k)$ \cite{CP16}. In fact, APG is the optimal algorithm for the convex problems of type \eqref{prob:min} when we only have access to the first order information of the smooth component $f$ of the objective. Convergence of the iterates of APG has been established \cite{CD15} for a specific range of parameters.
\begin{Algorithm}[Accelerated Proximal Gradient (APG)] \ \label{alg:APG}
    \begin{itemize}
        \item \key{Initialization:} $\xz=\xzm\in \spVar$, $\u\in\spPrm$, $0 < \sslow\leq\ssup < 2/L$.
        \item \key{Parameter:} $(\alpha_k)_{k\in\N} \in [\sslow, \ssup]$ and $(\beta_k)_{k\in\N} \in [0, 1]$.
        \item \key{Update $k\geq 0$:}
        \begin{equation} \label{itr:APG} \tag{APG}
            \begin{aligned}
                \yk &\coloneqq (1+\beta_k ) \xk - \beta_k \xkm \\
                \wk &\coloneqq \yk - \alpha_k \grad[\x] f(\yk, \u) \\
                \xkp &\coloneqq \bwd[k] (\wk, \u) \,.
            \end{aligned}
        \end{equation}
    \end{itemize}
\end{Algorithm}
Algorithm~\ref{alg:APG} shows the update steps of APG. In \eqref{itr:APG}, $\map{\bwd}{\spVar\times\spPrm}{\spVar}$ is defined as
\begin{equation} \label{eq:prox}
    \bwd (\w, \u) \coloneqq \argmin_{\x\in\spVar} \alpha g(\x, \u) + \frac{1}{2} \norm{\x - \w}^2 \,.
\end{equation}
A good choice of $\beta_k$ which also guarantees the convergence of the iterates is $(k-1)/(k+q)$ with $q>2$ \cite{CD15}. Henceforth, we will refer to Algorithm~\ref{alg:APG} with $\beta_k = 0$ as PGD and Algorithm~\ref{alg:APG}, in general or with $\beta_k > 0$ as APG.

\subsection{Linear Convergence of APG}
Liang et al.~\cite{LFP14, LFP17} showed that for partly smooth objectives, under the restricted positive definiteness and non-degeneracy assumptions, the iterates $\xk$ generated by APG eventually lie on the manifold $\manif$. In particular, there exists $K\in\N$ such that for all $k\geq K$, $\xk\in\manif$. Furthermore, the iterates converge to the solution at a linear rate. The following lemma formally summarizes this so-called finite activity identification and linear convergence of APG.

\begin{lemma}[Activity Identification and Linear Convergence of APG] \label{lem:APG}
Let $f$ and $g$ satisfy Assumption~\ref{ass:CPSO} and $(\xmin, \givenPrm)\in\manif\times\setPrm$ be such that Assumption~\ref{ass:ND} is satisfied. For $\alpha_k\in[\sslow, \ssup]$ and $\beta_k\in[0, 1]$, let the sequence $\seq[k\in\N]{\xk (\givenPrm)}$ generated by Algorithm~\ref{alg:APG} converges to $\xmin$. Then there exists $K\in\N$, such that $\xk (\givenPrm)\in\manif$ for all $k\geq K$. Moreover when Assumption~\ref{ass:RPD}\ref{itm:RPD-ii} is also satisfied, $\alpha_k \to \alpha_*$ and $\beta_k \to \beta_*$ such that $-1/(1 + 2\beta_*) < \lambda_{\min} (D_{\x} \pgd[*] (\xmin, \givenPrm))$, then $\xk (\givenPrm)$ converge linearly to $\xmin$ with rate $\rho (D_{\z} \apg[*] (\xmin, \xmin, \givenPrm) \projTan (\xmin, \xmin))$.
\end{lemma}

\begin{remark}
\begin{enumerate}[label=(\roman*)] \label{rem:APG}
    \item The active manifold here is $\manif$ and not $\manif\times\setPrm$ from Assumption~\ref{ass:CPSO}.
    \item Liang et al.~\cite{LFP17} consider a more general inertial Forward Backward or iFB algorithm, which encompasses different variants of FISTA. We stick to the more well-known version, that is, Algorithm~\ref{alg:APG}. However, it should be noted that our results rely on the analysis in \cite{LFP17} and with little effort, can be extended to iFB algorithm.
    \item \label{itm:FIP:upper}When the assumptions of Lemma~\ref{lem:APG} are satisfied and additionally, the iterates of Algorithm~\ref{alg:APG} are such that $\partial_{\x} g(\xk, \givenPrm) \subset \rbd (\partial_{\x} g(\xmin, \givenPrm)) $ whenever $\xk\notin\manif$, then $\xk\in\manif$ for all
    \begin{equation} \label{ineq:FIP:upper}
        k \geq \frac{\norm{\xz - \xmin}^2}{\sslow^2 \dist (-\grad[\x] f (\xmin, \givenPrm), \rbd (\partial_{\x} g (\xmin, \u)))^2} \,,
    \end{equation}
    which gives us a quantitative bound on the active set identification \cite[Proposition~3.6]{LFP17}.
    \item Under the assumptions of Lemma~\ref{lem:APG}, PGD can be shown to be locally faster than APG \cite[see Section~4.4]{LFP17}.
    \item Liang et al.~\cite{LFP17} argue that Assumption~\ref{ass:ND} is almost a necessary condition for finite activity identification of $\manif$. Fadili et al.~\cite{FMP18} study a more general identification property for mirror-stratifiable convex functions when such an assumption is not satisfied.
\end{enumerate}
\end{remark}

\subsection{Differentiation of Update Mapping of APG}

For applying AD, ID and FPAD from Section~\ref{ssec:intro:FPI}, we need the $C^1$-smoothness of the update mappings of PGD and APG. We first define the update mapping of PGD, that is, $\map{\pgd}{\spVar\times\spPrm}{\spVar}$ by
\begin{equation} \label{eq:FixMap:PGD}
    \pgd (\x, \u) \coloneqq \bwd (\x - \alpha \grad[\x] f (\x, \u), \u) \,,
\end{equation}
and the udpate mapping of APG, that is, $\map{\apg}{\spVar\times\spVar\times\spPrm}{\spVar\times\spVar}$ by
\begin{equation} \label{eq:FixMap:APG}
    \apg (\x_1, \x_2, \u) \coloneqq \left(\pgd \left(\x_1 + \beta (\x_1 - \x_2), \u \right), \x_1 \right) \,.
\end{equation}
\begin{remark} \label{rem:PGD:APG}
    \begin{enumerate}[label=(\roman*)]
        \item \label{itm:PGD:APG:FPI} The mappings $\pgd$ and $\apg$ yield more compact expressions for the update steps of PGD, that is, $\xkp \coloneqq \pgd[k] (\xk, \u)$ and APG, that is, $\zkp \coloneqq \apg[k] (\zk, \u)$ where $\zk\coloneqq (\xk, \xkm)$.
        \item \label{itm:PGD:APG:FPE} Let $f$ and $g$ be proper, lower semi-continuous, and convex and $f$ is $C^1$-smooth. Then for any $\alpha\in\R$ and $(\x, \u)\in\spVar\times\spPrm$ with $\x \in \argmin (f + g) (\cdot, \u)$ we have $\x = \pgd (\x, \u)$ and for any $\beta\in\R$ and $\z\coloneqq(\x_1, \x_2)\in\spVar\times\spVar$, we have
        \begin{equation*}
            \z = \apg (\z, \u) \quad \iff \quad \x_1 = \x_2 \textnormal{ and } \x_1 = \pgd (\x_1, \u) \,.
        \end{equation*}
    \end{enumerate}
\end{remark}
Since $\pgd$ and $\apg$ are defined through the proximal mapping $\bwd$ which by definition is a solution mapping of an optimization problem (see Equation~\ref{eq:prox}), we take advantage of Theorem~\ref{thm:IFT}.
\begin{lemma} \label{lem:prox:diff}
    Let $f$ and $g$ satisfy Assumption~\ref{ass:CPSO} and $(\xmin, \givenPrm)\in\manif\times\setPrm$ be such that Assumption~\ref{ass:ND} is satisfied. For any $\alpha_{*}\in(0, 2/L)$, there exists a neighbourhood $\neighbourhoodPRX[*]\subset\spVar\times\setPrm\times(0, 2/L)$ of $(\xmin - \alpha_{*} \grad[\x] f (\xmin, \givenPrm), \givenPrm, \alpha_{*})$ such that the mapping $(\w, \u, \alpha) \mapsto \bwd (\w, \u)$ defined in \eqref{eq:prox} is $C^1$-smooth on $\neighbourhoodPRX[*]$, $\bwd (\w, \u) \in \manif$ and
    \begin{equation} \label{eq:prox:derv}
        \begin{aligned}
            D_{\x} \bwd (\w, \u) &= \big ( \projTan (\x) + \wein[\x]{\cdot}{\projNor (\x) (\x - \w)} + \alpha \Hess[\manif] g(\x, \u) \big ) ^\dagger \\
            D_{\u} \bwd (\w, \u) &= - \alpha D_{\x} \bwd (\w, \u) D_{\u} \grad[\manif] g(\x, \u) \,,
        \end{aligned}
    \end{equation}
    for all $(\w, \u, \alpha)\in\neighbourhoodPRX[*]$ and $\x \coloneqq \bwd (\w, \u)$.
\end{lemma}
\findProofIn{lem:prox:diff}
\begin{remark} \label{rem:prox:ND}
    Assumption~\ref{ass:RPD} is not required in the hypothesis of Lemma~\ref{lem:prox:diff} for differentiating $\bwd$ near $(\xmin - \alpha \grad[\x] f(\xmin, \givenPrm), \givenPrm)$.
\end{remark}
Using Lemma~\ref{lem:prox:diff} and the chain rule, we can differentiate the update step of PGD.
\begin{corollary} \label{cor:PGD:Derv}
    Let $f$ and $g$ satisfy Assumption~\ref{ass:CPSO} and $(\xmin, \givenPrm)\in\manif\times\setPrm$ be such that Assumption~\ref{ass:ND} is satisfied. For any $\alpha_{*}\in(0, 2/L)$, there exists a neighbourhood $\neighbourhoodPGD[*]\subset\spVar\times\setPrm\times(0, 2/L)$ of $(\xmin, \givenPrm, \alpha_*)$ such that the mapping $(\x, \u, \alpha) \mapsto \pgd (\x, \u)$ defined in \eqref{eq:FixMap:PGD} is $C^1$-smooth on $\neighbourhoodPGD[*]$, $\pgd (\x, \u) \in \manif$ and
    \begin{equation}
        \begin{aligned}
            D_{\x} \pgd (\x, \u) &= \Dxbwd (\x, \u)^{\dagger} \left( \opid - \alpha \Hess[\x] f (\x, \u) \right) \\
            D_{\u} \pgd (\x, \u) &= -\alpha \Dxbwd (\x, \u)^{\dagger} D_{\u} \left( \grad[\x] f (\x, \u) + \grad[\manif] g (\pgd (\x, \u), \u) \right) \,,
        \end{aligned}
    \end{equation}
    for all $(\x, \u, \alpha)\in\neighbourhoodPGD[*]$ where $\Dxbwd (\x, \u) \coloneqq D_{\x} \bwd (\x - \alpha \grad[\x] f (\x, \u), \u)^{\dagger}$.
\end{corollary}
\findProofIn{cor:PGD:Derv}
Similarly for the mapping $\apg$ we have the following result
\begin{corollary} \label{cor:APG:Derv}
    Let $f$ and $g$ satisfy Assumption~\ref{ass:CPSO} and $(\xmin, \givenPrm)\in\manif\times\setPrm$ be such that Assumption~\ref{ass:ND} is satisfied. For any $\alpha_{*}\in(0, 2/L)$, there exists a neighbourhood $\neighbourhoodAPG[*]\subset\spVar\times\spVar\times\setPrm\times(0, 2/L)$ of $(\xmin, \xmin, \givenPrm, \alpha_*)$ such that the mapping $(\z, \u, \alpha) \mapsto \apg(\z, \u)$ defined in \eqref{eq:FixMap:APG} is $C^1$-smooth on $\neighbourhoodAPG[*]$, $ \apg(\x_1, \x_2, \u) \in \manif\times\spVar$ and
    \begin{equation} \label{eq:APG:Derv}
        \begin{aligned}
            D_{\z} \apg (\z, \u) &= \begin{bmatrix}
                (1+\beta) \Dxpgd (\z, \u) & -\beta \Dxpgd (\z, \u) \\
                \opid & 0
            \end{bmatrix} \\
            D_{\u} \apg (\z, \u) &= \begin{bmatrix}
                \Dupgd (\z, \u) \\
                0
            \end{bmatrix} \,,
        \end{aligned}
    \end{equation}
    for all $(\x_1, \x_2, \u, \alpha, \beta)\in\neighbourhoodAPG[*] \times [0, 1]$ where $\z\coloneqq (\x_1, \x_2)$, $\Dxpgd (\z, \u) \coloneqq D_{\x} \pgd (\x_1 + \beta (\x_1 - \x_2), \u)$ and $\Dupgd (\z, \u) \coloneqq D_{\u} \pgd (\x_1 + \beta (\x_1 - \x_2), \u)$.
\end{corollary}
\findProofIn{cor:APG:Derv}

\section{Approaches for Differentiating Fixed-Point Mappings of APG} \label{sec:PDM}

Now that we have the expressions for computing the derivative of the update mappings of PGD and APG, we are well-equipped to apply AD, ID and FPAD from Section~\ref{ssec:intro:FPI} to Algorithm~\ref{alg:APG}. In this section, we analyze the three differentiation techniques and compare them with \eqref{eq:IFT}. In many practical applications, analytic solution to the proximal mapping is known and is expressed as compositions of elementary functions. Consequently $\bwd$ defined in \eqref{eq:prox} and hence the operations in Algorithm~\ref{alg:APG} can be differentiated by using standard autograd packages. Therefore applying differentiation techniques to Algorithm~\ref{alg:APG} appears to be a more natural choice for estimating the derivative of solution mapping of \eqref{prob:min} as opposed to \eqref{eq:IFT}.

\subsection{Implicit Differentiation} \label{ssec:PDM:ID}
In this section we apply ID on the fixed-point equations of PGD and APG and recognize it to be more practical than simply evaluating the expression in \eqref{eq:IFT}.

\subsubsection{Implicit Differentiation of PGD} \label{sssec:PDM:ID:PGD}
We first apply ID on the fixed-point equation of PGD, that is, $\x = \pgd (\x, \u)$.
\begin{theorem} \label{thm:PGD:IFT}
    Let $f$ and $g$ satisfy Assumption~\ref{ass:CPSO} and $(\xmin, \givenPrm)\in\manif\times\setPrm$ be such that Assumptions~\ref{ass:RPD}\ref{itm:RPD-ii} and \ref{ass:ND} are satisfied. Then for any $\alpha\in[\sslow, \ssup]$, $\rho (D_{\x} \pgd (\xmin, \givenPrm) \projTan (\xmin)) < 1$. Additionally when Assumption~\ref{ass:RPD}\ref{itm:RPD-i} is also satisfied, the (possibly reduced) neighbourhood $\neighbourhoodPRM$ and the mapping $\argmap$ from Theorem~\ref{thm:IFT} satisfy $\x = \pgd (\x, \u)$ and
    \begin{equation} \label{eq:PGD:IFT}
        D \argmap (\u) = \left( \opid - D_{\x} \pgd (\x, \u) \projTan (\x) \right)^{-1} D_{\u} \pgd (\x, \u) \,,
    \end{equation}
    for all $\u\in \neighbourhoodPRM$ and $\x \coloneqq \argmap (\u)$.
\end{theorem}
\findProofIn{thm:PGD:IFT}
\begin{remark} \label{rem:PGD:IFT:Bad}
    Since the restricted positive definiteness assumption in \eqref{eq:RPD-ii} penalizes $\rho (D_{\x} \pgd (\xmin, \givenPrm) \projTan (\xmin))$ instead of $\rho (D_{\x} \pgd (\xmin, \givenPrm))$, the additional $\projTan (\argmap (\u))$ in \eqref{eq:PGD:IFT} may not be ignored. However, when the restricted positive definiteness is replaced with $\Hess[\x] f(\xmin, \givenPrm) \succ 0$, we can modify the above proof to obtain $\rho (D_{\x} \pgd (\xmin, \givenPrm)) < 1$ and simply apply Theorem~\ref{thm:IFT:FPE} without relying on \eqref{eq:RPD-i} and Theorem~\ref{thm:IFT}.
\end{remark}
We rely on Assumption~\ref{ass:RPD}\ref{itm:RPD-i} to ensure that a $C^1$-smooth mapping exists on a neighbourhood of $\givenPrm$ via Theorem~\ref{thm:IFT}. Without this assumption, it is not trivial to show that such a mapping exists by using Assumption~\ref{ass:RPD}\ref{itm:RPD-ii} alone and invoking Theorem~\ref{thm:IFT:FPE} on $\x = \pgd (\x, \u)$. However, since Assumption~\ref{ass:RPD}\ref{itm:RPD-i} only provides a sufficient condition for the existence of $C^1$-smooth solution mapping $\argmap$ on the neighbourhood of $\givenPrm$, we can still use \eqref{eq:PGD:IFT} to furnish the derivative of $\argmap$ when Assumption~\ref{ass:RPD}\ref{itm:RPD-ii} holds while Assumption~\ref{ass:RPD}\ref{itm:RPD-i} does not.
\begin{corollary} \label{cor:PGD:IFT}
    Let $f$ and $g$ satisfy Assumption~\ref{ass:CPSO} and $(\xmin, \givenPrm)\in\manif\times\setPrm$ be such that Assumptions~\ref{ass:RPD}\ref{itm:RPD-ii} and \ref{ass:ND} are satisfied. Moreover, for some neighbourhood $\neighbourhoodPRM$ of $\givenPrm$, let $\map{\argmap}{\neighbourhoodPRM}{\spVar}$ be a $C^1$-smooth mapping such that $\argmap (\givenPrm) = \xmin$ and $\neighbourhoodPRM\ni\u\mapsto\argmap (\u) = \argmin_{\x} F (\x, \u)$. Then for some neighbourhood $V\subset \neighbourhoodPRM$ of $\givenPrm$, $\argmap (V)\subset\manif$ and the derivative of $\rstDom{\argmap}{V}$ is given by \eqref{eq:PGD:IFT}.
\end{corollary}
\findProofIn{cor:PGD:IFT}

\subsubsection{Implicit Differentiation of APG} \label{sssec:PDM:ID:APG}
We similarly apply Implicit Differentiation on APG by making use of Theorem~\ref{thm:PGD:IFT}. However, before we state our result, we recall the definition of the orthogonal projection mapping for the manifold $\manif\times\manif$ onto the tangent space $T_{\z} (\manif \times \manif) = \spTan{\x_1} \times \spTan{\x_2}$ \cite[Proposition~3.20]{Bou23} at a given a point $\z\coloneqq (\x_1, \x_2)\in\manif\times\manif$, that is,
\begin{equation*}
    \projTan (\z) \coloneqq \proj{T_{\z} (\manif \times \manif)} = \begin{bmatrix}
        \projTan (\x_1) & 0 \\
        0 & \projTan (\x_2)
    \end{bmatrix} \,.
\end{equation*}
We also state a preliminary result which will be useful for providing the expression for the inverse of the Jacobian occurring in the Implicit Differentiation of APG.
\begin{lemma} \label{lem:spec_props:M}
    Let $R \in \spLin (\spVar, \spVar)$ and $\projTan \coloneqq \proj{\im R}$ be such that $\rho (R \projTan) < 1$. Then for any $\beta\in[0, 1]$ and $M_{\beta} \in \spLin (\spVar\times\spVar, \spVar\times\spVar)$ defined by
    \begin{equation*}
        M_{\beta} = \begin{bmatrix}
            (1 + \beta) R \projTan & -\beta R \projTan \\
            \projTan & 0
        \end{bmatrix} \,,
    \end{equation*}
    we have $\rho (M_{\beta}) < 1$ if and only if $-1/(1 + 2\beta) < \lambda_{\min} (R)$. Moreover, when $\rho (M_{\beta}) < 1$, we have
    \begin{equation*}
        \left( \opid - M_{\beta} \right)^{-1} = \begin{bmatrix}
            \left(\opid - R \projTan \right)^{-1} & -\beta\left(\opid - R \projTan \right)^{-1}R \projTan \\
            \projTan \left(\opid - R \projTan \right)^{-1} & \opid - \beta \projTan \left(\opid - R \projTan \right)^{-1}R \projTan
        \end{bmatrix} \,.
    \end{equation*}
\end{lemma}
\findProofIn{lem:spec_props:M}
\begin{remark} \label{rem:spec_props:M}
    \begin{enumerate}[label=(\roman*)]
        \item \label{itm:spec_props:M:R} It is crucial to understand $\rho (M_{\beta})$ and $\rho (R)$ for analyzing and comparing the convergence rates of the APG and PGD as well as those of their AD and FPAD counterparts as we will see in the next sections. We refer the reader to \cite[Section~4]{LFP17} for a detailed discussion on the spectral properties of $M_{\beta}$ and $R$.
        \item \label{itm:spec_props:M:beta=0} When $\beta = 0$, it is straightforward to show that any non-zero eigenvalue of $M_{\beta}$ is also an eigenvalue of $R\projTan$ suggesting that $\rho (M_{0}) = \rho (R)$.
    \end{enumerate}
\end{remark}
We now apply ID on the fixed-point equation of APG, that is, $\z = \apg (\z, \u)$.
\begin{theorem} \label{thm:APG:IFT}
    Let $f$ and $g$ satisfy Assumption~\ref{ass:CPSO} and $(\xmin, \givenPrm)\in\manif\times\setPrm$ be such that Assumptions~\ref{ass:RPD}\ref{itm:RPD-ii} and \ref{ass:ND} are satisfied. Then for any $\alpha\in[\sslow, \ssup]$ and $\beta\in[0, 1]$ with $-1/(1 + 2\beta) < \lambda_{\min} (D_{\x} \pgd (\xmin, \u))$, we have $\rho (D_{\z} \apg (\xmin, \xmin, \givenPrm) \projTan (\xmin, \xmin)) < 1$. Additionally when Assumption~\ref{ass:RPD}\ref{itm:RPD-i} is also satisfied, the (possibly reduced) neighbourhood $\neighbourhoodPRM$ and the mapping $\argmap$ from Theorem~\ref{thm:IFT} satisfy $\z = \apg (\z, \u)$ and
    \begin{equation} \label{eq:APG:IFT}
        \begin{bmatrix}
            D \argmap (\u) \\
            D \argmap (\u)
        \end{bmatrix} = \left( \opid - D_{\z} \apg (\z, \u) \projTan (\z)\right)^{-1} D_{\u} \apg (\z, \u) \,,
    \end{equation}
    for all $\u\in \neighbourhoodPRM$ and $\z \coloneqq (\argmap (\u), \argmap (\u))$.
\end{theorem}
\findProofIn{thm:APG:IFT}
\begin{remark}
    Following the approach employed for PGD in Corollary~\ref{cor:PGD:IFT}, an analogous argument can be made for APG by using the results of Theorem~\ref{thm:APG:IFT}.
\end{remark}

\subsubsection{Implicit Differentiation near the Solution} \label{sssec:PDM:ID:NS}
In general, we do not have access to the solution mapping $\argmap$ as in Theorem~\ref{thm:IFT} and only an approximate can be acquired through Algorithm~\ref{alg:APG}. Therefore, our only option is to compute the approximate of $D\argmap (\u)$ through \eqref{eq:PGD:IFT} and \eqref{eq:APG:IFT}. We first consider the case of PGD. The presence of projection $\projTan (\argmap (\u))$ in \eqref{eq:PGD:IFT} suggests two ways of evaluating this expression at some $\x$ near $\argmap (\u)$, that is, we replace $\projTan (\argmap (\u))$ with $\projTan (\x^{\prime})$ for either $\x^{\prime} = \x$ or $\x^{\prime} = \pgd (\x, \u)$. The former may seem more natural however, the latter is more useful because it does not require $\x$ to lie in $\manif$ because $\pgd (\x, \u)\in\manif$ from Corollary~\ref{cor:PGD:Derv}. Additionally it also does not require us to explicitly evaluate $\projTan (\x^{\prime})$ when using FPAD, as we will see in Section~\ref{ssec:PDM:FPAD}. In either case, when Assumptions~\ref{ass:CPSO}, \ref{ass:RPD}\ref{itm:RPD-ii} and \ref{ass:ND} are satisfied, we argue from Corollary~\ref{cor:PGD:Derv}, Theorem~\ref{thm:PGD:IFT} and Lemma~\ref{lem:FPE:neighbourhood} that for any $\alpha\in(0, 2/L)$, there exists a neighbourhood $\neighbourhoodPGD$ of $(\xmin, \givenPrm)$ (relative to $\manif\times\setPrm$ when $\x^{\prime}=\x$) such that $\rho (D_{\x} \pgd (\x, \u) \projTan (\x^{\prime})) < 1$ for all $(\x, \u)\in\neighbourhoodPGD$, which gives us an expression similar to \eqref{eq:FPE:IFT:approx}, that is,
\begin{equation} \label{eq:PGD:IFT:approx}
    \pgdPhi (\x, \u) \coloneqq \left( \opid - D_{\x} \pgd (\x, \u) \projTan (\x^{\prime}) \right)^{-1} D_{\u} \pgd (\x, \u) \,.
\end{equation}
Additionally for any $\beta\in[0, 1]$, we similarly evaluate \eqref{eq:APG:IFT} at $\z\coloneqq(\x_1, \x_2)$ near $\z^*\coloneqq(\xmin, \xmin)$ by replacing $\projTan (\z^*)$ with $\projTan (\z^{\prime})$ for either $\z^{\prime}=\z$ or $\z^{\prime}=\apg (\z, \u)$. This suggests that we must have $\x_1\in\manif$ for both cases. Moreover, we must have $\x_1+\beta(\x_1-\x_2)\in\manif$ for $\z^{\prime} = \z$ which holds only when $\manif$ is affine and requires $\x_1=\x_2$ otherwise. We therefore argue that there exists a neighbourhood $\neighbourhoodAPG \subset\manif\times\spVar\times\setPrm$ of $(\z^*, \givenPrm)$ when $\z^{\prime}=\apg (\z, \u)$ or a neighbourhood $\neighbourhoodAPG \subset\manif\times\setPrm$ of $(\xmin, \givenPrm)$ when $\z^{\prime}=\z$, such that for all $(\z, \u)\in\neighbourhoodAPG$ when $\z^{\prime}=\apg (\z, \u)$ or for all $(\x, \u) \in\neighbourhoodAPG$ with $\x_1=\x_2=\x$, $\rho (D_{\z} \apg (\z, \u) \projTan (\z^{\prime}) < 1$ and
\begin{equation} \label{eq:APG:IFT:approx}
    \apgPhi (\z, \u) \coloneqq \left( \opid - D_{\z} \apg (\z, \u) \projTan (\z^{\prime}) \right)^{-1} D_{\u} \apg (\z, \u) \,.
\end{equation}
By applying Lemma~\ref{lem:spec_props:M}, we observe that $\apgPhi (\z, \u)$ simplifies to
\begin{equation*}
    \apgPhi (\z, \u) = \begin{bmatrix}
        \pgdPhi (\x_1 + \beta (\x_1 - \x_2), \u) \\
        \pgdPhi (\x_1 + \beta (\x_1 - \x_2), \u)
    \end{bmatrix} \,.
\end{equation*}
In practice, we set $\z\coloneqq (\x, \x)$ for computing $\apgPhi (\z, \u)$ because we have only one approximation $\x$ of fixed-point. The mappings $\pgdPhi$ and $\apgPhi$ provide meaningful quantities only when the solution mapping $\argmap$ is $C^1$-smooth which is the case when Assumption~\ref{ass:RPD}\ref{itm:RPD-i} also holds. Moreover, from continuity, they converge to $\pgdPhi (\argmap (\u), \u)$ and $\apgPhi (\argmap (\u), \argmap (\u), \u)$ both of which correspond to $D\argmap (\u)$.

We try to simplify the expression for $\pgdPhi$ and acquire a quantity which is analogous to \eqref{eq:IFT}. The problem we face in doing so is that we may not have $\x = \pgd (\x, \u)$ because $\x$ is not necessarily the solution, that is, $\x\neq\argmap (\u)$ (from Theorem~\ref{thm:IFT}), which is crucial in proving Theorem~\ref{thm:PGD:IFT} and, in particular, verifying the equality \eqref{eq:PGD:IFT}. To understand that we expand the terms occurring in $\pgdPhi$:
\begin{equation} \label{eq:PGD:maps}
    \begin{aligned}
        \Dxbwd (\x, \u) &= \projTan (\y) + \wein[\y]{\cdot}{\projNor (\y) (\y - \w)} + \alpha \Hess[\manif] g (\y, \u)) \\
        R_{\alpha} (\x, \u) &= \Dxbwd (\x, \u)^{\dagger} \left( \opid - \alpha \Hess[\x] f (\x, \u) \right) \\
        S_{\alpha} (\x, \u) &= -\alpha \Dxbwd (\x, \u)^{\dagger} D_{\u} \left( \grad[\x] f (\x, \u) + \grad[\manif] g (\y, \u) \right) \\
        \pgdPhi (\x, \u) &= \left(\opid - R_{\alpha} (\x, \u) \projTan (\x^{\prime}) \right)^{-1} S_{\alpha} (\x, \u) \,,
    \end{aligned}
\end{equation}
where $\y\coloneqq \pgd (\x, \u)$ and $\w \coloneqq \x - \alpha \grad[\x] f(\x, \u)$. In the above equation, it is evident that $\y \neq \x$ complicates the simplification of $\pgdPhi$. However, if we slightly change the direction in the forward step and define a modified update for PGD, that is,
\begin{equation*}
    \begin{aligned}
        \bnu (\x, \u) &\coloneqq \proj{\partial_{\x} g (\x, \u)} (-\grad[\x] f (\x, \u)) \\
        \tilde{\mathcal A}_{\alpha} (\x, \u) &\coloneqq \bwd (\x + \alpha \bnu (\x, \u), \u) \,,
    \end{aligned}
\end{equation*}
we obtain $\x=\tilde{\mathcal A}_{\alpha} (\x, \u)$. The new direction $\bnu (\x, \u)$ is well-defined thanks to the projection theorem \cite[Theorem~3.16]{BC11}. In the result below, we show that this modified gradient step converges to the true gradient step under the given assumptions.
\begin{lemma} \label{lem:sg:seq}
Let $f$ and $g$ satisfy Assumption~\ref{ass:CPSO} and $(\xmin, \givenPrm)\in\manif\times\setPrm$ be such that $\xmin\coloneqq\argmin F(\cdot, \givenPrm)$. Then for any sequence $(\xk)_{k\in\N}$ in $\manif$ converging to $\xmin$, the sequence $( \bnu (\xk, \givenPrm))_{k\in\N}$ converges to $-\grad[\x] f (\xmin, \givenPrm)$.
\end{lemma}
\findProofIn{lem:sg:seq}
Similarly, we define modified mappings by replacing $\y$ and $\x^{\prime}$ with $\x$ and $\w\coloneqq\x - \alpha \grad[\x] f(\x, \u)$ with $\w\coloneqq\x + \alpha \bnu (\x, \u)$ in the respective mappings in \eqref{eq:PGD:maps}, that is,
\begin{equation} \label{eq:PGD:mod:maps}
    \begin{aligned}
        \tilde Q_{\alpha} (\x, \u) &= \projTan (\x) - \alpha \wein[\x]{\cdot}{\projNor (\x) \bnu (\x, \u)} + \alpha \Hess[\manif] g (\x, \u)) \\
        \tilde R_{\alpha} (\x, \u) &= \tilde Q_{\alpha} (\x, \u)^{\dagger} \left( \opid - \alpha \Hess[\x] f (\x, \u) \right) \\
        \tilde S_{\alpha} (\x, \u) &= -\alpha \tilde Q_{\alpha} (\x, \u)^{\dagger} D_{\u} \left( \grad[\x] f (\x, \u) + \grad[\manif] g (\x, \u) \right) \\
        \tilde\phi_{\alpha} (\x, \u) &= \left(\opid - \tilde R_{\alpha} (\x, \u) \projTan (\x) \right)^{-1} \tilde S_{\alpha} (\x, \u) \,.
    \end{aligned}
\end{equation}
Using Lemma~\ref{lem:sg:seq}, we show that $\tilde\phi_{\alpha}$ is a well-defined map on a neighbourhood of $(\xmin, \givenPrm)$ and it shares similarities with \eqref{eq:IFT}.
\begin{theorem} \label{thm:PGD:IFT:approx}
    Let $f$ and $g$ satisfy Assumption~\ref{ass:CPSO} and $(\xmin, \givenPrm)\in\manif\times\setPrm$ be such that Assumptions~\ref{ass:RPD} and \ref{ass:ND} are satisfied. Then for any $\alpha\in[\sslow, \ssup]$, $\tilde\phi_{\alpha}$ is well-defined on a neighbourhood $\neighbourhoodPGD$ of $(\xmin, \givenPrm)$ relative to $\manif\times\setPrm$ and $\tilde\phi_{\alpha} (\x, \u) \to D \argmap (\u)$ as $\x \to \argmap (\u)$ where the (possibly reduced) neighbourhood $\neighbourhoodPRM$ and the mapping $\argmap$ are from Theorem~\ref{thm:IFT}. In particular,
    \begin{equation} \label{eq:PGD:mod:IFT:approx}
        \tilde\phi_{\alpha} (\x, \u) = -\left ( \Hess[\manif] F(\x, \u) + \mathcal{W}_{\alpha} (\x, \u) \right )^\dagger D_{\u} \grad[\manif] F (\x, \u) \,,
    \end{equation}
    where $\mathcal{W}_{\alpha} (\x, \u)\coloneqq-\alpha\wein[\x]{\cdot}{\projNor (\x) (\bnu (\x, \u) + \grad[\x] f (\x, \u))}$.
\end{theorem}
\findProofIn{thm:PGD:IFT:approx}
\begin{remark} \label{rem:IDNS}
    \begin{enumerate}[label=(\roman*)]
        \item We can provide a similar analysis to obtain an analogous expression for APG on a neighbourhood $\neighbourhoodPGD \subset \manif\times\setPrm$ of $(\xmin, \givenPrm)$ defined by
        \begin{equation} \label{eq:APG:mod:IFT:approx}
            \tilde\phi_{\alpha, \beta} (\x, \x, \u) = \begin{bmatrix}
                \tilde\phi_{\alpha} (\x, \u) \\
                \tilde\phi_{\alpha} (\x, \u)
            \end{bmatrix} \,,
        \end{equation}
        for all $(\x, \u)\in\neighbourhoodPGD$.
        \item \label{itm:IDNS:affine} When $\manif$ is an affine manifold, $\mathcal{W}_{\alpha} (\x, \u)$ in \eqref{eq:PGD:mod:IFT:approx} vanishes and the expression for $\tilde\phi_{\alpha}$ simplifies to
        \begin{equation} \label{eq:PGD:IFT:approx:affine}
            \tilde\phi_{\alpha} (\x, \u) = - \Hess[\manif] F(\x, \u) ^\dagger D_{\u} \grad[\manif] F (\x, \u) \,.
        \end{equation}
    \end{enumerate}
\end{remark}

\subsection{Automatic Differentiation} \label{ssec:PDM:AD}

As we have seen in Section~\ref{ssec:intro:FPI}, another competitor of ID in the context of fixed-point iterations is AD.

\subsubsection{Convergence of Automatic Differentiation}

Given $\u\in\setPrm$ and the sequences $\seq[k\in\N]{\alpha_k}$ and $\seq[k\in\N]{\beta_k}$ with limits $\alpha_*$ and $\beta_*$ respectively, let $\seq[k\in\N]{\xk (\u)}$ be the iterates of Algorithm~\ref{alg:APG} such that $\xk(\u)\to\argmap (\u)$. As long as the update mapping $\apg$ is differentiable at $ (\xk (\u), \xkm (\u), \u)$, AD of Algorithm~\ref{alg:APG} yields meaningful quantities. This is possible because eventually, $(\xk (\u), \xkm (\u), \u, \alpha_k) \in\neighbourhoodAPG[*]$ (see Corollary~\ref{cor:APG:Derv}). In this situation, the assumptions which guarantee the local linear convergence of Algorithm~\ref{alg:APG} also guarantee the convergence of AD of Algorithm~\ref{alg:APG} as demonstrated by the following result.
\begin{theorem} \label{thm:APG:FAD}
Let $f$ and $g$ satisfy Assumption~\ref{ass:CPSO} and $(\xmin, \givenPrm)\in\manif\times\setPrm$ be such that Assumptions~\ref{ass:RPD} and \ref{ass:ND} are satisfied. Let $[\sslow, \ssup]\ni\alpha_k \to \alpha_*$ and $[0, 1] \ni \beta_k \to \beta_*$, such that $-1/(1 + 2\beta_*) < \lambda_{\min} (D_{\x} \pgd[*] (\xmin, \givenPrm))$, the iterates of Algorithm~\ref{alg:APG} $\xk(\givenPrm)\to\xmin$ such that $(\xk, \xkm, \givenPrm, \alpha_k) \in\neighbourhoodAPG[*]$ for all $k\geq0$ where $\neighbourhoodAPG[*]$ is from Corollary~\ref{cor:APG:Derv}. Then the sequence $\seq[k \in\N]{\xkd (\givenPrm)}$ generated by forward mode AD of Algorithm~\ref{alg:APG} converges to $D \argmap(\givenPrm) \ud$.
\end{theorem}
\findProofIn{thm:APG:FAD}

\begin{remark} \label{rem:APG:FAD}
\begin{enumerate}[label=(\roman*)]
    \item The earlier iterates of Algorithm~\ref{alg:APG} do not lie on $\manif$ in general and we must resort to late-start (see Section~\ref{sssec:FPI:AD}) because eventually, the iterates $\xk$ satisfy the hypothesis of the above theorem thanks to Lemma~\ref{lem:APG}. The question of when to late start the AD process is very hard to answer in general. Fortunately, even without late-starting, one can still recover the correct estimate of $D\argmap(\givenPrm)$ because, as mentioned before, the proximal mappings are given by analytical expressions in practice and the AD libraries always provide a finite quantity at points of non-differentiability \cite{BP20, BP21}. In this case, the sequence $\xkd(\givenPrm)$ may exhibit an odd behaviour in the beginning, but it still converges to $D\argmap (\u)\ud$ thanks to finite identification and Theorem~\ref{thm:LGFPI:conv}.
    \item Assumption~\ref{ass:ND} is necessary for proving Theorem~\ref{thm:IFT} which lays the foundation for establishing the differentiablity of $\pgd$ and $\apg$. The lack of this assumption calls for the use of generalized derivatives. Riis~\cite{Rii20} showed $\bwd$ to be piecewise smooth under further assumptions while Bolte et al.~\cite{BPV22} employed conservative analysis \cite{BP20, BP21} on the sequence $\xk$ when the update mapping is path differentiable.
    \item As mentioned in Section~\ref{ssec:intro:RW}, we generalize the results of \cite{Rii20}. Moreover, there are slight errors in \cite{Rii20} that we fix here. In \cite[Lemma~7.31]{Rii20} and afterwards in the proofs, the Weingarten term of \eqref{eq:prox:derv} is missing. Also in \cite[Proposition~2.7]{Rii20}, for linear convergence of $(d^k)_{k\in\N}$, the sequences $(A_k)_{k\in\N}$ and $(b^k)_{k\in\N}$ ought to converge linearly.
    \item Theorem~\ref{thm:APG:FAD} provides convergence guarantees of the derivative sequences for a wider class of functions as compared to those in \cite{BKM+21}.
\end{enumerate}
\end{remark}

\subsubsection{Linear Convergence of Automatic Differentiation}
Because $\pgd[k]$ and $\apg[k]$ depend on $k$, the linear convergence results from Section~\ref{sssec:FPI:AD} can not be easily translated. In fact, it is not trivial to show linear convergence for AD of Algorithm~\ref{alg:APG} under the given assumptions. However, if we use constant step size and assume local Lipschitz continuity of the involved derivatives, the convergence rate guarantees can be provided for AD of PGD. The rates for AD of APG are harder to obtain in the general setting.
\begin{theorem} \label{thm:APG:FAD:Lin:Conv}
Let $f$ and $g$ satisfy Assumption~\ref{ass:CPSO} and $(\xmin, \givenPrm)\in\manif\times\setPrm$ be such that Assumptions~\ref{ass:RPD} and \ref{ass:ND} are satisfied. Let $\alpha_k \coloneqq \alpha_* \in [\sslow, \ssup]$ and $\beta_k\coloneqq 0$, $\xk(\givenPrm)\to\xmin$ and $(\xk (\givenPrm), \givenPrm, \alpha_*) \in\neighbourhoodPGD[*]$ for all $k\geq0$ where $\neighbourhoodPGD[*]$ is from Corollary~\ref{cor:PGD:Derv}. The sequence $( \xkd (\givenPrm) )_{k \in\N}$ generated by forward mode AD of Algorithm~\ref{alg:APG} converges linearly with rate $k\rho (D_{\z} \pgd[*] (\xmin, \givenPrm) \projTan (\xmin))$, whenever $\grad[\manif] g$, $\Hess[\x] f$, $\Hess[\manif] g$, $D_{\u} \grad[\x] f$ and $D_{\u} \grad[\manif] g$ are locally Lipschitz continuous near $(\xmin, \givenPrm)$.
\end{theorem}
\findProofIn{thm:APG:FAD:Lin:Conv}
\begin{remark} \label{rem:APG:FAD:Lin:Conv}
    Under the given assumptions, the rate of convergence of AD of Algorithm~\ref{alg:APG}, that is, $k \rho (D_{\z} \pgd[*] (\xmin, \givenPrm) \projTan (\xmin))$ is borrowed from rate of convergence of the algorithm itself, that is, $\rho (D_{\z} \pgd[*] (\xmin, \givenPrm) \projTan (\xmin))$.
\end{remark}

\subsection{Fixed-Point Automatic Differentiation} \label{ssec:PDM:FPAD}

We have already seen that FPAD combines the nicer properties of ID and AD and is a useful practical tool for estimating the derivative of the solution mapping. We can similarly apply FPAD on Algorithm~\ref{alg:APG} by modifying the update rules for AD as dictated in Section~\ref{sssec:FPI:FPAD}. In this section, we only inspect the forward mode FPAD; the analysis for reverse mode is analogous. For any $\u\in\spPrm$, $\alpha\in[\sslow, \ssup]$, $\beta\in[0, 1]$ and an approximate fixed-point $\x\in\manif$, we perform the following updates for forward mode FPAD of Algorithm~\ref{alg:APG}
\begin{equation} \label{itr:FPI:FPFAD:APG}
    \zkph = D_{\z} \apg (\x, \x, \u) \zkh + D_{\u} \apg (\x, \x, \u) \ud \,.
\end{equation}
for $\ud\in\spPrm$ and $k\geq0$. FPAD resolves the memory overhead problem while still retaining the convergence rate (see Remark~\ref{rem:APG:FAD:Lin:Conv}) and the ease-of-implementation of AD. The following result provides superiority of FPAD of Algorithm~\ref{alg:APG} over its AD counterpart.
\begin{theorem} \label{thm:APG:FPFAD}
    Let $f$ and $g$ satisfy Assumption~\ref{ass:CPSO} and $(\xmin, \givenPrm)\in\manif\times\setPrm$ be such that Assumptions~\ref{ass:RPD} and \ref{ass:ND} are satisfied. Let $\alpha\in[\sslow, \ssup]$ and $\beta\in[0, 1]$, such that $-1/(1 + 2\beta) < \lambda_{\min} (D_{\x} \pgd (\xmin, \givenPrm))$, $\u\in \neighbourhoodPRM$ and $\x$ be sufficiently closed to $\argmap (\u)$ where the (possibly reduced) neighbourhood $\neighbourhoodPRM$ and the mapping $\argmap$ are from Theorem~\ref{thm:IFT}. Then the sequence $( \zkh )_{k \in\N}$ generated by \eqref{itr:FPI:FPFAD:APG} converges linearly to $\apgPhi (\x, \x, \u)$ defined in \eqref{eq:APG:IFT:approx}, with rate $\rho (D_{\z} \apg(\x, \x, \u) \projTan (\x, \x))$.
\end{theorem}
\findProofIn{thm:APG:FPFAD}
\begin{remark}
    \begin{enumerate}[label=(\roman*)]
        \item We can also perform reverse mode FPAD by replacing $\apg (\x, \x, \u)$ with $(\x, \x)$ and $-\grad[\x] f (\x, \u)$ with $\bnu (\x, \u)$ in the expanded expressions for $D \apg (\x, \x, \u)$ (see Section~\ref{sssec:PDM:ID:NS}). In that case, the FPAD sequence converges to $\tilde\phi_{\alpha} (\x, \x, \u)$ defined in \eqref{eq:APG:mod:IFT:approx} and \eqref{eq:PGD:mod:IFT:approx}. Additionally, when $\manif$ is an affine manifold, the FPAD sequence converges to a quantity resembling $D\argmap (\u)$ in \eqref{eq:IFT} (see also Remark~\ref{rem:IDNS}\ref{itm:IDNS:affine}).
        \item We do not need to explicitly compute the projection term $\projTan (\x, \x)$ which appears in $\apgPhi (\x, \x, \u)$ thus making FPAD a practical method for computing the derivative of the solution mapping of \eqref{prob:min}.
        \item Similarly, the forward mode FPAD of PGD given by
        \begin{equation*}
            \xkph = D_{\x} \pgd (\x, \u) \xkh + D_{\u} \pgd (\x, \u) \ud \,,
        \end{equation*}
        converges linearly to $\pgdPhi (\x, \u)$ with rate $\rho (D_{\x} \pgd(\x, \u) \projTan (\x))$.
    \end{enumerate}
\end{remark}

% ***********************
% >>>>> EXPERIMENTS <<<<<
% ***********************
\section{Experiments} \label{sec:Exp}

We now put our results to test by applying them on some practical applications. We first show, in Section~\ref{ssec:Exp:CRD}, the experimental validation of convergence and convergence rates of AD and FPAD applied to \eqref{itr:APG} proved in Sections~\ref{ssec:PDM:AD} and \ref{ssec:PDM:FPAD} on Lasso \cite{Tib96} and Group Lasso \cite{YL06}. We verify with our experiments that the convergence rate of FPAD of \eqref{itr:APG} is at least as good as that of AD of \eqref{itr:APG}. In Section~\ref{ssec:Exp:ID}, we solve a bilevel optimization problem where we learn the parameters of the regularizer for image denoising. The parameter learning in this case is the upper-level problem and solving it requires to compute the derivative of the solution mapping of the lower-level problem, that is, the denoising problem. We use PyTorch \cite{PGM+19} for both experiments mainly to differentiate the forward and the backward step in \eqref{itr:APG}. In our applications, the backward step is expressed analytically and therefore, no Riemannian gradient and Hessian terms are evaluated as in \eqref{eq:FixMap:PGD}. For convergence rate demonstration, the unrolling is carried out without using the autograd because we want to generate the AD and FPAD sequences and visualize their behaviour. The remaining details are deferred to \refSupp{sec:Exp:supp}.

\subsection{Convergence Rate Demonstration} \label{ssec:Exp:CRD}
Given matrices $A\in\R^{M\times N}$ and $B\in\R^{M\times L}$, a vector $\b\in\R^M$ and some $\lambda > 0$, we consider following two problems for the demonstration of convergence rates of AD and FPAD of Algorithm~\ref{alg:APG}, namely the lasso problem
\begin{equation} \label{prob:lasso} \tag{L}
    \min_{\x\in\R^N} F_\mathrm{l}(\x, A, \b, \lambda), \quad F_\mathrm{l}(\x, A, \b, \lambda) \coloneqq \frac{1}{2} \norm[2]{A\x - \b}^2 + \lambda \norm[1]{\x} \,,
\end{equation}
and the group lasso problem
\begin{equation} \label{prob:group:lasso} \tag{GL}
    \min_{X\in\R^{N\times L}} F_{\mathrm{gl}}(X, A, B, \lambda), \quad F_{\mathrm{gl}}(X, A, B, \lambda) \coloneqq \frac{1}{2} \norm[2]{AX - B}^2 + \lambda \norm[2,1]{X} \,.
\end{equation}
In \eqref{prob:lasso}, $\norm[2]{\cdot}$ and $\norm[1]{\cdot}$ correspond to $\ell_2$ and $\ell_1$ norms defined on $\R^N$ while in \eqref{prob:group:lasso}, $\norm[2]{\cdot}$ is the Frobenius norm and $\norm[2,1]{\cdot}$ is the $\ell_{2,1}$ group norm defined on $\R^{N\times L}$ by $\norm[2,1]{X} \coloneqq \sum_{i} \sqrt{\sum_{j} x_{i,j}^2}$. In our experiments, we set $u\coloneqq (A, \b, \lambda)$ for \eqref{prob:lasso} and $u\coloneqq (A, B, \lambda)$ for \eqref{prob:group:lasso}.

\begin{figure}[ht]
    \centering
    \begin{subfigure}[b]{0.98\textwidth}
        \centering
        \includegraphics[width=0.98\textwidth]{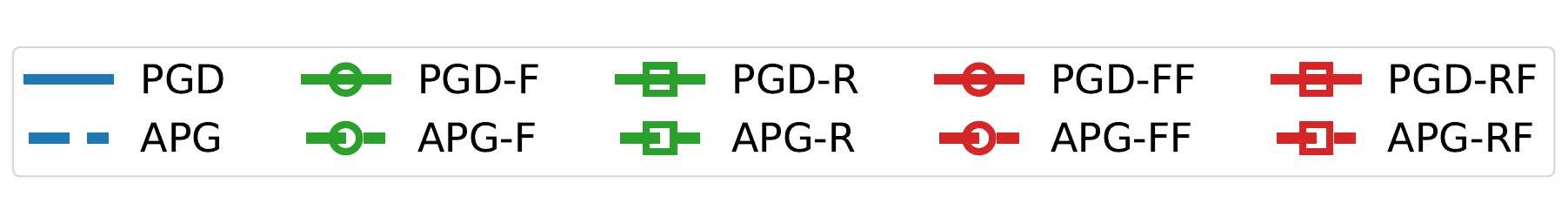}
        \label{sfig:cr:legs}
        \vspace*{-2.9ex}
    \end{subfigure}
    \begin{subfigure}[b]{0.475\textwidth}
        \centering
        \includegraphics[width=0.98\textwidth,height=0.8\textwidth]{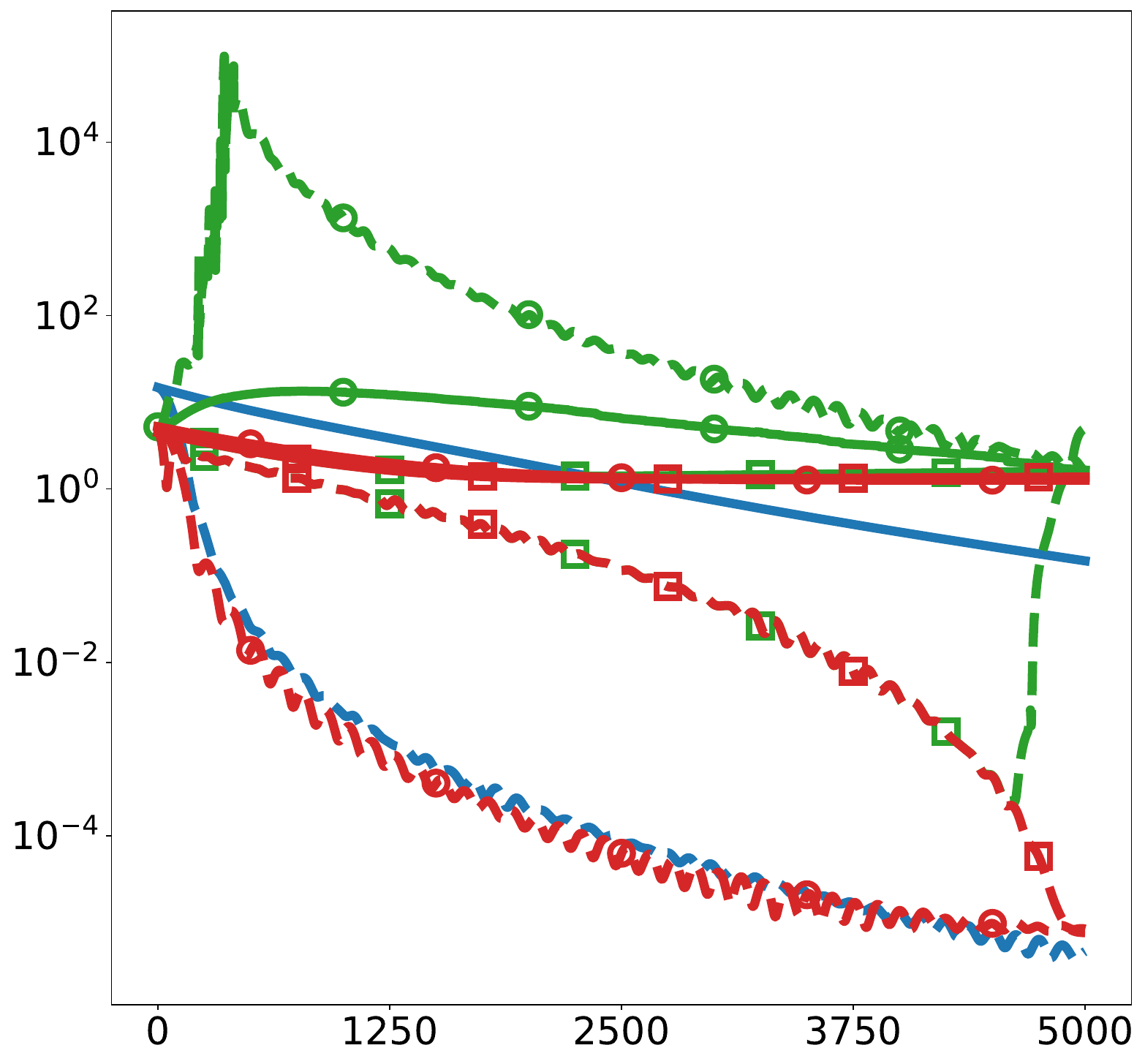}
        \label{sfig:cr:lasso}
    \end{subfigure}
    \begin{subfigure}[b]{0.475\textwidth}
        \centering
        \includegraphics[width=0.98\textwidth,height=0.8\textwidth]{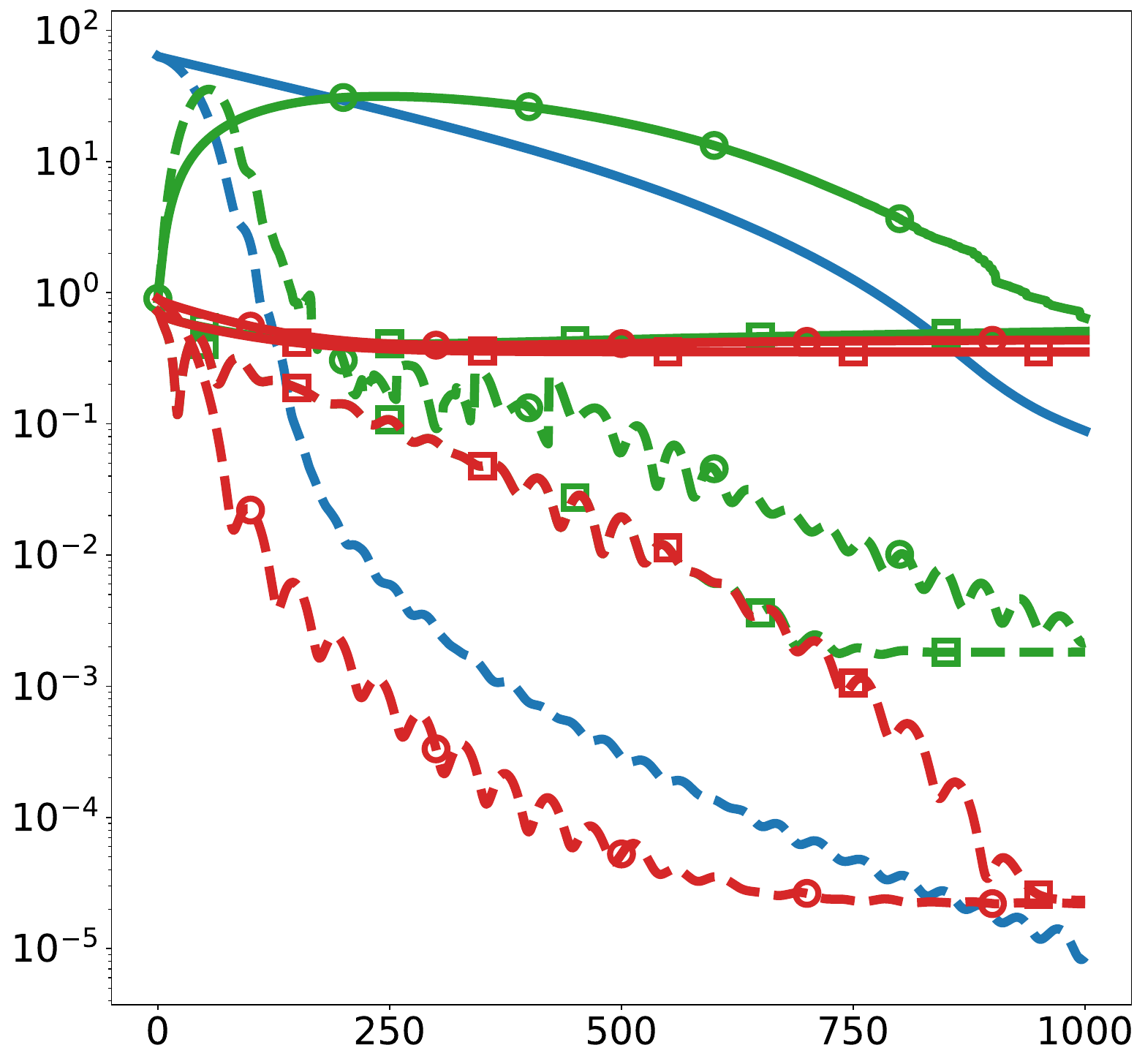}
        \label{sfig:cr:glasso}
    \end{subfigure}
    \caption{Error Plots for various sequences obtained from \eqref{prob:lasso} (left) and \eqref{prob:group:lasso} (right). The solid lines represent PGD and its derivative sequences while the dotted lines represent APG sequences. The lines without any markers correspond to the original iterates. The vertical lines mark the iteration at which the sequences enter the manifold (solid for PGD and dashed for APG). The lines with circular markers represent forward mode sequences while those with square markers represent Reverse Mode sequences. Markers for FPAD sequences are larger than their AD counterparts. The FPAD sequences are more stable and show faster convergence to the true derivative.}
    \label{fig:conv:rates}
\end{figure}

In Figure~\ref{fig:conv:rates}, we plot all the error sequences. For computing the error sequences, we run \eqref{itr:APG} for a long time to compute a very good estimate of $\xmin$ and compute $\xmind=D\argmap (\u)\ud$ and $\ub=\xminb D\argmap (\u)$ where $D\argmap (\u)$ is defined in \eqref{eq:IFT}. We use Conjugate Gradient Method \cite{HS52} to solve the linear systems for computing $\xmind$ and $\ub$. The left figure shows the sequences computed for \eqref{prob:lasso} while the right one shows those obtained from \eqref{prob:group:lasso}. In both plots, we note that APG algorithm (dashed blue lines) converges faster than PGD (solid blue lines). This acceleration effect is translated to the FPAD sequences; the sequences obtained by applying FPAD on APG (dashed red lines with circle for forward mode and square markers for Reverse Mode) converge more quickly as compared to those generated by applying FPAD on PGD (see corresponding solid red lines). For both problems, PGD requires more time to converge and therefore does not provide satisfying results for any derivative sequence (all the solid lines move slowly towards $0$). For \eqref{prob:group:lasso}, the forward and reverse mode AD of APG (dashed green lines) yield good results but after $K$ iterations, we observe that FPAD gives a better estimate (dashed red lines). For \eqref{prob:lasso}, the behaviour of forward mode AD on APG (dashed green line with circular marker) in the beginning and that of reverse mode AD on APG (dashed green line with square marker) becomes bad in the end is quite erratic. One possibility for this situation might be that \eqref{itr:APG} may not behave well in the initial iterations. Note that the derivative is a local quantity and therefore convergence of the generated sequence of derivatives may be hampered in initial iterates that are far away from this local neighbourhood.

% Conclusion
In the above experiment we saw that FPAD iterates behave as well as the original iterates and the convergence for the two types of sequences is very similar. We noted that the FPAD sequences also benefit from inertia just like the optimization algorithm they are being applied to. We further found out that AD applied in a naive way may not provide a good result in the end and the algorithm must be initially run for some time before starting to construct the computational graph. These findings together with the enormous memory-gain boost FPAD as a practical approach for computing derivative in such a setting.

\subsection{Learning Regularizers for Image Denoising} \label{ssec:Exp:ID}

We now demonstrate the effectiveness of FPAD by using it to solve a simple yet practical bilevel optimization problem. Let $\spImg\coloneqq\R^{N_x\times N_y\times N_c}$ and $\spFeat\coloneqq\R^{N_x\times N_y\times N_c\times N_f}$, we consider the following optimization problem
\begin{equation} \label{prob:ID}
    \min_{\recovered\in\spImg} \frac{1}{2}\norm[2]{\recovered - \noisy}^2 + \lambda \norm[\mathrm{r}]{\imgLinOp \recovered} \,,
\end{equation}
where $\lambda > 0$, $\imgLinOp\in\spLin (\spImg, \spFeat$) and $\norm[\mathrm{r}]{\cdot}$ is an appropriate regularization norm. The best way to solve the above problem is by applying APG to its Fenchel-Rockafellar dual
\[
    \min_{\edges\in\spFeat} \frac{1}{2}\norm[2]{\noisy - \imgLinOp^*\edges}^2 \quad \st \quad \norm[\mathrm{r},*]{\edges} \leq \lambda \,.
\]
From \cite[Proposition~15.24(x)]{BC11}, strong duality holds and the dual solution primal solution can be computed from the dual solution by $\recovered^* = \noisy - \imgLinOp^*\edges^*$. For the dual problem, Assumption~\ref{ass:CPSO} is satisfied. Assumption~\ref{ass:ND} is satisfied almost everywhere \cite[Theorem~3]{VDP+17}. However, Assumptions~\ref{ass:RPD} can be violated when $\dim (\ker \imgLinOp)$ is large relative to $\dim\spImg$. Here $\norm[\mathrm{r},*]{\cdot}$ is the dual norm of $\norm[\mathrm{r}]{\cdot}$. When $N_c\in\set{1,3}$, $N_f=2$, $\imgLinOp$ is a convolution operator computing the forward $x$ and $y$-derivative of its input image with Neumann boundary condition and $\norm[\mathrm{r}]{\cdot}$ is the $\ell_1$ norm or the $\ell_{2,1}$ group norm defined by $\norm[2,1]{\p} \coloneqq \sum_{i,j} \sqrt{\sum_{k,r} p_{i,j,k,r}^2}$, \eqref{prob:ID} corresponds to total variation image denoising \cite{ROF92} which aims to produce a denoised image $\recovered^*$ from a noisy image $\noisy$. We instead express $\imgLinOp$ as a weighted sum of any given basis filters and learn the weights $\Theta$ by solving the following bilevel learning problem.
\begin{equation} \tag{BL} \label{prob:BL}
\begin{aligned}
\min_{\Theta} \quad & \sum_{i=1}^{M_{\textnormal{tr}}} J_i (\Theta) \,, \quad J_i (\Theta) \coloneqq \frac{1}{2}\norm[2]{\recovered[i]^* (\Theta) - \ground[i]}^2 \\
\st \quad & \recovered[i]^* (\Theta) = \noisy - \imgLinOp^*\edges[i]^* (\Theta) \,, \\
& \edges[i]^* (\Theta) \in \argmin_{\norm[\mathrm{r},*]{\edges} \leq \lambda} \frac{1}{2}\norm[2]{\noisy - \imgLinOp^*\edges}^2 \,,\quad i = 1,\ldots,M_{\textnormal{tr}} \,.
\end{aligned}
\end{equation}
We solve \eqref{prob:BL} by using stochastic algorithms \cite{BCN18}, for example, Stochastic Gradient method. For back-propagating through the lower level problem, we make use of reverse mode FPAD of Algorithm~\ref{alg:APG}.

We train our model for $30$ epochs and use a training set of size $5$. The average PSNR value of training set improves from $18.38$ of noisy images to $21.53$ after only one epoch suggesting that the FPAD scheme is providing the right derivatives. The visual depiction of the learning process is provided in Figures~\ref{fig:ID:evol:supp} and \ref{fig:ID:test:supp} in the appendix. This example clearly shows the effectiveness of FPAD because if we had used AD instead, the memory overhead would have been of the order of $K N_x N_y N_c = 3.75\times10^{6}$ which equates to $28.6$ MiB as opposed to $N_x N_y N_c = 7500$ or $58.6$ KiB of FPAD.

\section{Conclusions}
\label{sec:Conc}

In this paper, we studied the problem of differentiating the solution mapping of a structured parametric optimization problem. We showed that when the objective of this problem is partly smooth and some other assumptions are satisfied, the classical derivative of the solution mapping can be estimated by Implicit and Automatic Differentiation of proximal splitting algorithms like PGD and APG. We showed that the memory overhead of the reverse mode AD can be overcome by using the Fixed-Point Automatic Differentiation technique. We also showed that in terms of convergence speed, FPAD outperforms AD both in theory and in practice. We also showed the working of FPAD by solving a bilevel optimization problem to learn the regularizers for image denoising.

% *******************
% >>>>> APPENDIX <<<<<
% *******************

\appendix

\section{Automatic Differentiation} \label{sec:AD:supp}
As the name suggests, Automatic or Algorithmic differentiation is a way of computing the derivatives of a function, given as a computer program, automatically. The key is the classic chain rule and a representation of a function as composition of a finite number of elementary functions like polynomials, logarithms, exponential and trigonometric functions etc. From classical Calculus, the derivatives of these elementary functions are known and therefore are combined via chain rule to obtain the derivative of the given function. There are two main modes of AD, namely the forward and the reverse mode which we briefly recall here. For more on AD and its applications, the reader is requested to look into \cite{GW08,BPR+18} and the references therein.

Let $\mathcal X$, $\mathcal Y$, and $\mathcal Z$ be Euclidean spaces. Consider a function $\map{h}{\mathcal X}{\mathcal Z}$ given by
\[
    h \coloneqq (h_3 h_2) \circ h_1 \,,\quad h_1 \textnormal{ composed with the product } h_3 h_2 \,
\]
where $\map{h_1}{\mathcal X}{\mathcal Y}$, $\map{h_2}{\mathcal Y}{\mathcal Z}$, and $\map{h_3}{\mathcal Y}{\R}$ are functions with known derivatives. We compute the derivative of $h$ at some $\x \in \mathcal X$ via the two modes of AD in Sections~\ref{ssec:AD:FM:supp} and \ref{ssec:AD:RM:supp}. We first use some intermediate variables to break down our example in the following manner
\begin{equation}
    \begin{aligned}
        \y &= h_1 (\x) \,, \\
        \w &= h_2 (\y) \,, \; a = h_3 (\y) \,, \\
        \z &= a\w \,.
    \end{aligned}
\end{equation}
\begin{figure}[ht]
	\centering
	\tikzset{font={\fontsize{10pt}{12}\selectfont}}
	\begin{tikzpicture}[->, semithick, node distance=3.5cm, >=stealth', auto]
	\tikzstyle{place} = [circle, thick, draw=black!75, fill=black!10, minimum size=6mm]
	\begin{scope}
	
	\node [place] (x) [label=above:$\x$] {};
	\node [place] (y) [right of=x, label=above:$\y$] {};
	\node [place] (w) [right of=y, yshift=+1.6cm, label=above:$\w$] {};
	\node [place] (a) [right of=y, yshift=-1.6cm, label=above:$a$] {};
	\node [place] (z) [right of=w, yshift=-1.6cm, label=left:$\times$, label=above:$\z$] {};
	
	\path[every node/.style={sloped,anchor=south,auto=false}]
	(x) edge node [above, pos=0.5] {$Dh_1 (\x)$} (y)
	(y) edge node [above, pos=0.5] {$Dh_2 (\y)$} (w)
	(y) edge node [above, pos=0.5] {$Dh_3 (\y)$} (a)
	(w) edge node [above, pos=0.5] {$a$} (z)
	(a) edge node [above, pos=0.5] {$\w$} (z);
	
	\end{scope}
	\end{tikzpicture}
	\caption{Computational Graph of $ h(\bm{x}) $.}
	\label{fig:CG:supp}
\end{figure}
This allows us to construct the computational graph for $h$ as shown in Figure~\ref{fig:CG:supp}. A node represents a variable while a directed edge shows the dependence of one variable on another (represented by the nodes it connects) through some elementary function or operation. In the figure, the nodes are labelled by the variables while the edges are denoted by the derivative (Jacobian) of the function they represent.

\subsection{Forward Mode AD} \label{ssec:AD:FM:supp}
The forward mode derivative $\zd$ of $\z$ is represented by placing a dot on it and is computed by starting with some $\xd\in\mathcal X$ and computing the directional derivative of $\z$ with respect to $\x$ along $\xd$, that is,
\begin{equation}
    \begin{aligned}
        \yd &= Dh_1 (\x) \xd \,, \\
        \dot\w &= h_2 (\y) \yd\,, \; \dot a = Dh_3 (\y) \yd \,, \\
        \zd &= a\dot\w + \dot a\w \,.
    \end{aligned}
\end{equation}
\begin{figure}[ht]
	\centering
	\tikzset{font={\fontsize{10pt}{12}\selectfont}}
	\begin{tikzpicture}[->, semithick, node distance=3.5cm, >=stealth', auto]
	\tikzstyle{place} = [circle, thick, draw=black!75, fill=black!10, minimum size=6mm]
	\begin{scope}
	
	\node [place] (x) [label=above:$\x$] {};
	\node [place] (y) [right of=x, label=above:$\y$] {};
	\node [place] (w) [right of=y, yshift=+1.6cm, label=above:$\w$] {};
	\node [place] (a) [right of=y, yshift=-1.6cm, label=above:$a$] {};
	\node [place] (z) [right of=w, yshift=-1.6cm, label=left:$\times$, label=above:$\z$] {};
	\node [dashed, place] (t) [left of=x, label=above:$t$] {};
	
	\path[every node/.style={sloped,anchor=south,auto=false}]
	(t) edge [dashed] node [above, pos=0.5] {$\xd=D\xi(t)$} (x)
	(x) edge node [above, pos=0.5] {$Dh_1 (\x)$} (y)
	(y) edge node [above, pos=0.5] {$Dh_2 (\y)$} (w)
	(y) edge node [above, pos=0.5] {$Dh_3 (\y)$} (a)
	(w) edge node [above, pos=0.5] {$a$} (z)
	(a) edge node [above, pos=0.5] {$\w$} (z);
	
	\end{scope}
	\end{tikzpicture}
	\caption{Depiction of Forward Mode AD of $ h $.}
	\label{fig:FAD:supp}
\end{figure}
We can also motivate the forward mode AD in the following manner. Assume that $\x$ is obtained by applying a differentiable function $\map{\xi}{\R}{\mathcal X}$ on some scalar $t\in\R$, that is, $\x=\xi(t)$. The forward mode derivative $\xd$ is then interpreted as the derivative of $\xi$ at $t$, that is, $\xd = D\xi(t)$ (see Figure~\ref{fig:FAD:supp}). In other words $\xd$ is the derivative of the variable $\x$ with respect to $t$, at $t$, that is, $\xd=\dydx{\x}{t}(t)$. In fact, all the forward mode derivatives $\yd$, $\dot\w$, $\dot a$ and $\zd$ are interpreted as the derivatives $\dydx{\y}{t}(t)$, $\dydx{\w}{t}(t)$, $\dydx{a}{t}(t)$ and $\dydx{\z}{t}(t)$. The intermediate variables $\y$, $\w$ and $a$ need not be stored since they can be computed alongside the derivatives $\yd$, $\dot\w$, $\dot a$ and $\zd$. However, even with no memory overhead, forward mode AD becomes impractical when the dimension of $\mathcal X$ is significantly larger than that of $\mathcal Z$, since, for example, computing the derivative $D h(\x)$ requires $\dim (\mathcal X)$ runs of the forward mode, once for each directional derivative.

\subsection{Reverse Mode AD} \label{ssec:AD:RM:supp}
The reverse mode derivatives are represented by bars on the variables. In this mode, we start at the output and move towards the input. That is, given $\bar\z\in\mathcal Z^*$, we compute $\bar\x\in\mathcal X_*$ by
\begin{equation} \label{eq:ADeg:RM:supp}
    \begin{aligned}
        \bar a &= \bar\z \w\,, \; \bar\w = a \bar\z \,, \\
        \bar\y &= \bar a Dh_3 (\y) + \bar\w Dh_2 (\y) \,, \\
        \bar\x &= \bar\y Dh_1 (\x) \,.
    \end{aligned}
\end{equation}
We can motivate the reverse mode analogously to the forward mode. Assume that $\z$ is fed into a scalar-valued function $\map{\zeta}{\mathcal Z}{\R}$ to obtain $s = \zeta(\z)$. The reverse mode derivative is then treated as the derivative of $\zeta$ at $\z$, that is, $\bar\z=D\zeta(\z)=\dydx{s}{\z}(\z)$ (Figure~\ref{fig:RAD:supp}). Similarly, $\bar\x$, $\bar a$, $\bar\w$ and $\bar\y$ are interpreted as the derivatives $\dydx{s}{\x}(\x)$, $\dydx{s}{a}(a)$, $\dydx{s}{\w}(\w)$ and $\dydx{s}{\y}(\y)$, respectively. This scheme is referred to as back-propagation in the Machine Learning community \cite{RHW86} and is useful when the dimension of the output space is very small, for example, a (scalar-valued) loss function. Its computation cost is independent of the size of the input thus making it the preferred mode in Deep Learning. However, since we are moving backward when computing the derivatives, we cannot compute the intermediate variables $\y$, $\w$ and $a$ in parallel anymore. We must store them beforehand to efficiently use this mode which causes a memory overhead, as can be seen by the dependence of each step in \eqref{eq:ADeg:RM:supp}.
\begin{figure}[ht]
	\centering
	\tikzset{font={\fontsize{10pt}{12}\selectfont}}
	\begin{tikzpicture}[->, semithick, node distance=3.5cm, >=stealth', auto]
	\tikzstyle{place} = [circle, thick, draw=black!75, fill=black!10, minimum size=6mm]
	\begin{scope}
	
	\node [place] (x) [label=above:$\x$] {};
	\node [place] (y) [right of=x, label=above:$\y$] {};
	\node [place] (w) [right of=y, yshift=+1.6cm, label=above:$\w$] {};
	\node [place] (a) [right of=y, yshift=-1.6cm, label=above:$a$] {};
	\node [place] (z) [right of=w, yshift=-1.6cm, label=left:$\times$, label=above:$\z$] {};
    \node [dashed, place] (s) [right of=z, label=above:$s$] {};
	
	\path[every node/.style={sloped,anchor=south,auto=false}]
	(z) edge [dashed] node [above, pos=0.5] {$\bar\z=D\zeta(\z)$} (s)
	(x) edge node [above, pos=0.5] {$Dh_1 (\x)$} (y)
	(y) edge node [above, pos=0.5] {$Dh_2 (\y)$} (w)
	(y) edge node [above, pos=0.5] {$Dh_3 (\y)$} (a)
	(w) edge node [above, pos=0.5] {$a$} (z)
	(a) edge node [above, pos=0.5] {$\w$} (z);
	
	\end{scope}
	\end{tikzpicture}
	\caption{Depiction of Reverse Mode AD of $ h $.}
	\label{fig:RAD:supp}
\end{figure}

\section{Computational Aspects}
In this section, we briefly look at the computational aspects when differentiating fixed-point iterations. In particular, we demonstrate the ease in implementing FPAD through a standard autograd library. We also compare the computation time and memory overhead of various techniques to compute the derivative of the solution mapping.

\subsection{Implementation of FPAD} \label{ssec:FPAD:Imp:supp}
FPAD is not only time and memory efficient but also very easy to implement. In fact, this can be achieved by using an automatic differentiation library \cite[Algorithm~3.1]{Chr94}. Below, we provide a code snippet for computing the derivative of the fixed-point of $\R^2\ni(x, u)\mapsto \fixMap (x, u) \coloneqq (x + u/x) / 2$ in PyTorch \cite{PGM+19}. These fixed-point iterations are the Babylonian method for computing $\sqrt u$. As suggested in \cite[Algorithm~3.1]{Chr94}, we first run the fixed-point iterations $K=10$ times to obtain a good estimate of the fixed-point. The computational graph is not generated during this process. Afterwards, we generate the computational graph for a single iteration, that is, \verb|y = iteration_map(x, u)| and apply the backward function on the output \verb|y| repeatedly. The default value of the property \verb|retain_graph=False| destroys the graph after a single backward pass and therefore must be set to \verb|True|. Calling the backward function repeatedly will accumulate the gradients of $y$ with respect to both $x$ (not desired!) and $u$ (desired!). Therefore we reinitialize the gradient of $y$ with respect $x$ at every iteration. At every call of \verb|y.backward(y_b, retain_graph=True)|, we are basically performing \eqref{itr:FPI:FPRAD} where \verb|x.grad|, \verb|u.grad| and \verb|y_b| are $\xkt$, $\unt$ and $\xkpt$ respectively.
\begin{lstlisting}[frame=single,language=Python,keywordstyle=\color{blue}]
import torch

def iteration_map(x, u):
    return 0.5 * (x + u / x)

# Setting Parameters.
u = torch.tensor(10.)
K = 10

# Computing the fixed point x^*.
x = u.clone()
for k in range(K):
    x = iteration_map(x, u)

# Computing dx^*/du by using FPAD.
y_b = torch.tensor(1.)
x = x.clone().detach().requires_grad_()
u = u.clone().detach().requires_grad_()
y = iteration_map(x, u)
for k in range(K):
    y.backward(y_b, retain_graph=True)
    y_b, x.grad = x.grad, None

print(float(x - torch.sqrt(u)))             # Value Error
print(float(u.grad - 0.5/torch.sqrt(u)))    # Derivative Error
\end{lstlisting}

\subsection{Time and Memory Complexity} \label{ssec:T&MC:supp}

We now compare the time and space complexity of computing the derivative of the fixed-point $\xmin$ of \eqref{itr:FPI} with respect to $\u$, that is, $D_{\u} \xmin$. Since, in the end, we want to minimize a function $\map{\ell}{\spPrm}{\R}$ defined through some mapping $\map{L}{\spVar}{\R}$ by $l(\u) = L(\xmin)$ in the context of bilevel optimization (see Sections~\ref{ssec:intro:apps} and \ref{ssec:Exp:ID}), we consider the computation complexity for estimating $\grad l (\u) = (D_{\u} \xmin) ^T \grad L(\xmin)$. Table~\ref{tab:comp:meths} summarizes the space and time complexity for various methods. From our discussion earlier, we note that the forward mode AD does not have a memory overhead as compared to the reverse mode AD while the computation time taken by the former is $P$ times that of the latter \cite[Chapter~4]{GW08}. On the other hand, while both the forward and the reverse mode FPAD take the same time as that of their AD counterparts, the reverse mode FPAD has no memory overhead, just like the forward mode AD and the forward mode FPAD. For using LU or Cholesky Factorization for Implicit Differentiation, we are required to compute and store the matrix representation of $D_{\x} \fixMap(\xmin, \u)$. However, since we are working with partly smooth functions -- where $\xmin$ lies on an $m$-dimensional manifold $\manif$ -- the storage needed is $\O (n^2)$ instead of $\O(N^2)$ with $m\leq n\leq N$. However, if we use conjugate gradient method to solve the linear system, both the time and memory complexity are reduced.

\begin{table}[]
    \centering
    \begin{tabular}{|c|c|c|}
        \hline
        \textbf{Method} & \textbf{Space Complexity} & \textbf{Time Complexity} \\
        \hline
        Forward AD & $\O(N)$ & $\O(KPt_{\fixMap} + Pt_L)$ \\
        \hline
        Reverse AD & $\O(KN)$ & $\O(Kt_{\fixMap} + t_L)$ \\
        \hline
        Forward FPAD & $\O(N)$ & $\O(KPt_{\fixMap} + Pt_L)$ \\
        \hline
        Reverse FPAD & $\O(N)$ & $\O(Kt_{\fixMap} + t_L)$ \\
        \hline
        ID (Factorization) & $\O(n^2)$ & $\O(n^3 + Nt_{\fixMap} + t_L)$ \\
        \hline
        ID (Conjugate Gradient) & $\O(N)$ & $\O(N t_{\fixMap} + t_L)$ \\
        \hline
    \end{tabular}
    \vspace{1em}
    \caption{Comparison of time and space complexity of various methods. $t_{\fixMap}$ and $t_L$ are the time complexities for performing one computation of mappings $\fixMap$ and $L$, respectively. In the last two columns, $m \leq n \leq N$ is the size of the (possibly) reduced system of equations arising from the fact that $\xmin$ lies on a manifold $\manif$ of dimension $m$.}
    \label{tab:comp:meths}
\end{table}

\section{Proofs of Section~\ref{sec:intro}}

\appSubSect{Lemma}{lem:subdiff:cartesian}
\vspace{2ex}
\begin{proof}
    The result is a straightforward application of \cite[Exercise~8.8, Corollary~10.11]{RW98}. ``$\subset$'' is trivial because $(\y, \v) \in \partial g (\x, \u)$ implies that $\y \in \partial_{\x} g (\x, \u)$ and $\v \in \partial_{\u} g (\x, \u) = \grad[\u] g (\x, \u)$. ``$\supset$'' follows because for any $\y\in\partial_{\x} g (\x, \u)$, we must have $(\y, \grad[\u] g (\x, \u))\in \partial g (\x, \u)$.
\end{proof}

\appSubSect{Lemma}{lem:spec_rad:inv_op}
\vspace{2ex}
\begin{proof}
    For any $\epsilon\in (0, 1-\rho)$, there exists a norm $\norm{\cdot}$ on $\spLin (\spVar, \spVar)$ consistent with a vector norm $\map{\norm{\cdot}}{\spVar}{\R}$ such that $\norm{B}\leq \rho+\epsilon$, where we abbreviate $\rho\coloneqq\rho (B)$. For any $\v\in\spVar$, $(I-B)\v = 0$ or equivalently $\v = B\v$ implies $\norm{\v} \leq \norm{B}\norm{\v} \leq (\rho+\epsilon) \norm{\v} $. But since $\rho+\epsilon < 1$, we conclude that $\v = 0$ which in turn implies that $\opid - B$ is invertible.
\end{proof}

\appSubSect{Theorem}{thm:LGFPI:conv}
\vspace{2ex}
\begin{proof}
    For any $\delta\in (0, 1-\rho)$, let $\norm{\cdot}$ on $\spLin (\spVar, \spVar)$ be a consistent norm with the vector norm $\map{\norm{\cdot}}{\spVar}{\R}$ such that $\norm{B}\leq \rho+\delta/2$. Following \cite[Proposition~2.7]{Rii20}, we set $\zk\coloneqq\xk - \x^*$ with $\x^* \coloneqq (I-B)^{-1}\b$, $\yk\coloneqq (B_k - B)\x^* + \bk - \b$ and $\wk\coloneqq (B_k - B)\zk $ which leads to
    \begin{equation*}
        \begin{aligned}
            \zkp &= (B_k \xk + \bk) - (B\x^* + \b) \\
            &= B \zk + \wk + \yk \,.
        \end{aligned}
    \end{equation*}
    Since $\wk = \o(\zk)$, we can find $K\in\N$ such that $\norm\wk / \norm\zk \leq \delta/2 $ for all $k\geq K$.
    By using the matrix and vector norms defined above, we obtain
    \begin{equation*}
        \begin{aligned}
            \norm{\zkp} &\leq (\rho + \delta/2) \norm\zk + \norm\wk + \norm\yk \\
            &\leq (\rho + \delta) \norm\zk + \norm\yk \,.
        \end{aligned}
    \end{equation*}
    Recursive expansion of the above expression yields
    \begin{equation} \label{ineq:LGFPI:zk:yk}
        \norm{\zkp} \leq (\rho + \delta)^{k-K+1} \norm{\zK} + \sum_{i=K}^{k} (\rho + \delta)^{k-i} \norm{\y\iter{i}} \,.
    \end{equation}
    Using the fact that $\yk\to 0$, for any $\gamma>0$, we can make $K$ large enough to obtain $\norm{\y\iter{i}} \leq (1-\rho-\delta) \gamma/2$ for all $i\geq K$. This allows for the following reduction
    \begin{equation*}
        \begin{aligned}
            \norm{\zkp} &\leq (\rho + \delta)^{k-K+1} \norm{\zK} + \sum_{i=K}^{k} (\rho + \delta)^{k-i} (1-\rho-\delta) \frac{\gamma}{2} \\
            &\leq (\rho + \delta)^{k-K+1} \norm{\zK} + (1-\rho-\delta) \frac{\gamma}{2} \sum_{i=0}^{\infty} (\rho + \delta)^i \\
            &= (\rho + \delta)^{k-K+1} \norm{\zK} + \frac{\gamma}{2} \,.
        \end{aligned}
    \end{equation*}
    Also, the quantity $(\rho + \delta)^{k-K+1} \norm{\zK}$ eventually becomes smaller than $\gamma/2$ as $k\geq K$ grows. In particular, there exists $N\geq K$ such that $(\rho + \delta)^{k-K+1} \norm{\zK} < \gamma/2$ and hence $\norm{\zkp} < \gamma$, for all $k\geq N$. This concludes the proof for convergence because $\gamma > 0$ was arbitrarily chosen.
    
    For the rate of convergence, we set $c_1(\delta)\coloneqq c(\delta) (\norm{\x^*} + 1)$ and use the bounds in \eqref{ineq:LGFPI:B:b:conv_rates} to expand \eqref{ineq:LGFPI:zk:yk} and obtain
    \begin{equation*}
        \begin{aligned}
            \norm{\zkp} &\leq (\rho + \delta)^{k-K+1} \norm{\zK} + \sum_{i=K}^{k} (\rho + \delta)^{k-i} c_1(\delta) (\rho + \delta)^{i-K} \\
            &\leq c_1(\delta) (k-K+1) (\rho + \delta)^{k-K} + (\rho + \delta)^{k-K+1} \norm{\zK} \\
            &\leq \frac{c_1(\delta)}{\rho + \delta} (k-K+1) (\rho + \delta)^{k-K+1} + \norm{\zK} (\rho + \delta)^{k-K+1} \,.
        \end{aligned}
    \end{equation*}
    Because all the norms on $\spVar$ are equivalent, the same bounds hold for the norm  induced from the inner product with possibly different constants which concludes our proof.
\end{proof}

\appSubSect{Corollary}{cor:FPI:FAD}
\vspace{2ex}
\begin{proof}
    The convergence of $\xkd (\givenPrm)$ follows from Assumption~\ref{ass:C1:contraction}, Theorem~\ref{thm:LGFPI:conv} and the fact that $D \fixMap (\xk (\givenPrm), \givenPrm) \to D \fixMap (\xmin, \givenPrm)$ . For linear convergence, we observe that the local Lipschitz continuity condition allows the linear convergence of both $D_{\x} \fixMap (\xk (\givenPrm), \givenPrm)$ and $D_{\u} \fixMap (\xk (\givenPrm), \givenPrm)$ with the same rate as that of $\xk$ because there exist $L > 0$ and $K\in\N$ such that for all $k\geq K$, we have $\norm{D\fixMap(\xk (\givenPrm), \givenPrm) - D\fixMap(\xmin, \givenPrm)} \leq L \norm{\xk (\givenPrm) - \xmin}$. The expression for the rate is furnished by combining the results of Theorem~\ref{thm:FPI} and Theorem~\ref{thm:LGFPI:conv}.
\end{proof}

\appSubSect{Corollary}{cor:FPI:FAD:Global}
\vspace{2ex}
\begin{proof}
    When $\fixMap$ is $C^1$-smooth on the whole space, we do not require $\xz$ to be close enough to $\xmin$ for the iterations in \eqref{itr:FPI} to be well-defined. Also, because the local Lipschitz continuity of $D\fixMap$ ensures that $\norm{D\fixMap(\xk (\givenPrm), \givenPrm) - D\fixMap(\xmin, \givenPrm)} \leq L \norm{\xk (\givenPrm) - \xmin}$ holds \textit{eventually} for all $k$, the convergence rate of $D_{\x} \fixMap (\xk (\givenPrm), \givenPrm)$ and $D_{\u} \fixMap (\xk (\givenPrm), \givenPrm)$ is still the same as in Corollary~\ref{cor:FPI:FAD}.
\end{proof}

\appSubSect{Corollary}{cor:FPI:FAD:Late}
\vspace{2ex}
\begin{proof}
    Once again, the sequence $D \fixMap (\x\iter{k+N} (\givenPrm), \givenPrm)$ converges to $D \fixMap (\xmin, \givenPrm)$ with rate same as that of $\xk (\givenPrm)$ due to local Lipschitz continuity which allows us to apply Theorem~\ref{thm:LGFPI:conv} under Assumption~\ref{ass:C1:contraction}.
\end{proof}

\appSubSect{Lemma}{lem:spec_rad:neu_ser}
\vspace{2ex}
\begin{proof}
    The recursive step in \eqref{itr:Neu:fwd} is easily verified by
    \begin{equation*}
        X_{r+1} = \sum_{i=0}^{r} B^i C = C + \sum_{i=1}^{r} B^i C = C + B \sum_{i=0}^{r-1} B^i C = C + B X_r \,,
    \end{equation*}
    For \eqref{itr:Neu:bwd}, we define $Y_r\coloneqq B^r$ and obtain
    \begin{equation*}
        X_{r+1} = \sum_{i=0}^{r} B^i C = B^{r} C + \sum_{i=0}^{r-1} B^i C = B^{r} C + X_r = Y_{r} C + X_r \,,
    \end{equation*}
    with $Y_0 \coloneqq I$ and $Y_{r+1} = Y_{r} B$. We get the convergence and the convergence rate, for any $X_0\in\spLin(\spVar, \spVar)$, by writing $X_{r+1} - X_* = B (X_r - X_*)$ from \eqref{itr:Neu:fwd} and choosing a norm $\norm{\cdot}$ on $\spLin (\spVar, \spVar)$ for any $\delta \in (0, 1-\rho)$.
\end{proof}

\section{Proofs of Section~\ref{sec:PS}}
\appSubSect{Theorem}{thm:IFT}
\vspace{2ex}
\begin{proof}
    Assumptions~\ref{ass:RPD}\ref{itm:RPD-i} and \ref{ass:ND} ensure that $\xmin$ is a strong critical point of $F(\cdot, \givenPrm)$ relative to $\manif$ (see Remark~\ref{rem:ND}\ref{itm:strong:crit}). The transversal embedding condition is also satisfied for $\manif\times\setPrm$ \cite[Assumption~5.1]{Lew02}. Thus \ref{itm:IFT:i} and \ref{itm:IFT:ii} follow from \cite[Theorem~5.7]{Lew02}. Because of Fermat's rule, we have $\grad[\manif] F(\argmap (\u), \u) = 0$ for every $\u\in\setPrm$. Differentiating this equation with respect to $\u$ by using the calculus from differential geometry \cite[Exercise~3.39, Example~5.19]{Bou23} and noting that \eqref{eq:RPD-i} is satisfied by $F$ at $(\argmap (\u), \u)$, we arrive at \ref{itm:IFT:iii}.
\end{proof}
    
\section{Proofs of Section~\ref{sec:APG}}
\appSubSect{Lemma}{lem:prox:diff}
\vspace{2ex}
\begin{proof}
    We identify the objective in \eqref{eq:prox} by mapping $\map{h}{\spVar\times\spVar\times\spPrm\times(0, 2/L)}{\eR}$, that is, $h(\y, \w, \u, \alpha) \coloneqq \alpha g(\y, \u) + \frac{1}{2} \norm{\y - \w}^2$ with parameters $\w$ and $\u$. To prove this result, we fix $\alpha_*\in (0, 2/L)$ and verify the assumptions similar to Assumptions~\ref{ass:CPSO}, \ref{ass:RPD}\ref{itm:RPD-i} and \ref{ass:ND} for $h$ at $(\bwd (\w^*, \givenPrm), \w^*, \givenPrm, \alpha_*)$ where $\w^*\coloneqq \xmin-\alpha_* \grad[\x] f (\xmin, \givenPrm)$ and then apply Theorem~\ref{thm:IFT}. It is straightforward to check the validity of Assumption~\ref{ass:CPSO}; $h$ is partly-smooth relative to $\manif\times\spVar\times\setPrm\times(0, 2/L)$. Assumption~\ref{ass:RPD}\ref{itm:RPD-i} is satisfied because $\xmin = \bwd (\w^*, \givenPrm)\in\manif$ and the smooth component, that is, $\frac{1}{2}\norm{\cdot - \w}^2$ is strongly convex. Finally from \eqref{eq:ND}, we note that $\w^* \in \xmin + \ri \alpha \partial_{\x} g (\xmin, \givenPrm)$, that is, Assumption~\ref{ass:ND} is also satisfied. To obtain the expressions for $D_{\x} \bwd$ and $D_{\u} \bwd$ at $(\w, \u, \alpha)$ sufficiently close to $(\xmin - \alpha \grad[\x] f (\xmin, \givenPrm), \givenPrm, \alpha_*)$, we use \eqref{eq:IFT}. That is, we set $\x = \bwd (\w, \u)$ and find that $\Hess[\manif] h(\x, \w, \u, \alpha) = \alpha \Hess[\manif] g(\x, \u) + \Hess[\manif] (1/2\norm{\cdot - \w}^2) (\x) = \alpha \Hess[\manif] g(\x, \u) + \projTan (\x) + \wein[\x]{\cdot}{\projNor (\x) (\x - \w)} $ from \eqref{eq:Riem:hess}, $D_{\w} \grad[\manif] h (\x, \w, \u, \alpha) = -\projTan (\x)$ and $D_{\u} \grad[\manif] h (\x, \w, \u, \alpha) = D_{\u} \grad[\manif] g(\x, \u)$.
\end{proof}

\appSubSect{Corollary}{cor:PGD:Derv}
\vspace{2ex}
\begin{proof}
    The proof follows by fixing  $\alpha_*\in (0, 2/L)$ and applying Lemma~\ref{lem:prox:diff}, the $C^1$-smoothness of the map $(\x, \u, \alpha) \mapsto (\x - \alpha \grad[\x] f (\x, \u), \u, \alpha)$ and the chain rule.
\end{proof}

\appSubSect{Corollary}{cor:APG:Derv}
\vspace{2ex}
\begin{proof}
    Given $\alpha_*\in (0, 2/L)$, we apply Corollary~\ref{cor:PGD:Derv} and choose a separable neighbourhood $\neighbourhoodPGD[*] = B_{\delta_{\x}} (\xmin) \times B_{\delta_{\u}} (\givenPrm) \times B_{\delta_{\alpha}} (\alpha_*)$. The proof then follows from the convexity of $B_{\delta_{\x}} (\xmin)$.
\end{proof}
    
\section{Proofs of Section~\ref{sec:PDM}}

\appSubSect{Theorem}{thm:PGD:IFT}
\vspace{2ex}
\begin{proof}
    Let $\map{\norm{\cdot}}{\spLin(\spVar, \spVar)}{\R}$ be the induced operator norm from the vector norm $\map{\norm{\cdot}}{\spVar}{\R}$. Since $D_{\x} \pgd (\x, \u)$ is symmetric because $\Dxbwd (\x, \u)$ is symmetric, the spectral radius may be replaced with the induced norm. Thus
    \begin{equation*}
        \begin{aligned}
            \norm{D_{\x} \pgd (\xmin, \givenPrm) \projTan (\xmin)} &= \norm{\Dxbwd (\xmin, \givenPrm)^{\dagger} (\opid - \alpha \Hess[\x] f(\xmin, \givenPrm)) \projTan (\xmin)} \\
            &= \norm{\Dxbwd (\xmin, \givenPrm)^{\dagger} (\projTan (\xmin) - \alpha \projTan (\xmin) \Hess[\x] f(\xmin, \givenPrm) \projTan (\xmin))} \\
            &\leq \norm{\Dxbwd (\xmin, \givenPrm)^{\dagger}} \norm{\projTan (\xmin) - \alpha \projTan (\xmin) \Hess[\x] f(\xmin, \givenPrm) \projTan (\xmin)} \,.
        \end{aligned}
    \end{equation*}
    The smallest eigenvalue $\mu_f$ of the symmetric linear operator $\projTan (\xmin) \Hess[\x] f (\xmin, \u) \projTan (\xmin)$ is positive from Assumption~\ref{ass:RPD}\ref{itm:RPD-ii}. By noting that $\mu_f I_{\spTan\xmin} \preceq \rstDom{\projTan (\xmin) \Hess[\x] f (\xmin, \u)}{\spTan\xmin} \preceq L I_{\spTan\xmin}$ and $0<\alpha_*<2/L$ we find that $\norm{\projTan (\xmin) - \projTan (\xmin) \alpha_* \Hess[\x] f(\xmin, \u) \projTan (\xmin)} \leq \max(\abs{1-\alpha_* \mu_f}, \abs{1-\alpha_* L})$ which is less than $1$ when $\mu_f > 0$.
    
    When Assumption~\ref{ass:RPD}\ref{itm:RPD-i} is satisfied, Theorem~\ref{thm:IFT} guarantees the existence of a neighbourhood $\neighbourhoodPRM$ of $\givenPrm$ and a $C^1$-smooth map $\argmap$ such that $\argmap (\u) = \argmin F(\cdot, \u)$ for all $\u\in \neighbourhoodPRM$. For any $\u\in \neighbourhoodPRM$, we set $\x \coloneqq \argmap (\u)$, $X\coloneqq D \argmap (\u) = -\Hess[\manif] F(\x, \u)^{\dagger} D_{\u} \grad[\manif] F(\x, \u)$, and $\w \coloneqq \x - \alpha \grad[\x] f(\x, \u)$. Thus, from Remark~\ref{rem:PGD:APG}\ref{itm:PGD:APG:FPE}, we have $\x = \pgd (\x, \u) = \bwd (\w, \u)$ and from Corollary~\ref{cor:PGD:Derv}, we have
    \begin{equation*}
        D_{\u} \pgd (\x, \u) = -\alpha\Dxbwd (\x, \u)^\dagger D_{\u} (\grad[\x] f (\x, \u) + \grad[\manif] g (\x, \u)) = -\alpha\Dxbwd (\x, \u)^\dagger D_{\u} \grad[\manif] F (\x, \u) \,,
    \end{equation*}
    and
    \begin{equation*}
        \begin{aligned}
            \left( \opid - D_{\x} \pgd (\x, \u)\projTan (\x) \right) X &= \left( \projTan (\x) - \Dxbwd (\x, \u)^\dagger \left(\opid - \alpha \Hess[\x] f(\x, \u) \right) \projTan (\x) \right) X \\
            &= \left( \projTan (\x) - \Dxbwd (\x, \u)^\dagger \left( \projTan (\x) - \alpha \projTan (\x) \Hess[\x] f(\x, \u) \projTan (\x) \right) \right) X \,.
        \end{aligned}
    \end{equation*}
    However, from the definition of Riemannian Hessian, we evaluate $\Dxbwd (\x, \u)$ by using Corollary~\ref{cor:PGD:Derv} and Lemma~\ref{lem:prox:diff}, that is,
    \begin{equation*}
        \begin{aligned}
            \Dxbwd (\x, \u) &= \projTan (\x) + \wein[\x]{\cdot}{\projNor (\x) (\x - \w)} + \alpha \Hess[\manif] g (\x, \u) \\
            &= \projTan (\x) + \alpha \wein[\x]{\cdot}{\projNor (\x)\grad[\x] f (\x, \u)} + \alpha \Hess[\manif] g (\x, \u) \\
            &= \projTan (\x) + \alpha \Hess[\manif] F (\x, \u) - \alpha \projTan (\x) \Hess[\x] f (\x, \u) \projTan (\x) \,,
        \end{aligned}
    \end{equation*}
    and arrive at the first equality, that is,
    \begin{equation*}
        \begin{aligned}
            \left( \opid - D_{\x} \pgd (\x, \u)\projTan (\x) \right) X &= \left( \projTan (\x) - \Dxbwd (\x, \u)^\dagger \left( \Dxbwd (\x, \u) - \alpha \Hess[\manif] F (\x, \u) \right) \right) X \\
            = \alpha \Dxbwd (\x, \u)^\dagger \Hess[\manif] F (\x, \u) &X = -\alpha \Dxbwd (\x, \u)^\dagger D_{\u} \grad[\manif] F(\x, \u) = D_{\u} \pgd (\x, \u) \,.
        \end{aligned}
    \end{equation*}
\end{proof}
    
\appSubSect{Corollary}{cor:PGD:IFT}
\vspace{2ex}
\begin{proof}
    For all $\u\in \neighbourhoodPRM$, we note that $\argmap (\u) = \pgd (\argmap (\u), \u)$ (see Remark~\ref{rem:PGD:APG}\ref{itm:PGD:APG:FPE}). From Corollary~\ref{cor:PGD:Derv}, Theorem~\ref{thm:PGD:IFT} and the continuity of $\rho$, we can find a neighbourhood $V\subset \neighbourhoodPRM$ of $\givenPrm$ such that for all $\u\in V$, $\pgd$ is differentiable at $(\u, \argmap (\u))$ and $\rho (D_{\x} \pgd (\argmap (\u), \u)) < 1$ which implies that $\argmap (\u) \in \manif$ from Corollary~\ref{cor:PGD:Derv}. The expression for the derivative is obtained by using chain rule on $\argmap (\u) = \pgd (\argmap (\u), \u)$ and noting that $\im D\argmap (\u) \subset \spTan{\argmap (\u)}$.
\end{proof}
    
\appSubSect{Lemma}{lem:spec_props:M}
\vspace{2ex}
\begin{proof}
    Notice that for any eigenvalue $\sigma$ of $M_{\beta}$ with corresponding eigenvector $(\v_1, \v_2)$, evaluating $M_{\beta} (\v_1, \v_2) = \sigma (\v_1, \v_2)$ implies either $\v_1 \in (\im R )^{\perp}$ and $\sigma = 0$ or $\v_1 \in \im R$ with $\v_1 = \sigma \v_2$ and $\sigma$ satisfies $\sigma^{2} - (1 + \beta) \lambda \sigma + \beta \lambda = 0$ with $\lambda$ being the eigenvalue of $R$ corresponding to eigenvector $\v_1$. We can therefore use \cite[Corollary~4.9]{LFP17} to establish that $\rho (M_{\beta}) < 1$. The expression for the inverse of $\opid - M_{\beta}$ is furnished by using \cite[Theorem~2.1]{LS02}.
\end{proof}
    
\appSubSect{Theorem}{thm:APG:IFT}
\vspace{2ex}
\begin{proof}
    Theorem~\ref{thm:IFT} asserts that under the given assumptions, a $C^1$-smooth solution map $\argmap$ of our problem exists on a neighbourhood $\neighbourhoodPRM$ of $\xmin$. For any $\u\in \neighbourhoodPRM$, we set $\x\coloneqq\argmap (\u)$ and $\z\coloneqq (\x, \x)$ and note that $\z$ is a fixed-point of $\apg (\cdot, \u)$ (see Remark~\ref{rem:PGD:APG}\ref{itm:PGD:APG:FPE}). Therefore, using Corollary~\ref{cor:APG:Derv}, we obtain
    \begin{equation*}
        D_{\z} \apg (\z, \u) \projTan (\z) = \begin{bmatrix}
        (1 + \beta) D_{\x} \pgd (\x, \u) \projTan (\x) & -\beta D_{\x} \pgd (\x, \u) \projTan (\x) \\
        \projTan (\x) & 0
        \end{bmatrix} \,,
    \end{equation*}
    and
    \begin{equation*}
        D_{\u} \apg (\z, \u) = \begin{bmatrix}
            D_{\u} \pgd (\x, \u) \\ 0
        \end{bmatrix} \,.
    \end{equation*}
    But since $\rho (D_{\x} \pgd (\x, \u) \projTan (\x)) < 1$ and $-1/(1 + 2\beta) < \underline{\lambda}$ we use Lemma~\ref{lem:spec_props:M} to conclude that $\rho (D_{\z} \apg (\z, \u) \projTan (\z)) < 1$ and the operator $\opid - D_{\z} \apg (\z, \u)\projTan (\z)$
    is invertible. By using Lemma~\ref{lem:spec_props:M} we obtain the expression for $\left(\opid - D_{\z} \apg (\z, \u)\projTan (\z)\right)^{-1} D_{\u} \apg (\z, \u)$, that is,
    \begin{equation*}
        \begin{bmatrix}
            \left( \opid - D_{\x} \pgd (\x, \u) \projTan (\x) \right)^{-1} D_{\u} \pgd (\x, \u) \\ \projTan \left( \opid - D_{\x} \pgd (\x, \u) \projTan (\x) \right)^{-1} D_{\u} \pgd (\x, \u)
        \end{bmatrix} \,.
    \end{equation*}
    Comparing this result with \eqref{eq:PGD:IFT} and noting that $\left( \opid - D_{\x} \pgd (\x, \u) \projTan (\x) \right)^{-1} \v \in \spTan \x$ whenever $\v\in \spTan{\x}$, we arrive at \eqref{eq:APG:IFT}.
\end{proof}
    
\appSubSect{Lemma}{lem:sg:seq}
\vspace{2ex}
\begin{proof}
    We set $\bnuk\coloneqq \bnu (\xk, \givenPrm)$ and assume that $\xk\neq \xmin$ because otherwise the proof is trivial. From subdifferential continuity of $g(\cdot, \givenPrm)$ relative to $\manif$, we note that $\partial_{\x} g (\xk, \givenPrm)$ converges to $\partial_{\x} g (\xmin, \givenPrm)$. Since $-\grad[\x] f (\xmin, \givenPrm) \in \partial_{\x} g (\xmin, \givenPrm)$, from the definition of convergence of $\partial_{\x} g (\xk, \givenPrm)$ \cite[Chapter~4]{RW98}, there exists a sequence $(\bmuk )_{k\in\N}$ converging to $-\grad[\x] f (\xmin, \givenPrm)$ such that for every $k\in\N, \bmuk\in\partial_{\x} g (\xk, \givenPrm)$. By the definition of projection, we note that for every $k\in\N$
    \begin{equation*}
        \begin{aligned}
            \dist (-\grad[\x] f (\xk, \givenPrm), \partial_{\x} g (\xk, \givenPrm)) &\coloneqq \inf_{\y\in\partial_{\x} g (\xk, \givenPrm)} \norm{\y + \grad[\x] f (\xk, \givenPrm)} \\
            = \norm{\bnuk + \grad[\x] f (\xk, \givenPrm)} &\leq \norm{\bmuk + \grad[\x] f (\xk, \givenPrm)} \,.
        \end{aligned}
    \end{equation*}
    But since $\grad[\x] f (\xk, \givenPrm) \to \grad[\x] f (\xmin, \givenPrm)$ and $\norm{\bmuk + \grad[\x] f (\xk, \givenPrm)} \to 0$, the convergence follows.
\end{proof}

\appSubSect{Theorem}{thm:PGD:IFT:approx}
\vspace{2ex}
\begin{proof}
    From Theorem~\ref{thm:IFT}, there exists a neighbourhood $\neighbourhoodPRM$ of $\givenPrm$ and a $C^1$-smooth solution mapping $\argmap$ on $\neighbourhoodPRM$ such that $\rho (D_{\x} \pgd (\argmap (\u), \u) \projTan (\argmap (\u)) < 1)$. For any $\u\in \neighbourhoodPRM$, we only substitute $\pgd (\x, \u)$ and $\x-\alpha\grad[\x] f(\x, \u)$ in all the expressions in \eqref{eq:PGD:maps} by $\x$ and $\x + \alpha \bnu (\x, \u)$ respectively to obtain the mappings in \eqref{eq:PGD:mod:maps}. Therefore, as $\x\to\argmap (\u)$, we have $\pgd (\x, \u) \to \argmap (\u)$ and $\bnu (\x, \u)\to-\grad[\x] f (\argmap (\u), \u)$, and from continuity, there exists $\epsilon>0$ such that $\tilde\phi_{\alpha}$ is well defined on a set $\tilde \neighbourhoodFM_{\alpha, \epsilon}\coloneqq\set{(\x, \u) \in \manif\times \neighbourhoodPRM\setsep \norm{\x-\argmap (\u)} < \epsilon}$ and $\tilde\phi_{\alpha}(\x, \u)\to\phi_{\alpha}(\argmap (\u), \u) = D\argmap (\u)$ for all $(\x,\u)\in \tilde \neighbourhoodFM_{\alpha, \epsilon}$.
    
    The expression on the right hand side of \eqref{eq:PGD:mod:IFT:approx} is also well-defined on $\tilde \neighbourhoodFM_{\alpha, \epsilon}$ for sufficiently small $\epsilon$ because $\Hess[\manif] F(\x, \u) + \mathcal{W}_{\alpha} (\x, \u)$ is bijective on $\spTan\x$, for all $(\x, \u) \in\tilde \neighbourhoodFM_{\alpha, \epsilon}$. To verify this claim, we observe that when $\x$ is sufficiently close to $\argmap(\u)$, $\Hess[\manif] F(\x, \u)$ is invertible due to Assumption~\ref{ass:RPD}\ref{itm:RPD-i} and $\norm{\bnu (\x, \u) + \grad[\x] f (\x, \u)} $ is sufficiently small due to Lemma~\ref{lem:sg:seq}, which in turn makes $\wein[\x]{\cdot}{\alpha \projNor (\x) (\bnu (\x, \u) + \grad[\x] f (\x, \u))}$ arbitrarily small thanks to the definition of the Weingarten map.
    
    Finally we let $X$ denote the right hand side of \eqref{eq:PGD:mod:IFT:approx} and write
    \begin{equation*}
        \begin{aligned}
            \left( \opid - \tilde R_{\alpha} (\x, \u)\projTan (\x) \right) X &= \left( \projTan (\x) - \tilde Q (\x, \u)^\dagger \left(\opid - \alpha \Hess[\x] f(\x, \u) \right) \projTan (\x) \right) X \\
            &= \left( \projTan (\x) - \tilde Q (\x, \u)^\dagger \left( \projTan (\x) - \alpha \projTan (\x) \Hess[\x] f(\x, \u) \projTan (\x) \right) \right) X \,.
        \end{aligned}
    \end{equation*}
    But from the definition of Riemannian Hessian
    \begin{equation*}
        \begin{aligned}
            \tilde Q (\x, \u) &= \projTan (\x) -\alpha \wein[\x]{\cdot}{\projNor (\x) \bnu (\x, \u)} + \alpha \Hess[\manif] g (\x, \u) \\
            &= \projTan (\x) + \alpha \wein[\x]{\cdot}{\projNor (\x)\grad[\x] f (\x, \u)} + \alpha \Hess[\manif] g (\x, \u) + \mathcal{W}_{\alpha} (\x, \u) \\
            &= \projTan (\x) + \alpha \Hess[\manif] F (\x, \u) - \alpha \projTan (\x) \Hess[\x] f (\x, \u) \projTan (\x) + \mathcal{W}_{\alpha} (\x, \u) \,,
        \end{aligned}
    \end{equation*}
    which gives us
    \begin{equation*}
        \begin{gathered}
            \left( \opid - \tilde R_{\alpha} (\x, \u)\projTan (\x) \right) X = \left( \projTan (\x) - \tilde Q_{\alpha} (\x, \u)^\dagger \left( \tilde Q_{\alpha} (\x, \u) - \alpha \Hess[\manif] F (\x, \u) + \mathcal{W}_{\alpha} (\x, \u) \right) \right) X \\
            = \alpha \tilde Q_{\alpha} (\x, \u)^\dagger \left(\Hess[\manif] F (\x, \u) + \mathcal{W}_{\alpha} (\x, \u)\right) X = -\alpha \tilde Q_{\alpha} (\x, \u)^\dagger D_{\u} \grad[\manif] F(\x, \u) = \tilde S_{\alpha} (\x, \u) \,.
        \end{gathered}
    \end{equation*}
\end{proof}
    
\appSubSect{Theorem}{thm:APG:FAD}
\vspace{2ex}
\begin{proof}
    By setting $\zk\coloneqq (\xk, \xkm)$, we get a more compact representation of APG (see Remark~\ref{rem:PGD:APG}\ref{itm:PGD:APG:FPI}) and obtain $\zk\to\z^*\coloneqq (\xmin, \xmin)$. The update step reads
    \begin{equation} \label{itr:APG:FAD:Comp}
        \zkpd \coloneqq D_{\z} \apg[k] (\zk, \u) \projTan (\zk) \zkd + D_{\u} \apg[k] (\zk, \u) \ud \,,
    \end{equation}
    where we added a redundant $\projTan (\zk)$ because $\zkd\in T_{\xk} \manif \times T_{\xkm}$. From Corollary~\ref{cor:APG:Derv}, there exists $\neighbourhoodAPG[*]$, such that the mapping $(\z, \u, \alpha) \mapsto \apg (\z, \u)$ is $C^1$-smooth on $\neighbourhoodAPG[*]$ and the mapping $(\z, \u, \alpha, \beta) \mapsto \apg (\z, \u)$ is continuous on $\neighbourhoodAPG[*]\times[0, 1]$. Since $(\zk, \givenPrm, \alpha_k)$ eventually lie in $\neighbourhoodAPG[*]$, we obtain $D_{\z} \apg[k] (\zk, \givenPrm) \to D_{\z} \apg[*] (\z^*, \givenPrm)$ and $D_{\u} \apg[k] (\zk, \givenPrm) \to D_{\u} \apg[*] (\z^*, \givenPrm)$. Finally, because $\rho (D_{\z} \apg[*] \projTan (\z^*)) < 1$ from Theorem~\ref{thm:APG:IFT}, the convergence $\zkd (\givenPrm) \to \left(D\argmap(\givenPrm) \ud, D\argmap(\givenPrm) \ud\right)$ follows from Theorem~\ref{thm:LGFPI:conv}.
\end{proof}

\appSubSect{Theorem}{thm:APG:FAD:Lin:Conv}
\vspace{2ex}
\begin{proof}
    The proof follows directly by applying the results from Section~\ref{sssec:FPI:AD} to AD of Algorithm~\ref{alg:APG} due to the local Lipschitz continuity of the derivative maps $D_{\x} \pgd[*] (\xmin, \givenPrm)$ and $D_{\u} \pgd[*] (\xmin, \givenPrm)$.
\end{proof}
\appSubSect{Theorem}{thm:APG:FPFAD}
\vspace{2ex}
\begin{proof}
    From Corollary~\ref{cor:APG:Derv} we have $\im D \apg (\x, \x, \u) \subset T_{(\x, \x)} (\manif \times \manif)$ which implies that $\zkh\in T_{(\x, \x)} (\manif \times \manif)$ for all $k\geq1$ which allows us to rewrite \eqref{itr:FPI:FPFAD:APG} as
    \begin{equation*}
        \zkph = D_{\z} \apg (\x, \x, \u) \projTan (\x, \x) \zkh + D_{\u} \apg (\x, \x, \u) \ud \,.
    \end{equation*}
    The proof then follows from Lemma~\ref{lem:spec_rad:neu_ser}.
\end{proof}

\section{Experimental Details} \label{sec:Exp:supp}
Here we provide the details of our setup for the two experiments shown in Section~\ref{sec:Exp}.

\subsection{Convergence Rate Demonstration} \label{ssec:Exp:CRD:supp}
Given a matrix $A\in\R^{M\times N}$ and a vector $\b\in\R^M$, the lasso problem \cite{Tib96} aims to find a sparse vector $\x^*\in\R^N$ which roughly solves $A\x = \b$, formulated as the optimization problem given in \eqref{prob:lasso} for some $\lambda>0$. The second part of the objective in \eqref{prob:lasso} controls the sparsity of the solution $\xmin$ of the problem where larger $\lambda$ makes the solution sparser. In a similar setting, given multiple vectors $\b_i \in \R^M$ for $i = 1,\ldots,L$, we are often interested in finding sparse solutions $\xmin_i \in \R^N$ which approximately satisfy $A\x_i = \b_i$ for all $i$ in such a way that all the solutions share the same sparsity pattern. This amounts to solving the problem in \eqref{prob:group:lasso} where we stack the vectors $\b_i$ as column vectors in the matrix $B\in\R^{M\times L}$. The problem of finding the group sparsity pattern is called group lasso \cite{YL06} in the literature.

Since the $\ell_1$ norm and the $\ell_{2,1}$ group norm are partly smooth \cite{VDP+17}, the two problems presented above fit the description of \eqref{prob:min} and satisfy Assumptions~\ref{ass:CPSO} and \ref{ass:RPD} for all $\u$ when $A$ is full-rank and $M\geq N$. We verify Assumption~\ref{ass:ND} experimentally by first solving the problem for given $\u$ and check if \eqref{eq:ND} is satisfied or not. As an example, we explain this procedure for \eqref{prob:lasso}. Let $\xmin$ be the solution of \eqref{prob:lasso} for some given $\u$. The vector $\v \coloneqq -\grad[\x] f(\x, \u) / \lambda = A^T (\b - A\xmin) / \lambda$ must be a subgradient of $\norm[1]{\cdot}$ at $\xmin$. Therefore for all $i$, $x^*_i \neq 0$ simply implies $v_i = 0$ and for \eqref{eq:ND} to hold at $(\xmin, \u)$, $v_i$ must lie in the open interval $(-1, 1)$ whenever $x^*_i = 0$.

For both problems, we perform Algorithm~\ref{alg:APG} and its forward and reverse mode AD and FPAD (five algorithms in total) twice --- once with $\beta_k=0$ (PGD) and once with $\beta_k>0$ (APG) --- to generate $20$ sequences in total; $10$ for each problem (see Figure~\ref{fig:conv:rates}). For simplicity, we do not construct the full Jacobian of the iterates of PGD and APG for AD and FPAD. We instead just perform the AD and FPAD algorithms for some $\ud\in\R^N$ and $\xminb\in\R^{1\times N}$ (a row vector). The Riemannian Hessian of $\norm[1]{\cdot}$ vanishes everywhere while its Riemannian Gradient is $(\sign (x_1), \ldots, \sign (x_N))$ \cite[Example~5.2.1]{Lia16} where $\map{\sign}{\R}{\set{-1, 0, 1}}$ is defined by
\[
    \sign(t) \coloneqq \begin{cases} 
    -1, & \textnormal{if} \ t < 0 \,; \\
    \;0, & \textnormal{if} \ t = 0 \,; \\
    +1, & \textnormal{if} \ t > 0 \,.
   \end{cases} 
\]
The affine manifold $\manif$ and $\spTan\xmin$ are both given by $\set{\x\in\R^N\setsep\supp{\x}\subset\supp{\xmin}}$, where $\supp{\x}$ denotes the support of a vector $\x\in\R^N$ defined as $\supp{\x} \coloneqq \set{i\in\N \setsep x_i \neq 0}$.
In other words, for all $\x$,
\begin{equation} \label{eq:lasso:manif}
    \x\in\manif=\spTan\xmin \quad\iff\quad (x^*_i = 0 \implies x_i = 0) \,.
\end{equation}
This allows us to compute the Riemannian Gradient and Hessian of $F_{\mathrm{l}}$ as
\begin{equation*}
    \begin{aligned}
        \grad[\manif] F_{\mathrm{l}} (\x, \u) &= \proj{\spTan\xmin} (A^T(A\x - \b)) + \lambda \grad[\manif] \norm[1]{\cdot} (\x) \\
        \Hess[\manif] F_{\mathrm{l}} (\x, \u) &= \proj{\spTan\xmin} \circ A^TA \circ \proj{\spTan\xmin} \,.
    \end{aligned}
\end{equation*}
Thus from \eqref{eq:lasso:manif}, for all $i, j \notin \supp{\xmin}$, $(\grad[\manif] F_{\mathrm{l}} (\x, \u))_i = 0$ and $(\Hess[\manif] F_{\mathrm{l}} (\x, \u))_{i,j} = 0$ and we can solve a reduced linear system of $\abs{\supp{\xmin}}$ equations. The derivative $D_{\u} \grad[\manif] F(\x, \u)$ can be computed by using autograd package. Similarly, for given parameter $\u$, let $X_*$ be the solution of \eqref{prob:group:lasso}. Then, $\manif = \spTan{X_*} = \set{X\in\R^{N\times L} \setsep \supp X \subset \supp X_*} $ where we define the support of a matrix $\R^{N\times L}\ni X\coloneqq [\x_1,\ldots,\x_N]^T$ by $\supp X \coloneqq \set{i\in\N\setsep \x_i \neq 0}$. Equivalently, we have
\begin{equation} \label{eq:glasso:manif}
    X\in\manif=\spTan{X_*} \quad\iff\quad (\x^*_i = 0 \implies \x_i = 0) \,.
\end{equation}
The Riemannian gradient and Hessian of $\norm[2,1]{\cdot}$ is given by
\begin{equation*}
    \begin{aligned}
        \big( \grad[\manif] \norm[2,1]{\cdot} (X) \big)_i &= \begin{cases}
            0,& \text{if } \x_i = 0\,; \\ \frac{\x_i}{\norm[2]{\x_i}},& \text{if } \x_i \neq 0\,,
        \end{cases} \\
        \big( \Hess[\manif] \norm[2,1]{\cdot} (X) \big)_{i,j} &= \begin{cases}
            \frac{1}{\norm[2]{\x_i}} (\opid - \frac{\x_i\x_i^T}{\norm[2]{\x_i}^2}),& \text{if } i=j \text{ and } \x_i \neq 0 \,; \\ 0,& otherwise \,, 
        \end{cases}
    \end{aligned}
\end{equation*}
which leads us to the following expressions:
\begin{equation*}
    \begin{aligned}
        \grad[\manif] F_{\mathrm{gl}} (X, \u) &= \proj{\spTan\xmin} (A^T(AX - B)) + \lambda \grad[\manif] \norm[2,1]{\cdot} (X) \\
        \Hess[\manif] F_{\mathrm{gl}} (X, \u) &= \proj{\spTan\xmin} \circ A^TA \circ \proj{\spTan\xmin} + \lambda \Hess[\manif] \norm[2,1]{\cdot} (X) \,.
    \end{aligned}
\end{equation*}
From the expressions above, it is easy to verify that the conditions of linear rate of convergence in Theorem~\ref{thm:APG:FAD} are satisfied by both $F_{\mathrm{l}}$ and $F_{\mathrm{gl}}$.

For \eqref{prob:lasso}, we set $M=1000$, and $N=250$ and for \eqref{prob:group:lasso}, we set $M=1000$, $N=100$, and $L=40$. We generate $A$ by sampling each $a_{i, j} $ from a uniform distribution with parameters $0$ and $1$ and $\b$ and $B$ by sampling each $b_i$ and $b_{i,j}$ from a normal distribution with mean $0$ and standard deviation $1$. We run each algorithm on \eqref{prob:lasso} for $K=5000$ iterations and on \eqref{prob:group:lasso} for $K=1000$ iterations. The computed forward mode sequences are $(\xkd)_{k\in[K]}$ and $(\xkh)_{k\in[K]}$ while the reverse mode sequences are $(\unb)_{n\in[K]}$ and $(\unt)_{n\in[K]}$. For a fair comparison, we use the same initializations for all algorithms. For performing FPAD on an optimization algorithm, we use its own output as approximate fixed-point. For example, we run PGD on \eqref{prob:lasso} for $K=5000$ iterations and use this $\xK$ for performing Forward and Reverse Mode FPAD on PGD. We use a constant step size for all iterations $\alpha_k = 1/\norm[\mathrm{op}]{A^TA}$, where $\norm[\mathrm{op}]{\cdot}$ is the operator norm and for APG, we set $\beta_k \coloneqq (k-1) / (k+q)$ with $q=5$ and for FPAD on APG, we set $\beta\coloneqq\beta_K$. For computing $\xmin$, we run \eqref{itr:APG} for a large number of iterations until it satisfies $\norm[2]{\xmin - \pgd (\xmin, \u)} < 10^{-12}$.

\subsection{Learning Regularizers for Image Denoising} \label{ssec:Exp:ID:supp}
\begin{figure}[!t]
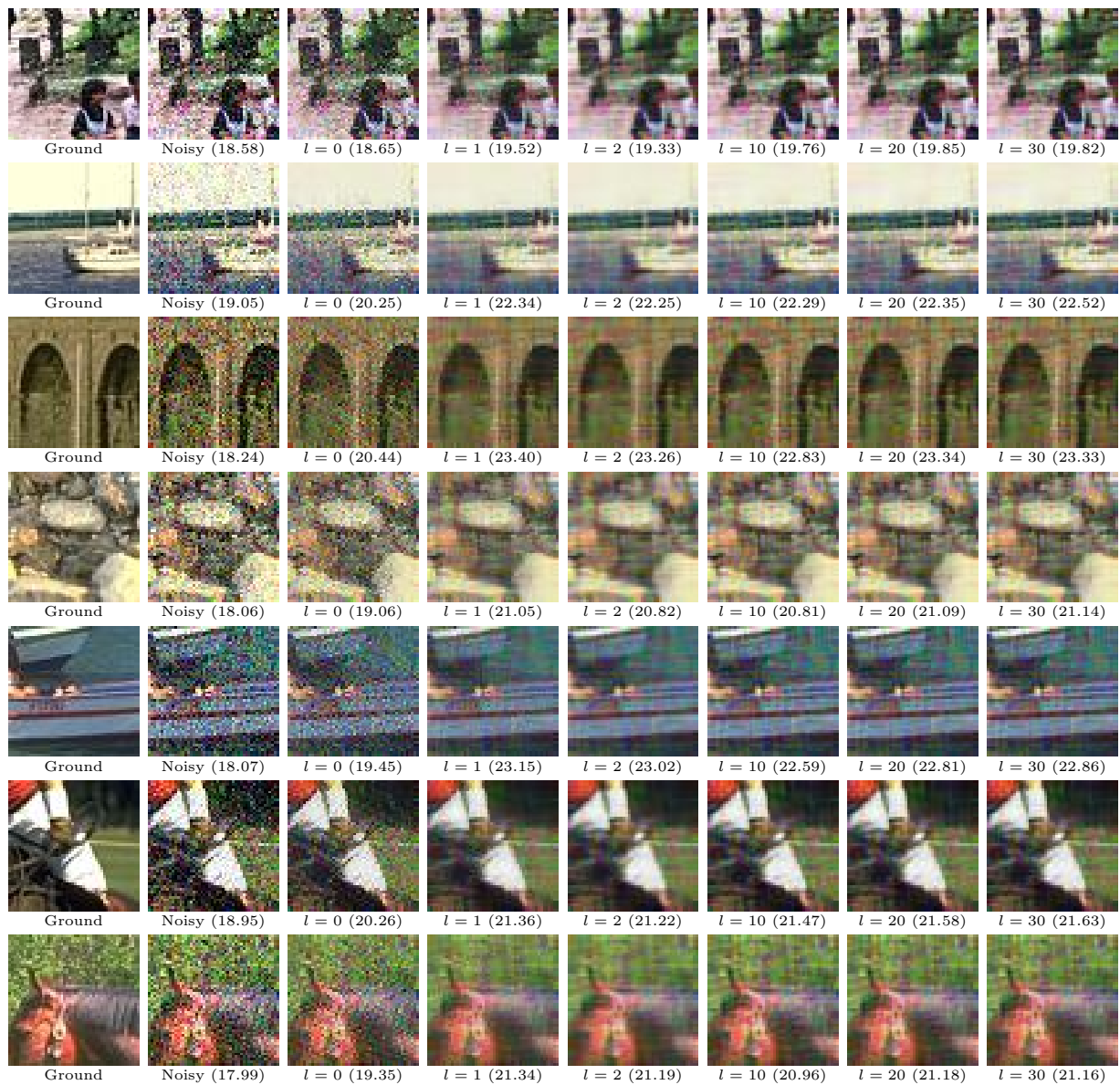

    \captionsetup[subfigure]{font=tiny,aboveskip=1pt,labelformat=empty}
    \centering
    
    \subfig{\evolwidth}{figs/denoising/train/}{Clean}{1}{}
    \subfig{\evolwidth}{figs/denoising/train/}{Noisy}{1}{\ ($18.58$)}
    \subfigevol{\evolwidth}{train}{1}{0}{$18.65$}
    \subfigevol{\evolwidth}{train}{1}{1}{$19.52$}
    \subfigevol{\evolwidth}{train}{1}{2}{$19.33$}
    \subfigevol{\evolwidth}{train}{1}{10}{$19.76$}
    \subfigevol{\evolwidth}{train}{1}{20}{$19.85$}
    \subfigevol{\evolwidth}{train}{1}{30}{$19.82$}
    \hfill%
    
    % \vspace{-1ex}
    \subfig{\evolwidth}{figs/denoising/train/}{Clean}{2}{}
    \subfig{\evolwidth}{figs/denoising/train/}{Noisy}{2}{\ ($19.05$)}
    \subfigevol{\evolwidth}{train}{2}{0}{$20.25$}
    \subfigevol{\evolwidth}{train}{2}{1}{$22.34$}
    \subfigevol{\evolwidth}{train}{2}{2}{$22.25$}
    \subfigevol{\evolwidth}{train}{2}{10}{$22.29$}
    \subfigevol{\evolwidth}{train}{2}{20}{$22.35$}
    \subfigevol{\evolwidth}{train}{2}{30}{$22.52$}
    \hfill%
    
    % \vspace{-1ex}
    \subfig{\evolwidth}{figs/denoising/train/}{Clean}{3}{}
    \subfig{\evolwidth}{figs/denoising/train/}{Noisy}{3}{\ ($18.24$)}
    \subfigevol{\evolwidth}{train}{3}{0}{$20.44$}
    \subfigevol{\evolwidth}{train}{3}{1}{$23.40$}
    \subfigevol{\evolwidth}{train}{3}{2}{$23.26$}
    \subfigevol{\evolwidth}{train}{3}{10}{$22.83$}
    \subfigevol{\evolwidth}{train}{3}{20}{$23.34$}
    \subfigevol{\evolwidth}{train}{3}{30}{$23.33$}
    \hfill%
    
    % \vspace{-1ex}
    \subfig{\evolwidth}{figs/denoising/train/}{Clean}{4}{}
    \subfig{\evolwidth}{figs/denoising/train/}{Noisy}{4}{\ ($18.06$)}
    \subfigevol{\evolwidth}{train}{4}{0}{$19.06$}
    \subfigevol{\evolwidth}{train}{4}{1}{$21.05$}
    \subfigevol{\evolwidth}{train}{4}{2}{$20.82$}
    \subfigevol{\evolwidth}{train}{4}{10}{$20.81$}
    \subfigevol{\evolwidth}{train}{4}{20}{$21.09$}
    \subfigevol{\evolwidth}{train}{4}{30}{$21.14$}
    \hfill%
    
    % \vspace{-1ex}
    \subfig{\evolwidth}{figs/denoising/train/}{Clean}{5}{}
    \subfig{\evolwidth}{figs/denoising/train/}{Noisy}{5}{\ ($18.07$)}
    \subfigevol{\evolwidth}{train}{5}{0}{$19.45$}
    \subfigevol{\evolwidth}{train}{5}{1}{$23.15$}
    \subfigevol{\evolwidth}{train}{5}{2}{$23.02$}
    \subfigevol{\evolwidth}{train}{5}{10}{$22.59$}
    \subfigevol{\evolwidth}{train}{5}{20}{$22.81$}
    \subfigevol{\evolwidth}{train}{5}{30}{$22.86$}
    \hfill%
    
    % \vspace{-1ex}
    \subfig{\evolwidth}{figs/denoising/valid/}{Clean}{2}{}
    \subfig{\evolwidth}{figs/denoising/valid/}{Noisy}{2}{\ ($17.99$)}
    \subfigevol{\evolwidth}{valid}{2}{0}{$19.35$}
    \subfigevol{\evolwidth}{valid}{2}{1}{$21.34$}
    \subfigevol{\evolwidth}{valid}{2}{2}{$21.19$}
    \subfigevol{\evolwidth}{valid}{2}{10}{$20.96$}
    \subfigevol{\evolwidth}{valid}{2}{20}{$21.18$}
    \subfigevol{\evolwidth}{valid}{2}{30}{$21.16$}
    \hfill%
    
    % \vspace{-1ex}
    \subfig{\evolwidth}{figs/denoising/valid/}{Clean}{1}{}
    \subfig{\evolwidth}{figs/denoising/valid/}{Noisy}{1}{\ ($18.95$)}
    \subfigevol{\evolwidth}{valid}{1}{0}{$20.26$}
    \subfigevol{\evolwidth}{valid}{1}{1}{$21.36$}
    \subfigevol{\evolwidth}{valid}{1}{2}{$21.22$}
    \subfigevol{\evolwidth}{valid}{1}{10}{$21.47$}
    \subfigevol{\evolwidth}{valid}{1}{20}{$21.58$}
    \subfigevol{\evolwidth}{valid}{1}{30}{$21.63$}
    \hfill%
    
    \caption{Performance of our model depicted over the epochs on train (first five rows) and test (last two rows) images. The first two columns correspond to clean $\ground[i]$ and the noisy images $\noisy[i]$, respectively, while the remaining six columns correspond to the recovered images $\recovered[i]\iter K (\Theta_l)$ for $l\in\set{0,1,2,10,20,30}$. We also show PSNR value of each image relative to the ground truth.}
    \label{fig:ID:evol:supp}
\end{figure}

Consider \eqref{prob:ID} as the image denoising problem where $N_c\in\set{1,3}$, $N_f=2$, $\imgLinOp$ computes the $x$ and $y$-derivatives of the input image and $\norm[\mathrm{r}]{\cdot}$ is the $\ell_1$ norm or the $\ell_{2,1}$ group norm. The first term in the objective, also known as the data or the fidelity term makes sure that $\recovered^*$ is not too far away from $\x$. The second term, called the regularization term, ensures the sparsity of the derivative of the image which relies on the assumption that a natural image is ``sort of'' piecewise constant and is controlled by regularization parameter $\lambda$. Chambolle and Pock~\cite{CP16} list various variants of this model, for instance, image deblurring, inpainting and zooming and list appropriate algorithms to solve the problem, which our experiment easily generalizes to.

Instead of simply using $\imgLinOp$ as the discrete derivative operator, we instead try to learn it. This problem has also been well-studied and numerous strategies have been proposed \cite[see, for example,]{RB05, CPR+13, KP13, KEK+20}. We take inspiration from \cite{CPR+13} and define $\imgLinOp$ by $\imgLinOp\recovered \coloneqq (\lrndFil{1} \recovered, \ldots, \lrndFil{N_f} \recovered)$ where each $\lrndFil{r}$ is a linear combination of $N_b$ basis filters $\basFil{s} \in \spLin (\spImg, \spImg)$, that is, $\lrndFil{r} = \sum_s \theta_{r, s} \basFil{s}$ for some $\Theta\in\R^{N_f\times N_b}$. Since $\imgLinOp$ is linear in $\Theta$, we can ignore $\lambda$ by simply absorbing it into $\Theta$. An example of $\basFil{s}$ is the $2$D DCT basis \cite{ANR74}. Our goal is to learn the weight matrix $\Theta$. In particular, given a training dataset $\trainset\coloneqq\set{(\noisy[i], \ground[i]) \setsep i=1,\ldots,M_{\textnormal{tr}}}$, we solve the bilevel problem in \eqref{prob:BL}.

\begin{figure}[ht]
    \captionsetup[subfigure]{font=tiny,aboveskip=1pt,labelformat=empty}
    \centering
    
    \subfig{0.32}{figs/denoising/test/}{Clean}{1}{}
    \subfig{0.32}{figs/denoising/test/}{Noisy}{1}{\ ($18.75$)}
    \subfig{0.32}{figs/denoising/test/}{Recovered}{1}{\ ($24.79$)}
    \hfill%
    
    \subfig{0.32}{figs/denoising/test/}{Clean}{2}{}
    \subfig{0.32}{figs/denoising/test/}{Noisy}{2}{\ ($18.92$)}
    \subfig{0.32}{figs/denoising/test/}{Recovered}{2}{\ ($25.4$)}
    \hfill%
    
    \caption{Performance of our model after training for $30$ epochs on two unseen (larger) images. The columns respectively show the clean $\ground$, noisy $\noisy$ and recovered images $\recovered\iter K (\Theta_{30})$.}
    \label{fig:ID:test:supp}
\end{figure}

Algorithm~\ref{alg:BL} lists the steps for solving \eqref{prob:BL}. The idea is to randomly sample $(\noisy[i], \ground[i])$ from $\trainset$ to compute $\grad J_i (\Theta) = D_{\Theta}\recovered[i]^* (\Theta) (\recovered[i]^* (\Theta) - \ground[i])$ and use the update $\Theta\leftarrow\Theta - \tau \grad J_i (\Theta)$. In practice, stochastic algorithms consist of epochs where each epoch involves randomly shuffling the training dataset and drawing the samples $(\noisy[i], \ground[i])$ sequentially from the shuffled set until it is exhausted \cite{GBC16}. We will use reverse mode FPAD of Algorithm~\ref{alg:APG} to compute $\grad J_{i} (\Theta)$.
\begin{Algorithm} \ \label{alg:BL}
    \begin{itemize}
        \item \key{Initialization:} $\Theta_0 \in\R^{N_f \times N_b}$.
        \item \key{Parameter:} $K\in\N$, $(\tau_l)_{l\in\N_0}$.
        \item \key{Define:} $l\coloneqq0$, $i\coloneqq0$ and $\dataset\coloneqq\trainset$.
        \item \key{Update:}
        \begin{itemize}
            \item Draw $(\noisy[i], \ground[i])$ from $\dataset$ and set
            \begin{equation*}
                \dataset\leftarrow\dataset\backslash\set{(\noisy[i], \ground[i])} \,.
            \end{equation*}
            \item Run Algorithm~\ref{alg:APG} with $\edges[i]\iter 0\coloneqq \imgLinOp\noisy[i]$ on
            \begin{equation} \label{prob:ID:FRDual}
                \min_{\edges\in\spFeat} \frac{1}{2}\norm[2]{\noisy[i] - \imgLinOp (\Theta)^*\edges}^2 \quad \st \quad \norm[\mathrm{r},*]{\edges} \leq 1 \,,
            \end{equation}
            for $K$ iterations to obtain $\edges[i]\iter K$.
            \item Compute $\recovered[i] (\Theta)\iter K = \noisy[i] - \imgLinOp^* (\Theta) \edges[i]\iter K(\Theta)$ and $J_{i}\iter K (\Theta) = \frac{1}{2} \norm[2]{\recovered[i]\iter K (\Theta) - \ground[i]}^2$.
            \item Compute $\grad J_{i}\iter K (\Theta)$ by some standard autograd package. Run reverse mode FPAD of Algorithm~\ref{alg:APG} for $K$ iterations for Back-propagation through \eqref{prob:ID:FRDual}.
            \item Update $\Theta_{i+1} \coloneqq \Theta_{i} - \tau_l \grad J_i\iter K (\Theta_{l})$.
            \item When
            $\dataset=\emptyset$; increment $l \leftarrow l + 1$, and replenish $\dataset\coloneqq\trainset$.
        \end{itemize}
    \end{itemize}
\end{Algorithm}

For our experimental setup, we use a training dataset of $5$ images and a test dataset of $2$ images, each being a colored image patch normalized to $[0, 1]$ with size $50\times50\times3$, taken from \cite{MFT+01}. To each clean image, we add a Gaussian noise with standard deviation $40/255$ to generate the noisy images. We run Algorithm~\ref{alg:BL} for $30$ epochs and define the learning rate sequence $\tau_l \coloneqq 10^{-4} / (\floor{l/4}+1)$, that is, we start with a learning rate of $10^{-4}$ and decay it linearly after every $4$ epochs. We use the built-in SGD optimizer of PyTorch for this purpose and set the momentum to $0.75$. We learn $N_f = 24$ filters by using $5\times5$ DCT bases without the constant basis vector, giving us $N_b = 24$ basis filters and a weight matrix $\Theta$ of size $24\times24$. For initialization, we draw each entry of $\Theta$ from a uniform distribution with parameters $0$ and $0.01$. We run Algorithm~\ref{alg:APG} and its reverse mode FPAD for $K=500$ iterations each. In the two algorithms, we set $\beta_k\coloneqq (k-q)/(k+1)$ and $\beta\coloneqq (K-q)/(K+1)$ with $q=5$ and use a constant step size $\alpha_k\coloneqq 1/L^2$ where $L$ is an upper bound on the operator norm of $\imgLinOp$. This upper bound can be computed by noting that the operator norm of a convolution operator with kernel matrix $C$ is bounded from above by $\norm[1]{C}$. In our experiments, we use PyTorch to implement the forward pass which involves running Algorithm~\ref{alg:APG} to obtain $\edges[i]\iter K$ for any given $\Theta$ followed by computing the recovered image $\recovered[i]\iter K$ and the batch loss $J_i\iter K$. Algorithm~\ref{alg:APG} is implemented as a function with custom gradient where we use reverse mode FPAD to implement the backward pass.

\begin{figure}[ht]
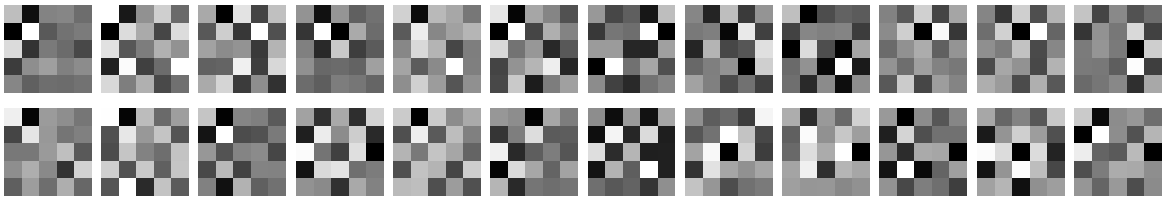

    \captionsetup[subfigure]{font=tiny,aboveskip=1pt,labelformat=empty}
    
    \centering
    \foreach \filt in {1,...,12}{%
        \subfigfilt{\filtwidth}{\filt}
    }%
    \hfill%
    
    \vspace{1ex}
    \centering
    \foreach \filt in {13,...,24}{%
        \subfigfilt{\filtwidth}{\filt}
    }%
    \hfill%
    
    \caption{$N_f=24$ learned filters after $30$ epochs.}
    \label{fig:ID:filters}
\end{figure}

In Figure~\ref{fig:ID:evol:supp}, we show the evolution of our model over epochs by showing its performance on train images (first five rows) and test (last two rows) images. Before training, the model outputs an image (third column) which is closer to the noisy image (second column). This happens because, the entries of $\Theta$ are very small and the regularization has close to no impact on the output. However, after only the first epoch (fourth column), the model starts to improve and is already close to the clean image (first column). We can also see how only with such a small dataset, the model learns to perform really well on unseen images like two images of same size as that of the training images (bottom two rows in Figure~\ref{fig:ID:evol:supp}) and two larger images of size $481\times381\times3$ (the two rows in
Figure~\ref{fig:ID:test:supp}).

Figure~\ref{fig:ID:filters} displays the $N_f=24$ learned filters, that is, the $24$ kernels each corresponding to some $\lrndFil{i}$ for $i=1,\ldots,N_f$.

\section*{Acknowledgments}
Sheheryar Mehmood and Peter Ochs are supported by the German Research Foundation
(DFG Grant OC 150/4-1).

% \section{Proofs}

% ************************
% >>>>> bibliography <<<<<
% ************************
{\small
\bibliographystyle{ieee}
\bibliography{main}
}

\end{document}

%% file: commands.tex
\newcommand{\R}{\mathbb R}             % real numbers
\newcommand{\N}{\mathbb N}             % natural numbers
\newcommand{\eR}{\overline{\R}}        % extended real numbers

\newcommand{\st}{\textrm{s.t.} \ }          % subject to

\DeclareMathOperator*{\Argmin}{Argmin}

\DeclareMathOperator*{\argmin}{argmin}

\newcommand{\setsep}{\,:\,}

\newcommand{\iter}[1]{^{(#1)}}

\let\phitmp\phi
\let\phi\varphi
\let\varphi\phitmp
\let\epstmp\eps
\let\eps\varepsilon
\let\varepsilon\epstmp

\newcommand{\set}[1]{\lbrace #1 \rbrace}                    % generic set
\newcommand{\ri}{\mathrm{ri}\,}                             % relative interior of a set
\newcommand{\rbd}{\mathrm{rbd}\,}                           % relative boundary of a set
\def\evolwidth{0.115}
\newcommand{\im}{\operatorname{im}}                         % image of linear map
\renewcommand{\ker}{\operatorname{ker}}                     % kernel of linear map
\newcommand{\seq}[2][]{(#2)_{#1}}                            % sequence
\newcommand{\abs}[1]{\vert #1 \vert}                         % absolute value
\newcommand{\norm}[2][]{\Vert #2 \Vert_{#1}}                 % norm
\newcommand{\scal}[2]{\left\langle #1,#2 \right\rangle}      % scalar product
\newcommand{\map}[3]{#1\colon #2 \to #3}                     % standard mapping
\newcommand{\sign}{\mathrm{sign}}                            % signum function
\newcommand{\dist}{\mathrm{dist}}                            % distance function
\newcommand{\opid}{I}                                        % identity operator

\newcommand{\grad}[1][]{\nabla_{#1}}                              % gradient
\newcommand{\Hess}[1][]{\nabla^2_{#1}}                     % Hessian
\newcommand{\dydx}[2]{\frac{d#1}{d#2}}                  % Leibniz Notation for Derivative

\newcommand{\dom}{\mathrm{dom}\,}                       % effective domain
\renewcommand{\o}{\mathrm{o}}                           % little-o notation
\renewcommand{\O}{\mathcal O}                          % big-o notation
\newcommand{\smap}[3]{#1\colon #2 \rightrightarrows #3} % set valued mapping
\newcommand{\proj}[1]{\mathrm{proj}_{#1}}               % projection
\newcommand{\floor}[1]{\lfloor #1 \rfloor}              % floor function

\newcommand{\spLin}{\mathcal L}                       % space of linear operators
\newcommand{\spSob}[2][]{\ensuremath{H^{#2\ifthenelse{\equal{#1}{\empty}}{}{,#1}}}}
\newcommand{\setPrm}{\Omega}                            % Open set in parameter space
      
\newcommand{\spVar}{\mathcal X}                         % variable space
\newcommand{\spPrm}{\mathcal U}                         % parameter space
\newcommand{\spImg}{\mathcal I}                         % Image space
\newcommand{\spFeat}{\mathcal J}                        % Feature space
\newcommand{\manif}{\mathcal M}                         % Manifold in variable space
\newcommand{\spPar}[1]{\mathrm{par}\,#1}                % Parallel space to set.
\newcommand{\spTan}[1]{T_{#1} \manif}                   % Tangent space to manifold.
\newcommand{\spNor}[1]{N_{#1} \manif}                   % Normal space to manifold.

\newcommand{\x}{\bm x}                                  % optimization variable
\newcommand{\xmin}{\x^*}                                % optimization variable
\newcommand{\xz}{\x ^{(0)}}                             % k-th iterate
\newcommand{\xzm}{\x ^{(-1)}}                           % (k-1)-th iterate
\newcommand{\xk}{\x ^{(k)}}                             % k-th iterate
\newcommand{\xkp}{\x ^{(k+1)}}                          % (k+1)-th iterate
\newcommand{\xkm}{\x ^{(k-1)}}                          % (k-1)-th iterate
\newcommand{\xK}{\x ^{(K)}}                             % k-th iterate
\newcommand{\xd}{\dot\x}                                % k-th iterate
\newcommand{\xmind}{\dot\x^*}                           % optimization variable
\newcommand{\xzd}{\dot\x ^{(0)}}                        % k-th iterate
\newcommand{\xkd}{\dot\x ^{(k)}}                        % k-th iterate
\newcommand{\xkpd}{\dot\x ^{(k+1)}}                     % (k+1)-th iterate
\newcommand{\xKd}{\dot\x ^{(K)}}                        % k-th iterate
\newcommand{\xzh}{\hat\x ^{(0)}}                        % k-th iterate
\newcommand{\xkh}{\hat\x ^{(k)}}                        % k-th iterate
\newcommand{\xkph}{\hat\x ^{(k+1)}}                     % (k+1)-th iterate
\newcommand{\xminb}{\bar\x^*}                           % optimization variable
\newcommand{\xkb}{\bar\x ^{(k)}_K}                      % k-th iterate
\newcommand{\xkpb}{\bar\x ^{(k+1)}_K}                   % k-th iterate
\newcommand{\xKb}{\bar\x ^{(K)}_K}                      % k-th iterate
\newcommand{\xkt}{\tilde\x ^{(k)}}                      % k-th iterate
\newcommand{\xkpt}{\tilde\x ^{(k+1)}}                   % k-th iterate
\newcommand{\xKt}{\tilde\x ^{(K)}}                      % k-th iterate
\newcommand{\z}{\bm z}                                  % Concatenated variable
\newcommand{\zk}{\z ^{(k)}}                             % k-th extrapolation variable
\newcommand{\zkp}{\z ^{(k+1)}}                          % (k+1)-th iterate
\newcommand{\zK}{\z ^{(K)}}                             % k-th iterate
\newcommand{\zd}{\dot\z}                                % Concatenated variable
\newcommand{\zkd}{\dot\z ^{(k)}}                        % k-th iterate
\newcommand{\zkpd}{\dot\z ^{(k+1)}}                     % (k+1)-th iterate
\newcommand{\zkh}{\hat\z ^{(k)}}                        % k-th iterate
\newcommand{\zkph}{\hat\z ^{(k+1)}}                     % (k+1)-th iterate
\newcommand{\sslow}{\underline\alpha}
\newcommand{\ssup}{\overline\alpha}

\renewcommand{\u}{\bm u}                                % parameter
\newcommand{\ud}{\dot\u}                                % dotted parameter
\newcommand{\ub}{\bar\u}                                % barred parameter
\newcommand{\uzb}{\bar\u^{(0)} _K}                      % 0-th barred parameter
\newcommand{\unb}{\bar\u^{(n)} _K}                      % n-th barred parameter
\newcommand{\unpb}{\bar\u^{(n+1)} _K}                   % (n+1)-th barred parameter
\newcommand{\uKb}{\bar\u^{(K)} _K}                      % K-th barred parameter
\newcommand{\ubK}{\bar\u_K}                             % K-th barred parameter
\newcommand{\uzt}{\tilde\u^{(0)}}                       % 0-th tilde'ed parameter
\newcommand{\unt}{\tilde\u^{(n)}}                       % n-th tilde'ed parameter
\newcommand{\unpt}{\tilde\u^{(n+1)}}                    % (n+1)-th tilde'ed parameter
\newcommand{\givenPrm}{\u^*}

\newcommand{\y}{\bm y}                                  % Extrapolation variable
\newcommand{\yk}{\y ^{(k)}}                             % k-th extrapolation variable
\newcommand{\yd}{\dot\y}                                % k-th extrapolation variable
\newcommand{\w}{\bm w}                                  % Input to prox
\newcommand{\wk}{\w ^{(k)}}                             % k-th extrapolation variable
\newcommand{\p}{\bm p}                                % Input to prox

\renewcommand{\b}{\bm b}                                % Input to prox
\newcommand{\bk}{\b ^{(k)}}                             % k-th extrapolation variable

\renewcommand{\v}{\bm v}                                % optimization variable
\newcommand{\bnu}{\bm\nu}
\newcommand{\bnuk}{\bnu^{(k)}}
\newcommand{\bmu}{\bm\mu}
\newcommand{\bmuk}{\bmu^{(k)}}

\newcommand{\argmap}{\psi}
\newcommand{\Argmap}{\Psi}
\newcommand{\supp}[1]{\mathrm{supp} (#1)}
\newcommand{\wein}[3][]{\mathfrak W_{#1} (#2, #3)}
\newcommand{\bwd}[1][]{\mathrm{P}_{\alpha_{#1}g}}                  % (k+1)-th iterate
\newcommand{\pgd}[1][]{\mathcal{A}_{\alpha_{#1}}}
\newcommand{\apg}[1][]{\mathcal{A}_{\alpha_{#1}, \beta_{#1}}}

\newcommand{\Dxbwd}[1][]{Q_{\alpha_{#1}}}
\newcommand{\Dxpgd}[1][]{R_{\alpha_{#1}, \beta_{#1}}}
\newcommand{\Dupgd}[1][]{S_{\alpha_{#1}, \beta_{#1}}}
\newcommand{\pgdPhi}{\phi_{\alpha}}
\newcommand{\apgPhi}{\phi_{\alpha, \beta}}

\newcommand{\DRSecs}[2]{Section~#1[#2]}
\newcommand{\projTan}{\Pi}
\newcommand{\projNor}{\Pi^\perp}
\newcommand{\rstDom}[2]{\left.{#1}\right|_{#2}}                     % function with restricted domain.

\newcommand{\fixMap}[1][]{\mathcal A_{#1}}
\newcommand{\neighbourhood}{\mathcal N}
\newcommand{\neighbourhoodPRM}{U}
\newcommand{\neighbourhoodFM}{V}
\newcommand{\neighbourhoodPRX}[1][]{W_{\alpha_{#1}}}
\newcommand{\neighbourhoodPGD}[1][]{V_{\alpha_{#1}}}
\newcommand{\neighbourhoodAPG}[1][]{V_{\alpha_{#1}}}

\newcommand{\noisy}[1][]{\bm x_{#1}}
\newcommand{\ground}[1][]{\bm t_{#1}}
\newcommand{\recovered}[1][]{\bm y_{#1}}
\newcommand{\edges}[1][]{\bm p_{#1}}
\newcommand{\imgLinOp}{\mathcal D}
\newcommand{\lrndFil}[1]{A_{#1}}
\newcommand{\basFil}[1]{B_{#1}}
\newcommand{\dataset}{\mathcal T}
\newcommand{\trainset}{\dataset_{\textnormal{tr}}}

\newcommand{\appSecLabel}[1]{ssec:#1:proof}
\newcommand{\findProofIn}[1]{\begin{proof}The proof is in Section~\ref{\appSecLabel{#1}} in the appendix.\end{proof}}
\newcommand{\appSubSect}[2]{\subsection{Proof of {#1}~\ref{#2}} \label{\appSecLabel{#2}}}
\newcommand{\refSupp}[1]{Section~\ref{#1} in the appendix}

\newcommand{\subfigevol}[5]{%
    \begin{subfigure}{#1\textwidth}
        \centering
        \includegraphics[width=\linewidth]{figs/denoising/#2/recovered-#3_#4.jpg}
        \caption{$l=#4$\ (#5)}
    \end{subfigure}%
}

\newcommand{\subfig}[5]{%
    \begin{subfigure}{#1\textwidth}
        \centering
        \includegraphics[width=\linewidth]{#2#3-#4.jpg}
        \caption{#3#5}
    \end{subfigure}%
}

\newcommand{\subfigfilt}[2]{%
    \begin{subfigure}{#1\textwidth}
        \centering
        \includegraphics[width=\linewidth]{figs/denoising/filters/filter-#2.jpg}
    \end{subfigure}%
}

\def\evolwidth{0.115}
\def\filtwidth{0.07}

\newcommand{\crefItms}[2]{%
  \hyperref[#2]{\namecref{#1}~\labelcref*{#1}\ref*{#2}}%
}